\documentclass{amsart}
\usepackage{amssymb}
\usepackage{amscd}
\usepackage{verbatim}
\usepackage{epsfig}

\newcommand\calH{\mathcal{H}}

\newcommand{\rad}{\mathrm{rad}}

\newcommand\bdef{{x}}
\newcommand\bdefh{{\hat x}}
\newcommand\Ad{\operatorname{Ad}}
\newcommand\Hso{H_{\mathrm{ss},o}}
\newcommand\Diffso{\Diff_{\mathrm{ss},o}}
\begin{document}

\newcommand\Mand{\ \text{and}\ }
\newcommand\Mwith{\ \text{with}\ }
\newcommand\Mfor{\ \text{for}\ }
\newcommand\Mst{\ \text{such that}\ }
\newcommand\Mor{\ \text{or}\ }
\newcommand\Mif{\ \text{if}\ }
\newcommand\Miff{\ \text{iff}\ }
\newcommand\Mthen{\ \text{then}\ }
\newcommand\nin{\notin}
\newcommand\identity{\operatorname{id}}
\newcommand\Id{\operatorname{Id}}
\newcommand\Real{\mathbb{R}}
\newcommand\RR{\mathbb{R}}
\newcommand\CC{\mathbb{C}}
\newcommand\AAA{\mathbb{A}}
\newcommand\olM{\overline{M}}
\newcommand\olX{\overline{X}}
\newcommand\pos{\Real^+}
\newcommand\Rnp{\Real\setminus\{0\}}
\newcommand\nzero{\setminus\{0\}}
\newcommand\Cx{\mathbb{C}}
\newcommand\Cxp{\Cx^+}
\newcommand\Cxm{\Cx^-}
\newcommand\Nat{\mathbb{N}}
\newcommand\halfNat{{\frac{1}{2}}\mathbb{N}}
\newcommand\intgr{\mathbb{Z}}
\newcommand\im{\operatorname{Im}}
\newcommand\re{\operatorname{Re}}
\newcommand\loc{{\operatorname{loc}}}
\newcommand\sign{\operatorname{sign}}
\newcommand\codim{\operatorname{codim}}
\newcommand\End{\operatorname{End}}
\newcommand\Ker{\operatorname{Ker}}
\newcommand\Hom{\operatorname{Hom}}
\newcommand\tr{\operatorname{tr}}
\newcommand\Tr{\operatorname{Tr}}
\newcommand\ideal{{\mathcal I}}
\newcommand\Span{\operatorname{span}}
\newcommand\image{\operatorname{image}}
\newcommand\Range{\operatorname{Ran}}
\newcommand\Graph{\operatorname{graph}}
\newcommand\slim{\operatornamewithlimits{s-lim}}
\newcommand\frakg{\mathfrak g}
\newcommand\frakk{\mathfrak k}
\newcommand\frakp{\mathfrak p}
\newcommand\gll{\mathfrak{gl}}
\newcommand\sll{\mathfrak{sl}}
\newcommand\sol{\mathfrak{so}}
\newcommand\GL{\operatorname{G\ell}}
\newcommand\SL{\operatorname{SL}}
\newcommand\SO{\operatorname{SO}}
\newcommand\On{\operatorname{O}}
\newcommand\pa{\partial}
\newcommand\del{\partial}
\newcommand\Rn{\Real^n}
\newcommand\Rm{\Real^m}
\newcommand\RN{\Real^N}
\newcommand\RtN{\Real^{2N}}
\newcommand\RM{\Real^M}
\newcommand\HH{\mathbb{H}}
\newcommand\sphere{\mathbb{S}}
\newcommand\Sn{\sphere^{n-1}}
\newcommand\Sm{\sphere^{m-1}}
\newcommand\Snp{\sphere^n_+}
\newcommand\Smp{\sphere^m_+}
\newcommand\SN{\sphere^{N-1}}
\newcommand\SNp{\sphere^N_+}
\newcommand\circlep{\sphere^1_+}
\newcommand\Phom{P_{h}}
\newcommand\Shom{S_{h}}
\newcommand\distance{\operatorname{dist}}
\newcommand\cl{\operatorname{cl}}
\newcommand\interior{\operatorname{int}}
\newcommand\ind{\mathrm{ind}}
\newcommand\Fa{\operatorname{Fa}}
\newcommand\ff{\operatorname{ff}}
\newcommand\mf{\operatorname{mf}}
\newcommand\cf{\operatorname{cf}}
\newcommand\scf{\operatorname{sf}}
\newcommand\lf{\operatorname{lf}}
\newcommand\rf{\operatorname{rf}}
\newcommand\indfam{{\mathcal K}}
\newcommand\ev{{\lambda}}
\newcommand\al{\alpha}
\newcommand\ad{\operatorname{ad}}
\newcommand\la{\lambda}
\newcommand\fraka{{\mathfrak a}}
\newcommand\tfraka{\widetilde{\fraka}}
\newcommand\ofraka{\overline{\fraka}}
\newcommand\frakb{{\mathfrak b}}
\newcommand\frakw{{\mathfrak w}}
\newcommand\calA{{\mathcal A}}
\newcommand\calB{{\mathcal B}}
\newcommand\calD{{\mathcal D}}
\newcommand\calR{{\mathcal R}}
\newcommand\calS{{\mathcal S}}
\newcommand\calO{{\mathcal O}}
\newcommand\calJ{{\mathcal J}}
\newcommand\calM{{\mathcal M}}
\newcommand\calN{{\mathcal N}}
\newcommand\calX{{\mathcal X}}
\newcommand\calY{{\mathcal Y}}
\newcommand\calF{{\mathcal F}}
\newcommand\calG{{\mathcal G}}
\newcommand\calT{{\mathcal T}}
\newcommand\calC{{\mathcal C}}
\newcommand\calP{{\mathcal P}}
\newcommand\calU{{\mathcal U}}
\newcommand\calV{{\mathcal V}}
\newcommand\calCt{{\tilde {\mathcal C}}}
\newcommand\calW{{\mathcal W}}
\newcommand\calK{{\mathcal K}}
\newcommand\calCL{{\mathcal C}_{\text L}}
\newcommand\calCR{{\mathcal C}_{\text R}}
\newcommand\Cinf{{\mathcal C}^{\infty}}
\newcommand\dist{{\mathcal C}^{-\infty}}
\newcommand\dCinf{\dot{\mathcal C}^\infty}
\newcommand\ddist{\dot\dist}
\newcommand\Cj{{\mathcal C}^j}
\newcommand\Linf{L^{\infty}}
\newcommand\phg{{\text{phg}}}
\newcommand\bcon{{\mathcal A}}
\newcommand\bconc{{\mathcal A}_{\text{phg}}}
\newcommand\Sch{{\mathcal S}}
\newcommand\temp{\Sch^{\prime}}
\newcommand\Diff{\operatorname{Diff}}
\newcommand\Diffb{\operatorname{Diff}_{\text{b}}}
\newcommand\Diffc{\operatorname{Diff}_{\text{c}}}
\newcommand\Diffsc{\operatorname{Diff}_{\text{sc}}}
\newcommand\DiffI{\operatorname{Diff}_{\text{I}}}
\newcommand\DiffIq{\operatorname{Diff}_{\text{I},q}}
\newcommand\sing{\mathrm{sing}}
\newcommand\reg{\mathrm{reg}}
\newcommand\ot{\overline{\tau}}
\newcommand\osig{\overline{\sigma}}
\newcommand\supp{\operatorname{supp}}
\newcommand\ssupp{\operatorname{sing\ supp}}
\newcommand\csupp{\operatorname{cone\ supp}}
\newcommand\esupp{\operatorname{ess\ supp}}
\newcommand\Fr{{\mathcal F}}
\newcommand\Frinv{\Fr^{-1}}
\newcommand\bop{{\mathcal B}}
\newcommand\spec{\operatorname{spec}}
\newcommand\pspec{\spec_{pp}}
\newcommand\cspec{\spec_{c}}
\newcommand\FIO{{\mathcal I}}
\newcommand\SP{\operatorname{RC}}
\newcommand\RC{\operatorname{RC}}
\newcommand\Symc{S_c}
\newcommand\Symca{S_c^{\alpha}}
\newcommand\Symczero{S_c^{0,...,0}}
\newcommand\sci{{}^{\text{sc}}}
\newcommand\sct{\sci T^*}
\newcommand\scdt{\sci \dot T^*}
\newcommand\dS{\dot S^*}
\newcommand\dT{\dot T^*}
\newcommand\dSreg{\dot\Sigma_{\text reg}}
\newcommand\scct{\sci\bar{T}^*}
\newcommand\Csc{C_{\text{sc}}}
\newcommand\SNpscd{(\SNp)^2_{\text{sc}}}
\newcommand\scdiag{\Delta_{\text{sc}}}
\newcommand\projscl{\pi^L_{\text{sc}}}
\newcommand\projscr{\pi^R_{\text{sc}}}
\newcommand\scHL{\sci H^{2,0}_{|\zeta|^2-\lambda^2}}
\newcommand\scHrg{\sci H^{2,0}_{\sqrt{g}}}
\newcommand\Hsc{H_{\text{sc}}}
\newcommand\WF{\operatorname{WF}}
\newcommand\WFp{\operatorname{WF^{\prime}}}
\newcommand\WFsc{\operatorname{WF}_{\text{sc}}}
\newcommand\WFscp{\operatorname{WF_{sc}^{\prime}}}
\newcommand\WFC{\operatorname{WF}_C}
\newcommand\WFCi{\operatorname{WF}_{C_i}}
\newcommand\elliptic{\operatorname{ell}}
\newcommand\Psop{\operatorname{\Psi}}
\newcommand\Psiscrs{\operatorname{\Psi_{sc}^{-2,\infty}}}
\newcommand\Psiscr{\operatorname{\Psi_{sc}^{-2,0}}}
\newcommand\Psiscrm{\operatorname{\Psi_{sc}^{0,2}}}
\newcommand\PsiscHam{\operatorname{\Psi_{sc}^{2,0}}}
\newcommand\Psisci{\operatorname{\Psi_{sc}^{*,*}}}
\newcommand\Psiscid{\operatorname{\Psi_{sc}^{0,0}}}
\newcommand\Psiscis{\operatorname{\Psi_{sc}^{0,\infty}}}
\newcommand\Psiscsi{\operatorname{\Psi_{sc}^{-\infty,0}}}
\newcommand\Psiscs{\operatorname{\Psi_{sc}^{-\infty,\infty}}}
\newcommand\Psiscalg{\operatorname{\Psi_{sc}^{\infty,-\infty}}}
\newcommand\nullHam{{\mathcal N}}
\newcommand\charD{\Sigma_{\Delta-\lambda^2}}
\newcommand\charLap{\Sigma_{\Delta-\lambda}}
\newcommand\Snl{\Sn_{\lambda}}
\newcommand\SNl{\SN_{\lambda}}
\newcommand\gammat{\tilde\gamma}
\newcommand\gammasc{\gamma}
\newcommand\Tau{\mathcal{T}}
\newcommand\taut{\tilde\tau}
\newcommand\taub{\bar\tau}
\newcommand\Nout{N^+_{\lambda}}
\newcommand\Nin{N^-_{\lambda}}
\newcommand\Nio{N^{\pm}_{\lambda}}
\newcommand\El{E_{\lambda}}
\newcommand\Elt{\tilde E_{\lambda}}
\newcommand\Eil{E^i_{\lambda}}
\newcommand\Ejl{E^j_{\lambda}}
\newcommand\Eajl{E^{\alpha_j}_{\lambda}}
\newcommand\Eilt{\tilde E^i_{\lambda}}
\newcommand\Np{N^+}
\newcommand\Nm{N^-}
\newcommand\Npm{N^{\pm}}
\newcommand\Fin{F^-(\lambda)}
\newcommand\Fini{F^-_i(\lambda)}
\newcommand\Fout{F^+(\lambda)}
\newcommand\Fouti{F^+_i(\lambda)}
\newcommand\Foutj{F^+_j(\lambda)}
\newcommand\Rout{R^+_{\lambda}}
\newcommand\Routl{R^+_{\lambda^2}}
\newcommand\Routsgnl{R^{\sign\lambda}_{\lambda^2}}
\newcommand\Rin{R^-_{\lambda}}
\newcommand\Rinl{R^-_{\lambda^2}}
\newcommand\Rinsgnl{R^{-\sign\lambda}_{\lambda^2}}
\newcommand\Rio{R^{\pm}_{\lambda}}
\newcommand\Riol{R^{\pm}_{\lambda^2}}
\newcommand\Roi{R^{\mp}_{\lambda}}
\newcommand\Roil{R^{\mp}_{\lambda^2}}
\newcommand\Riob{R^{\pm}}
\newcommand\Roib{R^{\mp}}
\newcommand\Tio{T^{\pm}}
\newcommand\Tiob{T^{\pm}_{\ff}}
\newcommand\Toi{T^{\mp}}
\newcommand\Toib{T^{\mp}_{\ff}}
\newcommand\TIiob{T_I^{\pm}}
\newcommand\Rinb{R^-}
\newcommand\Rinbsgnl{R^{-\sign\lambda}}
\newcommand\Tin{T^-}
\newcommand\Tinb{T^-_{\ff}}
\newcommand\TIinb{T^-_I}
\newcommand\Routb{R^+}
\newcommand\Routbsgnl{R^{\sign\lambda}}
\newcommand\Tout{T^+}
\newcommand\Toutb{T^+_{\ff}}
\newcommand\TIoutb{T^+_I}
\newcommand\Rlkf{(|\xib|^2-(\lambda-i0)^2)^{-1}}
\newcommand\Rlk{\rho_0(\lambda)}
\newcommand\Rmlk{\rho_0(-\lambda)}
\newcommand\Rpmlk{\rho_0(\pm\lambda)}
\newcommand\Rlka{\rho_1(\lambda)}
\newcommand\Rlkb{\rho_2(\lambda)}
\newcommand\Rilk{\rho_i(\lambda)}
\newcommand\reduced{\natural}
\newcommand\Rlf{R_0(\lambda)}
\newcommand\Rla{R_1(\lambda)}
\newcommand\Rlb{R_2(\lambda)}
\newcommand\Ril{R_i(\lambda)}
\newcommand\Rlj{R_j(\lambda)}
\newcommand\Rlft{R_0(\lambda)}
\newcommand\Rflambda{R_0^{\reduced}(\sigma)}
\newcommand\RV{R^{\reduced}_V}
\newcommand\Rfsigma{R_0^{\reduced}(\sigma)}
\newcommand\Rfsigmah{R_0^{\reduced}(\sigma^{1/2})}
\newcommand\Rfzero{R_0^{\reduced}(0)}
\newcommand\RlV{R^{\reduced}_V(\sigma)}
\newcommand\RlVi{R^{\reduced}_{V_i}(\sigma)}
\newcommand\RlVt{R_V(\lambda)}
\newcommand\RlVtL{{R}_V^L(\lambda)}
\newcommand\RlVtR{{R}_V^R(\lambda)}
\newcommand\RlVit{{R}_{V_i}(\lambda)}
\newcommand\RlVta{{R}_V^{(1)}(\lambda)}
\newcommand\RlVtk{{R}_V^{(k)}(\lambda)}
\newcommand\RlVatV{{R}_{V_{\alpha}}(\lambda)V_{\alpha}}
\newcommand\RlVatVa{{R}_{V_{\alpha_1}}(\lambda)V_{\alpha_1}}
\newcommand\RlVatVb{{R}_{V_{\alpha_2}}(\lambda)V_{\alpha_2}}
\newcommand\RlVatVk{{R}_{V_{\alpha_k}}(\lambda)V_{\alpha_k}}
\newcommand\RlVatVkk{{R}_{V_{\alpha_{k+1}}}(\lambda)V_{\alpha_{k+1}}}
\newcommand\RlVaptV{{R}_{V_{\alpha'}}(\lambda)V_{\alpha'}}
\newcommand\RlVapptV{{R}_{V_{\alpha''}}(\lambda)V_{\alpha''}}
\newcommand\RlVajtV{{R}_{V_{\alpha_j}}(\lambda)V_{\alpha_j}}
\newcommand\RlVaktV{{R}_{V_{\alpha_k}}(\lambda)V_{\alpha_k}}
\newcommand\RlVakktV{{R}_{V_{\alpha_{k+1}}}(\lambda)V_{\alpha_{k+1}}}
\newcommand\Tl{T(\lambda)}
\newcommand\Tlt{\tilde\Tl}
\newcommand\Tltp{\tilde T'(\lambda)}
\newcommand\Tltpp{\tilde T''(\lambda)}
\newcommand\Tli{T_i(\lambda)}
\newcommand\Tlit{\tilde\Tli}
\newcommand\Tlip{T_i'(\lambda)}
\newcommand\Tlipp{T_i''(\lambda)}
\newcommand\Tlj{T_j(\lambda)}
\newcommand\Tla{T_{\alpha}(\lambda)}
\newcommand\Tlaa{T_{\alpha_1}(\lambda)}
\newcommand\Tlab{T_{\alpha_2}(\lambda)}
\newcommand\Tlak{T_{\alpha_k}(\lambda)}
\newcommand\Tlakt{\tilde\Tlak}
\newcommand\Tlaj{T_{\alpha_j}(\lambda)}
\newcommand\Tlajj{T_{\alpha_{j+1}}(\lambda)}
\newcommand\Tlajp{T_{\alpha_j}'(\lambda)}
\newcommand\Tlajpt{\tilde\Tlajp}
\newcommand\Tlajt{\tilde\Tlaj}
\newcommand\Tlakk{T_{\alpha_{k+1}}(\lambda)}
\newcommand\Tlakkp{T_{\alpha_{k+1}}'(\lambda)}
\newcommand\Tlap{T_{\alpha'}(\lambda)}
\newcommand\Tlapt{\tilde\Tlap}
\newcommand\Tlapp{T_{\alpha''}(\lambda)}
\newcommand\Tkl{T^{(k)}(\lambda)}
\newcommand\Tcl{T^{\flat}(\lambda)}
\newcommand\Fl{F(\lambda)}
\newcommand\BlVt{\tilde B_V(\lambda)}
\newcommand\KBlVt{K_{\BlVt}}
\newcommand\BlVaat{B_{V_{\alpha_1}}(\lambda)}
\newcommand\BV{B_V}
\newcommand\Bone{B_1}
\newcommand\Btwo{B_2}
\newcommand\Bthree{B_3}
\newcommand\Banyj{B_j}
\newcommand\PlV{P_V(\lambda)}
\newcommand\PlVc{P_V^{\flat}(\lambda)}
\newcommand\Pl{P_0(\lambda)}
\newcommand\SVl{S_V(\lambda)}
\newcommand\Sjr{S_j^{\reduced}}
\newcommand\Rkp{{\mathcal R}^k_+}
\newcommand\Rkm{{\mathcal R}^k_-}
\newcommand\Rkpm{{\mathcal R}^k_{\pm}}
\newcommand\Phys{{\mathcal P}}
\newcommand\Pc{\overline{\mathcal P}}
\newcommand\pip{\pi^{\perp}}
\newcommand\pipa{\pi_\partial}
\newcommand\gammapa{\gamma_\partial}
\newcommand\pipah{\hat\pi_\partial}
\newcommand\pit{\tilde\pi}
\newcommand\xit{\tilde\xi}
\newcommand\zetat{\tilde\zeta}
\newcommand\etat{\tilde\eta}
\newcommand\sigmat{\tilde\sigma}
\newcommand\sigmahat{\hat\sigma}
\newcommand\thetat{\tilde\theta}
\newcommand\psit{\tilde\psi}
\newcommand\phit{\tilde\phi}
\newcommand\chit{\tilde\chi}
\newcommand\rhot{\tilde\rho}
\newcommand\xib{\bar\xi}
\newcommand\zetab{\bar\zeta}
\newcommand\thetab{\bar\theta}
\newcommand\etab{\bar\eta}
\newcommand\iotal{\iota_{\lambda}}
\newcommand\rhoat{\rhot_{\alpha_1}}
\newcommand\Lambdat{\tilde\Lambda}
\newcommand\Lambdati{\tilde\Lambda^{\text{in}}}
\newcommand\Lambdato{\tilde\Lambda^{\text{out}}}
\newcommand\Lambdatp{\tilde\Lambda^{\text{prop}}}
\newcommand\Lambdai{\Lambda^{\text{in}}}
\newcommand\Lambdao{\Lambda^{\text{out}}}
\newcommand\poles{\Lambda'}
\newcommand\rpoles{\Lambda_p}
\newcommand\thresholds{\Lambda}
\newcommand\Vt{\tilde V}
\newcommand\It{\tilde I}
\newcommand\half{{\frac{1}{2}}}
\newcommand\sigmah{\sigma^{1/2}}
\newcommand\bX{\partial X}
\newcommand\bXb{\partial \Xb}
\newcommand\Deltabt{\tilde\Delta_0}
\newcommand\strip{\Omega_T}
\newcommand\Kf{K^{\flat}}
\newcommand\Gs{G^{\sharp}}
\newcommand\Gt{\tilde G}
\newcommand\Osc{\sci\Omega}
\newcommand\OSc{{}^\Scl\Omega}
\newcommand\Osch{\sci\Omega^{\half}}
\newcommand\Oscmh{\sci\Omega^{-\half}}
\newcommand\Isc{I_{sc}}
\newcommand\os{{\text{os}}}
\newcommand\Qzl{Q^0_{-\lambda}}
\newcommand\Lie{{\mathcal L}}
\newcommand\bl{{\text b}}
\newcommand\scl{{\text{sc}}}
\newcommand\sccl{{\text{scc}}}
\newcommand\Scl{{\text{Sc}}}
\newcommand\ScLl{{\text{Sc,L}}}
\newcommand\ScRl{{\text{Sc,R}}}
\newcommand\Sccl{{\text{Scc}}}
\newcommand\sus{{\text{sus}}}
\newcommand\ssl{{\text{ee}}}
\newcommand\bzl{{\text{b0}}}
\newcommand\XXb{X^2_\bl}
\newcommand\XXbt{\Xt^2_\bl}
\newcommand\XXsc{X^2_\scl}
\newcommand\XXsct{\Xt^2_\scl}
\newcommand\XXSc{X^2_\Scl}
\newcommand\XXSct{\Xt^2_\Scl}
\newcommand\XXScL{X^2_\ScLl}
\newcommand\XXScR{X^2_\ScRl}
\newcommand\MMsc{M^2_\scl}
\newcommand\Deltab{\Delta_\bl}
\newcommand\Deltasc{\Delta_\scl}
\newcommand\DeltaSc{\Delta_\Scl}
\newcommand\DeltaScL{\Delta_\ScLl}
\newcommand\DeltaScR{\Delta_\ScRl}
\newcommand\prs{\sigma}
\newcommand\Nsc{N_\scl}
\newcommand\Nscp{N_{\scl,p}}
\newcommand\Nff{N_{\ff}}
\newcommand\Nffz{N_{\ff,0}}
\newcommand\Nffzp{N_{\ff,0,p}}
\newcommand\Nffl{N_{\ff,l}}
\newcommand\Nffml{N_{\ff,-l}}
\newcommand\Nmf{N_{\mf}}
\newcommand\Nmfz{N_{\mf,0}}
\newcommand\Nmfl{N_{\mf,l}}
\newcommand\Nmfml{N_{\mf,-l}}
\newcommand\ffb{\operatorname{bf}}
\newcommand\Ffb{\operatorname{bf'}}
\newcommand\ffsc{\operatorname{sf}}
\newcommand\ffSc{\operatorname{sf_C}}
\newcommand\Ffsc{\operatorname{sf'}}
\newcommand\rff{\rho_{\ff}}
\newcommand\rmf{\rho_{\mf}}
\newcommand\rffb{\rho_{\ffb}}
\newcommand\rffsc{\rho_{\ffsc}}
\newcommand\rFfsc{\rho_{\Ffsc}}
\newcommand\rffSc{\rho_{\ffSc}}
\newcommand\rinf{\rho_{\infty}}
\newcommand\CL{C_L}
\newcommand\CR{C_R}
\newcommand\betab{\beta_\bl}
\newcommand\betasc{\beta_\scl}
\newcommand\betaSc{\beta_\Scl}
\newcommand\BetaSc{\bar\beta_\Scl}
\newcommand\betaScL{\beta_\ScLl}
\newcommand\betaScR{\beta_\ScRl}
\newcommand\ScT{{}^\Scl T^*}
\newcommand\SccT{{}^\Scl \bar T^*}
\newcommand\ScS{{}^\Scl S^*}
\newcommand\Tb{{}^\bl T}
\newcommand\Tss{{}^\ssl T}
\newcommand\Tsc{{}^\scl T}
\newcommand\TSc{{}^\Scl T}
\newcommand\CSc{C_\Scl}
\newcommand\Lambdasc{{}^\scl\Lambda}
\newcommand\XXXb{X^3_\bl}
\newcommand\XXXsc{X^3_\scl}
\newcommand\XXXSc{X^3_\Scl}
\newcommand\XXXScO{X^3_{\Scl,O}}
\newcommand\XXXScF{X^3_{\Scl,F}}
\newcommand\XXXScS{X^3_{\Scl,S}}
\newcommand\XXXScC{X^3_{\Scl,C}}
\newcommand\KDsc{\operatorname{KD^{\half}_\scl}}
\newcommand\KDSc{\operatorname{KD^{\half}_\Scl}}
\newcommand\KDScEF{\operatorname{KD^{E,F}_\Scl}}
\newcommand\Oh{\operatorname{\Omega^{\half}}}
\newcommand\WFSc{\WF_\Scl}
\newcommand\WFtSc{\WF_{\text 3sc}}
\newcommand\WFScmf{\WF_{\Scl,\mf}}
\newcommand\WFScff{\WF_{\Scl,\ff}}
\newcommand\WFScs{\WF_{\Scl,\prs}}
\newcommand\WFScp{\WF'_\Scl}
\newcommand\WFScmfp{\WF'_{\Scl,\mf}}
\newcommand\WFScffp{\WF'_{\Scl,\ff}}
\newcommand\WFScsp{\WF'_{\Scl,\prs}}
\newcommand\Diffscc{\Diff_\sccl}
\newcommand\DiffSc{\Diff_\Scl}
\newcommand\Diffss{\Diff_\ssl}
\newcommand\DiffScc{\Diff_\Sccl}
\newcommand\DiffscI{\Diff_{\scl,\text{I}}}
\newcommand\VscI{\Vf_{\scl,\text{I}}}
\newcommand\DiffsV{\operatorname{Diff}_{\sus(V)}}
\newcommand\DiffsVsc{\operatorname{Diff}_{\sus(V),\scl}}
\newcommand\DiffsVCsc{\operatorname{Diff}_{\sus(V)-C,\scl}}   
\newcommand\Psisc{\Psop_\scl}
\newcommand\Psiscc{\Psop_\sccl}
\newcommand\Psiss{\Psop_\ssl}
\newcommand\Psisch{\Psop_{\scl,h}}
\newcommand\Psiscch{\Psop_{\sccl,h}}
\newcommand\PsiSc{\Psop_\Scl}
\newcommand\PsiScph{\Psop_{\Scl,\phi}}
\newcommand\PsiScra{\Psop_{\Scl,\rho^\sharp_a}}
\newcommand\PsiScc{\Psop_\Sccl}
\newcommand\PsiSccml{\Psop^{m,l}_\Sccl}
\newcommand\PsiScxx{\Psop^{*,*}_\Scl}
\newcommand\PsiScml{\Psop^{m,l}_\Scl}
\newcommand\PsiScmz{\Psop^{m,0}_\Scl}
\newcommand\PsiScmmz{\Psop^{-m,0}_\Scl}
\newcommand\PsiSckz{\Psop^{k,0}_\Scl}
\newcommand\PsiScmmml{\Psop^{-m,-l}_\Scl}
\newcommand\Psiscmkk{\Psop^{-k,k}_\scl}
\newcommand\Psiscmmmkk{\Psop^{-m-k,k}_\scl}
\newcommand\Psiscmoo{\Psop^{-1,1}_\scl}
\newcommand\Psiscmz{\Psop^{m,0}_\scl}
\newcommand\Psiscmmz{\Psop^{-m,0}_\scl}
\newcommand\PsiSckmkl{\Psop^{km,kl}_\Scl}
\newcommand\PsiScmplp{\Psop^{m',l'}_\Scl}
\newcommand\PsiScmmpllp{\Psop^{m+m',l+l'}_\Scl}
\newcommand\Psiscml{\Psop^{m,l}_\scl}
\newcommand\PsiScid{\Psop^{0,0}_\Scl}
\newcommand\PsiSczo{\Psop^{0,1}_\Scl}
\newcommand\PsiScmii{\Psop^{-\infty,\infty}_\Scl}
\newcommand\PsiScmiz{\Psop^{-\infty,0}_\Scl}
\newcommand\PsiScmoo{\Psop^{-1,1}_\Scl}
\newcommand\PsisCid{\Psop^{0,0}_{\scl-C}}
\newcommand\PsisC{\Psop_{\scl-C}}
\newcommand\Psiinf{\Psop_{\infty}}
\newcommand\Psiinfid{\Psop_{\infty}^0}
\newcommand\PsiFinf{\Psop_{\infty-\Fr}}
\newcommand\PsisVscml{\Psop^{m,l}_{\sus(V),\scl}}
\newcommand\PsisVsc{\Psop_{\sus(V),\scl}}
\newcommand\PsisVpsc{\Psop_{\sus(V_p),\scl}}
\newcommand\PsisVCSc{\Psop_{\sus(V)-C,\scl}}
\newcommand\SFinf{S_{\infty-\Fr}}
\newcommand\YsVC{Y^2_{\sus(V)-C,\scl}}
\newcommand\ffYsc{\ffsc_{\sus(V)}}
\newcommand\SXC{S(X;C)}
\newcommand\Ios{I_{\text{os}}}
\newcommand\pbL{\pi^2_{\bl,{\text L}}}
\newcommand\pbR{\pi^2_{\bl,{\text R}}}
\newcommand\pscL{\pi^2_{\scl,{\text L}}}
\newcommand\pscR{\pi^2_{\scl,{\text R}}}
\newcommand\PbO{\pi^3_{\bl,{\text O}}}
\newcommand\PscO{\pi^3_{\scl,{\text O}}}
\newcommand\PScO{\pi^3_{\Scl,{\text O}}}
\newcommand\PScF{\pi^3_{\Scl,{\text F}}}
\newcommand\PScC{\pi^3_{\Scl,{\text C}}}
\newcommand\PScS{\pi^3_{\Scl,{\text S}}}
\newcommand\pScL{\pi^2_{\Scl,{\text L}}}
\newcommand\pScR{\pi^2_{\Scl,{\text R}}}
\newcommand\CLF{\CL^F}
\newcommand\CLO{\CL^O}
\newcommand\CLS{\CL^S}
\newcommand\CLC{\CL^C}
\newcommand\DeltaYb{\Delta_{\bl,Y}}
\newcommand\DeltaYsc{\Delta_{\sus-\scl}}
\newcommand\diag{\operatorname{diag}}
\newcommand\Diag{\operatorname{Diag}}
\newcommand\Vf{{\mathcal V}}
\newcommand\Vb{{\mathcal V}_{\bl}}
\newcommand\Vsc{{\mathcal V}_{\scl}}
\newcommand\Vss{{\mathcal V}_{\ssl}}
\newcommand\VSc{{\mathcal V}_{\Scl}}
\newcommand\VfI{\Vf_{\text{I}}}
\newcommand\VfIq{\Vf_{\text{I},q}}
\newcommand\scH{{}^\scl H}
\newcommand\scHg{\scH_g}
\newcommand\Hss{H_\ssl}
\newcommand\Hbz{H_\bzl}
\newcommand\xh{\hat x}
\newcommand\Yh{\hat Y}
\newcommand\Zh{\hat Z}
\newcommand\Yb{\bar Y}
\newcommand\hb{\bar h}
\newcommand\xih{\hat\xi}
\newcommand\etah{\hat\eta}
\newcommand\muh{\hat\mu}
\newcommand\mub{\bar\mu}
\newcommand\nub{\bar\nu}
\newcommand\mubh{\widehat{\bar\mu}}
\newcommand\yb{\bar y}
\newcommand\zb{\bar z}
\newcommand\ub{\bar u}
\newcommand\Qb{\bar Q}
\newcommand\Wbp{{\bar W}^\perp}
\newcommand\Wp{W^\perp}
\newcommand\Kt{\tilde K}
\newcommand\Wt{\tilde W}
\newcommand\Ut{\tilde U}
\newcommand\yt{\tilde y}
\newcommand\ut{\tilde u}
\newcommand\vt{\tilde v}
\newcommand\ft{\tilde f}
\newcommand\htil{\tilde h}
\newcommand\St{\tilde S}
\newcommand\Pt{\tilde P}
\newcommand\Rt{\tilde R}
\newcommand\qt{\tilde q}
\newcommand\Qt{\tilde Q}
\newcommand\Xb{\overline{X}}
\newcommand\lambdat{\tilde\lambda}
\newcommand\betat{\tilde\beta}
\newcommand\epst{\tilde\epsilon}
\newcommand\ep{\epsilon}
\newcommand\bt{\tilde b}
\newcommand\Xt{\widetilde X}
\newcommand\Mt{\widetilde M}
\newcommand\Mb{\overline{M}}
\newcommand\At{\tilde A}
\newcommand\Et{\tilde E}
\newcommand\Ht{\tilde H}
\newcommand\at{\tilde a}
\newcommand\Ct{\tilde C}
\newcommand\pih{\hat\pi}
\newcommand\Rh{\hat R}
\newcommand\Ah{\hat A}
\newcommand\Bh{\hat B}
\newcommand\Ch{\hat C}
\newcommand\Gh{\hat G}
\newcommand\Hh{\hat H}
\newcommand\Qh{\hat Q}
\newcommand\Ph{\hat P}
\newcommand\Nh{\hat N}
\newcommand\Sh{\hat S}
\newcommand\Gcal{{\mathcal G}}
\newcommand\GcalC{{\mathcal G}_C}
\newcommand\Jcal{{\mathcal J}}
\newcommand\JcalC{{\mathcal J}_C}
\setcounter{secnumdepth}{3}
\newtheorem{lemma}{Lemma}[section]
\newtheorem{prop}[lemma]{Proposition}
\newtheorem{thm}[lemma]{Theorem}
\newtheorem{cor}[lemma]{Corollary}
\newtheorem{result}[lemma]{Result}
\newtheorem*{thm*}{Theorem}
\newtheorem*{prop*}{Proposition}
\newtheorem*{conj*}{Conjecture}
\numberwithin{equation}{section}
\theoremstyle{remark}
\newtheorem{rem}[lemma]{Remark}
\theoremstyle{definition}
\newtheorem{Def}[lemma]{Definition}
\newtheorem*{Def*}{Definition}
\def\signature#1#2{\par\noindent#1\dotfill\null\\*
{\raggedleft #2\par}}
\newcommand\ess{\mbox{ess}}

\renewcommand{\theenumi}{\roman{enumi}}
\renewcommand{\labelenumi}{(\theenumi)}

\title[Analytic continuation of the resolvent] 
{Analytic continuation \\ of the resolvent of the Laplacian \\
on symmetric spaces of noncompact type}
\author[Rafe Mazzeo and Andras Vasy]{Rafe Mazzeo and Andr\'as Vasy}
\address{R.\ M.: Department of Mathematics, Stanford University, Stanford,
CA 94305}
\email{mazzeo@math.stanford.edu}
\address{A.\ V.: Department of Mathematics, Massachusetts Institute of
Technology, MA 02139}
\email{andras@math.mit.edu}
\date{August 5, 2003}

\begin{abstract}
Let $(M,g)$ be a globally symmetric space of noncompact type, of
arbitrary rank, and $\Delta$ its Laplacian. We prove the existence
of a meromorphic continuation of the resolvent $(\Delta-\ev)^{-1}$ 
across the continuous spectrum to a Riemann surface multiply covering 
the plane. The methods are purely analytic and are adapted from
quantum $N$-body scattering.
\end{abstract}

\maketitle

\section{Introduction}
A basic problem in geometric scattering theory is to carry out a refined analysis of 
the resolvent of the Laplacian on various classes of complete manifolds with regular 
geometry at infinity. The symmetric spaces of noncompact type comprise a natural 
class of manifolds to understand from this point of view because their asymptotic 
geometry is so well understood. An added attraction is that the analytic properties
of the Laplacians on these spaces are closely connected to representation theory and 
number theory.  In this paper we continue our program, initiated in \cite{Mazzeo-Vasy:Resolvents}, 
to extend the methods and results of geometric scattering theory to this setting. 
More specifically, 
let $M = G/K$ be a symmetric space of noncompact type, with $\mbox{rank}(M) = n$, and denote 
by $\Delta = \Delta_M$ its Laplace-Beltrami operator with respect to some choice of 
invariant metric. We do not assume that $M$ is irreducible, so any such 
metric is obtained by fixing a constant multiple of the Killing form on each irreducible
factor. As $M$ is complete, $\Delta$ is self-adjoint.
The resolvent of the Laplacian is the operator $R(\ev) = (\Delta - \ev)^{-1}$, initially
defined when $\ev \in \Cx \setminus [0,\infty)$ as a bounded operator on $L^2(M)$. 
In this paper we prove that $R(\ev)$ continues meromorphically to a larger set. 
The existence of this continuation is classical when $M$ is a Euclidean space, 
and is also well known for rank one symmetric spaces and their 
geometric generalizations, e.g.\ conformally compact spaces 
\cite{Mazzeo-Melrose:Meromorphic} and their complex analogues 
\cite{Epstein-Melrose-Mendoza}; it is also known in the case of higher rank 
{\it complex} symmetric spaces, but surprisingly, its existence for higher 
rank {\it real} symmetric spaces is only known indirectly \cite{FlJ}. Recently we used techniques from
microlocal analysis to prove this continuation in the two simplest rank $2$ situations: when 
$M$ is a product of hyperbolic spaces \cite{Mazzeo-Vasy:Resolvents} and when $M = 
\SL(3)/\SO(3)$ \cite{Mazzeo-Vasy:Sl3}, \cite{Mazzeo-Vasy:Sl3-Analytic}, and our goal 
in this paper is to extend that construction to the general case.
Let $G_o(\ev)$ denote the Green function, i.e.\ the Schwartz kernel 
of $R(\ev)$. This is our main result:

\begin{thm}\label{thm:main}
The Green function $G_o(\ev)$ continues meromorphically as a distribution
to a Riemann surface $\widetilde{\calY}_{\pi/2}$
(see Definition~\ref{Def:calY-pi-2}),
ramified at a sequence of points corresponding to translates of the poles 
of the meromorphic continuation of $G_o(\ev)$ on symmetric spaces of lower rank. 
\end{thm}

It is then natural to ask whether these poles exist. Here we show that 
they lie in a compact set in the complement of any cone containing 
a singular direction; in fact, an estimate which implies this plays an
important role in the proof of the existence of the continuation. However,
we conjecture that this continuation has no poles at all on $\widetilde{\calY}_{\pi/2}$;
see the remark at the end of the last section.

\begin{figure}[ht]
\begin{center}
\mbox{\epsfig{file=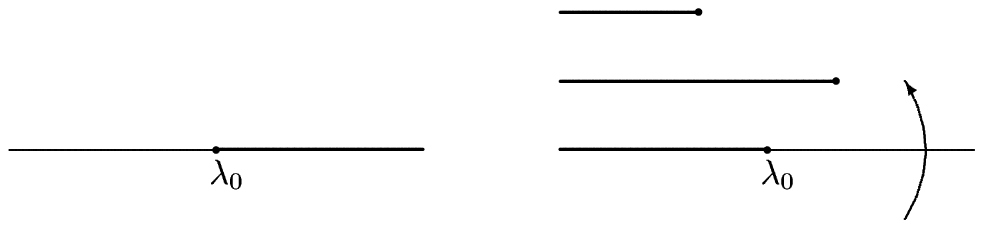}}
\end{center}
\caption{Part of the domain of the analytic continuation of $G_o(\ev)$.}
\label{fig:acrp1}
\end{figure}

We sketch part of $\widetilde{\calY}_{\pi/2}$ in Figure~\ref{fig:acrp1}. The thick line in
picture on the left shows the spectrum of $\Delta$ (inside $\Cx$). It is a half 
line $[\ev_0,\infty)$, and the Green function $G_o(\ev)$ is defined a priori 
for $\ev \in \Cx \setminus [\ev_0,\infty)$. The thin line indicates the rest of 
the real axis.

The picture on the right shows the analytic continuation of $G_o(\ev)$ 
across the ray $[\ev_0,\infty)$ from below; 
it is defined {\em outside} the thick half lines, one of which is
$(-\infty,\ev_0]$. The thin line is again the rest of the real axis.
Thus, for $\ev\in\Cx$ with $\im\ev<0$, $G_o(\ev)$ is defined identically in
the two pictures, but for $\im\ev>0$, on the right hand side, $G_o(\ev)$
lives on a different sheet of the Riemann surface, whose projection to
$\Cx$ is shown. The ramification points
are indicated by the thickened points; the conjecture then is that none of
these exist except $\ev_0$.

We proceed by induction on the rank of the symmetric space. The two key 
ingredients of the proof are complex scaling, and the construction of a 
parametrix, i.e.\ an approximate inverse, for the complex scaled
$K$-radial Laplacian. 
This method is closely related to the analogous problem in $N$-body scattering, 
where it was introduced by Balslev and Combes 
\cite{Balslev-Combes:Spectral} and extended by C.\ G\'erard \cite{Gerard:Distortion}. Indeed, technically the only reason we cannot use
the N-body results directly is that if we identify $\Delta$ acting
on $K$-invariant functions with a differential operator on a flat $\calA
=\exp(\fraka)$, and hence on $\fraka$,
the $L^2$ space on $\fraka$ is
not the Euclidean one, and the first order terms are singular
at the walls of the Weyl chambers.

Complex scaling in this setting is induced by dilation along geodesic rays
from $o$. These are the maps $\Phi_\theta$ that, for $\theta\in\RR$,
send any point $\gamma(t)$ on any geodesic $\gamma$ with $\gamma(0)=o$ to
the point $\gamma(e^\theta t)$. Complex scaling extends these analytically
in $\theta$ to a domain in the complex plane. The virtue of the scaling is
that, for complex values of $\theta$, the essential spectrum of the scaled 
radial Laplacian is (almost) a rotation of the essential spectrum of the Laplacian, and
this allows the analytic continuation of the resolvent. We define and describe
the scaling here in \S 5, and we refer to the introduction of 
\cite{Mazzeo-Vasy:Sl3-Analytic} for a brief description of this procedure
for the Laplacian on the hyperbolic plane.

Although the other ingredient, 
the parametrix construction, is fundamentally microlocal, we minimize 
the explicit use of microlocal techniques, which is possible because of
the essentially `soft' nature of such an analytic continuation result,
and because there are finitely many
local `product models' for the scaled radial Laplacian $\Delta_{\rad,\theta}$,
i.e.\ locally (in certain neighbourhoods of infinity) this operator has the 
form $A\otimes\Id+\Id\otimes B$ modulo decaying error terms.
More delicate questions concerning the precise asymptotic behaviour of the 
Green function may be approached using an elaboration of the same construction,
as in \cite{Mazzeo-Vasy:Resolvents}, \cite{Mazzeo-Vasy:Sl3}, but do require
more attention to the microlocal aspects; we shall return to this elsewhere.

Although our analysis seems to make essential use of various compactifications 
of $M$, in fact these are not truly essential. Rather, they are very 
helpful in the construction of certain partitions of unity, on the support 
of which $\Delta_{\rad,\theta}$ is particularly well approximated by product 
models. Such partitions of unity could also be described by requiring various
homogeneity properties, but in the further development of the scattering theory
on symmetric spaces, e.g.\ in the study of the asymptotics of the Green function,
these compactifications play a central role.  

We would also like to underline that it is crucial that the product models
for $\Delta_{\rad,\theta}$ are valid in {\em conic} subsets of $\fraka$
-- in the language of compactifications, this is the reason we
use a partition of unity and cutoffs
on the radial (or geodesic) compactification
$\hat\fraka$. The conic cutoffs give decaying error terms in the
parametrix construction; this would not be the case if we localized at
finite distances from Weyl chamber walls.

Finally, this work would not be complete without commenting on its
relationship to the meromorphic continuation of Harish-Chandra's $c$\,-function.
The $c$\,-function is known to have a meromorphic extension to the
flat $\fraka^*_\Cx$, and its restriction to the vectors in $\fraka^*$ with
length $\sqrt{\ev-\ev_0}$ can be thought of as a `scattering matrix'
by analogy with both the rank-one case and $N$-body scattering.
Now, in the latter settings, the poles of the meromorphic continuation
of the scattering matrix (considered as an operator) and the resolvent
coincide -- one might expect that if the $c$\,-function is analytic 
on the rotation of this sphere around $0$ in $\fraka^*_\Cx$ by angle 
$\arg\sqrt{\ev-\ev_0}$, then the continuation of the resolvent does not have
a pole at $\ev$, and conversely. There is an explicit formula for the
$c$\,-function, see \cite[Chapter~IV, Theorem~6.14]{Helgason:Groups},
and it is apparent from it that this requirement on the $c$\,-function
is never satisfied in the higher rank setting (since some inner product may
vanish). This phenomenon already occurs in the 
complex case, when the formula is simpler
\cite[Chapter~IV, Theorem~5.7]{Helgason:Groups}. We expect, however,
that, suitably renormalized and considered as an operator, the 
meromorphic structure of the $c$\,-function can be related to that of 
the resolvent. 

We are very grateful to Gilles Carron, Sigurdur Helgason, Lizhen Ji,
Richard Melrose, David Vogan and Maciej Zworski
for helpful discussions and their encouragement.
R.\ M.\ is partially supported by NSF grant \#DMS-0204730;
A.\ V.\ is partially supported by NSF grant \#DMS-0201092 
and a Fellowship from the Alfred P.\ Sloan Foundation.

\section{Compactifications of $\fraka$ and the radial Laplacian}
In this section, we begin by reviewing some well-known facts about the Lie-theoretic 
algebra and global geometry of the symmetric space $M$; we refer to \cite{Helgason:DS}, 
\cite{Helgason:Groups} for a comprehensive development and all proofs, and also to 
\cite{Eberlein:Geometry} for a detailed summary 
from a more geometric point of view. Of central importance here is the flat
$\calA=\exp(\fraka)$; $\fraka$
is a Euclidean space of dimension $\mbox{rank}(M)$, and
it is the ultimate locus
of our analysis. We shall systematically identify $\fraka$ with
its exponential, and will usually work on $\fraka$ rather than $\calA$,
since it is more customary to use 
linear coordinates rather than their exponentials. 
We go on to define two compactifications of this flat, $\ofraka$ and the
larger one $\tfraka$, which play a central role in our approach. Motivation for these 
definitions is provided by the specific form of the radial Laplacian $\Delta_\rad$ on $M$, 
which is introduced and discussed along the way. We conclude by showing that the radial 
Laplacian on symmetric spaces of lower rank appear in the restrictions of this operator to 
boundary faces of $\tfraka$. 

\subsection{Geometry of flats}
Suppose $M = G/K$, and let $\frakg = \frakk + \frakp$ be the Cartan decomposition. 
Thus $\frakk$ is the Lie algebra of $K$ and $\frakp$ its orthogonal complement with
respect to the Kiling form, which is identified with $T_oM$ ($o$ will always denote 
the identity coset). We also fix a maximal abelian subspace $\fraka \subset 
\frakp$; this is always of the form $\frakp \cap \frakg_0$, where $\frakg_0$
is a maximal abelian subalgebra (called a Cartan subalgebra) in $\frakg$, and
conversely, any such intersection is a maximal abelian subspace in $\frakp$. 
The number $n :=\dim \fraka$ is called the rank of $M$, and $\exp\fraka := \calA$ 
is a totally geodesic flat submanifold which is maximal with respect to this property, 
and is called a flat. It is isometric to $\RR^n$.  

A key example, to which we shall refer back repeatedly throughout this paper
for purposes of illustration, is $\calM_{n+1} = \SL(n+1)/\SO(n+1)$. Here $\frakg = 
\sll(n+1)$ consists of all ($n+1$)-by-($n+1$) matrices of trace zero, and $\frakk 
= \sol(n+1)$ and $\frakp$ consist of all such matrices which are skew-symmetric, 
respectively symmetric. We may take $\fraka$ to be the subspace of diagonal matrices 
of trace zero. Denoting these diagonal entries by $t_i$, $i = 1, \ldots, n+1$, then 
the diagonal matrices $A_i$, $i = 1, \ldots, n$, with $t_i = 1$, $t_{i+1} = -1$ and 
all other $t_j = 0$ comprise the standard basis of $\fraka$. We identify 
$\calM_{n+1}$ with the space of positive definite symmetric matrices via the identification
$\SL(n+1) \ni B \to \sqrt{B^t B}$.  The flat $\calA = \exp(\fraka)$ consists of 
diagonal matrices with positive entries $\lambda_1, \ldots, \lambda_{n+1}$ and 
determinant $1$.

Since $\fraka$ is abelian, there is a simultaneous diagonalization for the commuting family 
of symmetric
homomorphisms $\ad\,H$, $H \in \fraka$, on $\frakg$.
A simultaneous eigenvector $X$ 
satisfies $(\ad\, H)(X) = \al(X)$ for every $H \in \fraka$, for some element 
$\al \in \fraka^*$; the set of linear forms which arise in this way constitute the (finite) 
set of (restricted)
roots $\Lambda$ for $\frakg$, and the space of eigenvectors associated to each
$\al \in \Lambda$ is the `root space' $\frakg_\al$. Thus in particular $0 \in \Lambda$
and its root space $\frakg_0$ is the Cartan subalgebra above (i.e.\ if we fix $\fraka$ 
first, then a Cartan subalgebra is uniquely associated in this way), and $\frakg = 
\oplus_{\al \in \Lambda} \frakg_\al$. We shall always use the restriction 
of the Killing form of $\frakg$ to $\frakp$ as the inner product $\langle \cdot , 
\cdot \rangle$ (rather than allowing for different scalar multiples of the Killing
form on different factors in a decomposition into irreducible subalgebras). 
This determines the root vectors $H_\al \in \fraka$ by the relationship 
$\al(H) = \langle H, H_\al \rangle$ for all $H \in \fraka$. We also fix a partition 
$\Lambda = \Lambda^+ \cup \Lambda^-$, $\Lambda_- = -\Lambda_+$, into positive and negative 
roots. There is a subset $\Lambda^+_\ind \subset \Lambda^+$ of indecomposible (or simple) 
positive roots which is a basis for $\fraka^*$ (so in particular, $\#\Lambda^+_\ind = n$) 
such that for any $\al \in \Lambda$,
\[
\al = \sum_{\al_j \in \Lambda^+_\ind} n_j \al_j, \qquad \mbox{where all}\ n_j \in 
{\mathbb Z}\quad
\mbox{and} \quad \left\{\begin{array}{rcl} & \mbox{all }n_j \geq 0 \qquad &\mbox{if}\ \al 
\in \Lambda^+ \\ & \mbox{all }n_j \leq 0 \qquad &\mbox{if}\ \al \in \Lambda^-. \end{array}
\right.
\]
Of particular importance is the element
\begin{equation}
\rho = \frac12 \sum_{\al \in \Lambda^+} m_\al\,\al \in \fraka^*,
\label{eq:casimir}
\end{equation}
where $m_\al = \dim \frakg_\al$, and its metrically dual vector $H_\rho \in \fraka$. 

Each $\al \in \Lambda$ determines a hyperplane $W_\al = \al^{-1}(0) \subset \fraka$,
called the Weyl chamber wall associated to $\al$, and by definition
\[
\fraka_\reg = \fraka \setminus \left(\bigcup_{\al \in \Lambda} W_\al\right)
\]
is called the set of {\it regular} vectors; the components of this set are called
(open) Weyl chambers, and the distinguished component 
\[
C^+ = \{H \in \fraka: \al(H) > 0, \quad \forall \, \al \in \Lambda^+\},
\]
is called the positive Weyl chamber. We also define
\[
W_{\al,\reg} = W_\al \setminus \left(\bigcup_{\beta\neq\alpha}
(W_\beta \cap W_\al) \right).
\]
As already indicated,
we shall systematically identify each of these sets with their corresponding exponentials 
in $\calA$: in particular, set $\calA_\reg = \exp(\fraka_\reg)$, $\exp(W_\al) = 
\calW_\al$, $\calW_{\al,\reg} = \exp(W_{\al,\reg})$ and $\exp(C^+) = \calC^+$.

The orthogonal reflections across the Weyl chamber walls generate a finite group, called
the Weyl group $W$. Alternately, $W$ is the quotient $N(\fraka)/Z(\fraka)$
of the normalizer by the centralizer of $\fraka$ with respect to the adjoint action 
$\mbox{Ad}$ of $K$ on $\frakg$. The Weyl group acts simply transitively on the set 
of Weyl chambers. 

Returning again to the special case $M = \calM_{n+1}$, the root set $\Lambda$ consists
of all $\al_{ij}$, where for the diagonal matrix $T = \mbox{diag}(t_1, \ldots, t_{n+1})$, 
$\al_{ij}(T) = t_i - t_j$. We take $\Lambda^+= \Lambda^+_\ind = \{\al_{i+1\, i}, 1 \leq i 
\leq n\}$; so that the positive Weyl chamber $C^+$ consists of all traceless diagonal 
matrices $A$ with all $t_1 < t_2 \ldots < t_{n+1}$, while $\calC^+$ consists of all 
unimodular diagonal matrices such that $0 < \la_1 < \ldots < \la_{n+1}$. The centralizer 
$Z(\fraka)$ in $\SO(n+1)$ is the set of diagonal matrices with entries equal to $\pm 1$, 
while the normalizer $N(\fraka)$ in $\SO(n+1)$ is the set of signed permutation matrices, 
and so the Weyl group $W$ is identified with the symmetric group 
$S_{n+1}$, and acts by permutations on the entries of the diagonal matrices.

$G$ acts on $M = G/K$ by left multiplication. The Cartan decomposition states
that $G = K \cdot \calA \cdot K$, and in stronger form, $G = K \cdot 
\overline{\calC^+} \cdot K$. Moreover, for $g\in G$, with $g=k_1 a k_2$, the element
$a\in\overline{\calC^+}$, as well as $H \in C^+$ satisfying
$a=\exp H$, are uniquely determined; we write $H=H(g)$. This induces a map on $M$, 
so for $p=gK\in M$, $H(p)=H(g)$.

The geodesic exponential map $\exp:\frakp\to M$ is a diffeomorphism.
Moreover, $k\cdot\exp(X)=\exp(\Ad(k) X)$ for $k\in K$, $X\in\frakp$.

Letting $G_{\reg}=K\calA_{\reg}K = K \calC^+ K$ and $M_{\reg} =G_{\reg}\cdot o$, 
then $M_{\reg}$ is diffeomorphic to $K'\times \calC^+$, where $K' = K/Z(\calA)$, 
see \cite[Ch. IX, Corollary 1.2]{Helgason:DS}. In fact, $K'$ acts freely 
on $\calA_\reg$, but if $X \in \calA \setminus \calA_\reg$, then the isotropy group
$K^X \subset K$ is strictly larger than $Z(\calA)$. Fixing a root $\al$, then 
all the isotropy groups $K^X$ for $X \in \calW_{\al,\reg}$ are the same, and 
we denote this common group by $K^\al$. There is a 
larger subgroup $K^\calW \subset K$ which maps $\calA \setminus \calA_\reg$ to itself 
(and hence permutes the Weyl chamber walls). The entire symmetric space is obtained as 
the quotient of $K' \times \overline{\calC^+}$ by the diagonal Weyl group action. 

Following the last paragraph, we see that elements of $\calC^\infty(M)^K$,
the space of smooth $K$-invariant 
functions on $M$, restrict to elements of $\calC^\infty(\calA)^W$,
the space of smooth $W$-invariant
functions on $\calA$; we later show in Proposition~\ref{prop:A-M-W-K}
that this map is an isomorphism.
More generally, we shall use the notation that if $E$ is any space of 
functions (on $M$ or $\calA$ or any other related space) and if $\Gamma$ is
a group on the underlying space, then $E^\Gamma$ is the subspace of $\Gamma$-invariant
elements. 

\subsection{The radial Laplacian} 
Before proceeding with further geometric considerations, we now introduce the radial 
Laplacian $\Delta_\rad$, which is simply the restriction the full Laplacian $\Delta_M$ 
to $K$-invariant functions (or distributions) on $M$. $\Delta_\rad$ is our principal 
object of study in this paper, and the main task ahead of us is the construction of 
parametrices for $(\Delta_\rad - \lambda)^{-1}$.

Rather than thinking of the radial Laplacian as an operator on $M$, acting
on a restricted space of functions, it is more useful to realize $\Delta_\rad$ 
as an operator acting on essentially arbitrary functions on a lower dimensional 
manifold. This is done by restricting to functions on a submanifold
transverse to the orbits of $K$ on $M$, and the simplest choice is to
restrict to the regular part of the flat $\calA_\reg$, which we identify 
with $\fraka_{\reg}$. Of course, we will then have to investigate the extension
of this operator to the entire flat. 

There is an elegant expression for the radial Laplacian on $\fraka_\reg$:
\begin{equation}
\Delta_\rad =\Delta_{\fraka}+ \frac12
\sum_{\alpha\in\Lambda}(m_\alpha\, \coth\alpha)\, H_\alpha,
\label{eq:radlap}
\end{equation}
where $\Delta_\fraka$ is the standard Laplacian on the vector space
$\fraka$, $m_\alpha 
= \dim \frakg_\alpha$ and $H_\al$ is the root vector associated to the root $\al$, as
defined in \S 2.1. Noting that $m_\alpha=m_{-\alpha}$, $\coth(-\alpha)=-\coth\alpha$ and 
$H_{-\alpha} =-H_{\alpha}$, we also have
\begin{equation}
\Delta_\rad =\Delta_{\fraka}+\sum_{\alpha\in\Lambda^+}(m_\alpha\, \coth\alpha)\, 
H_\alpha,
\label{eq:radlap2}
\end{equation}
which is the expression found in
\cite[Ch. II, Proposition~3.9]{Helgason:Groups}.
It is clear from (\ref{eq:radlap}) that the action of $W$ on $\fraka_\reg$ leaves 
$\Delta_\rad$ invariant. The singularities in the coefficients of these first order 
terms along the Weyl chamber walls might seem to complicate the process of
extending this operator to all of $\fraka$, and indeed this would be the
case if we were to try to let $\Delta_\rad$ act on $\Cinf(\fraka)$, for example.
However, this difficulty disappears if we restrict to $W$-invariant functions.
Indeed, we shall prove in the next section that $\calC^\infty(M)^K$ is naturally
identified with $\calC^\infty (\fraka)^W$, and so (tautologically) $\Delta_\rad$
extends to this latter space, and then also to $W$-invariant distributions, etc.
As a first step toward this identification, we prove the

\begin{lemma}\label{lemma:Delta-rad-L}
The operator $\Delta_\rad:\Cinf(\fraka_\reg)^W\to\Cinf(\fraka_\reg)^W$
induces a map $L:\Cinf(\fraka)^W\to\Cinf(\fraka)^W$ via
the inclusion $\iota:\fraka_\reg\hookrightarrow\fraka$. That is,
if $f\in\Cinf(\fraka)^W$, then $\Delta_{\rad}\iota^*f=\iota^*g$ for
some $g\in\Cinf(\fraka)^W$, and $g=Lf$ is uniquely determined by $f$.
\end{lemma}

\begin{proof}
By the density of $\fraka_{\reg}$ in $\fraka$ and the smoothness of $g$, it is clear 
that $g$ will be unique once we know it exists. To prove its existence, note first
that $\Delta_{\fraka}$ commutes with any reflection on $\fraka$, hence is
invariant by the action of $W$, and so maps $\Cinf(\fraka)^W$ to itself.
Thus it suffices to prove that the same is true for each of the summands 
$\coth\alpha\,H_\alpha$, $\alpha\in\Lambda^+$. For any $\beta \in \Lambda^+$, let 
$R^\beta$ denote the reflection across the wall $W_\beta$, and $\Cinf(\fraka)^{R_\beta}$ 
the space of functions invariant by this reflection. Writing
\[
\coth \al \, H_\al = (\al \coth \al) \frac{1}{\al}H_\al,
\]
then, since both $\al$ and $\coth \al$ are simultaneously either fixed
or taken to their negatives by any $R_\beta$, we have 
$\al \coth \al \in \Cinf(\fraka)^{R_\beta}$ for every $\beta$. Thus
we reduce at last to proving that for each $\al$ and $\beta$, $\al^{-1}H_\al$ 
maps $\Cinf(\fraka)^{R_\beta}$ to itself. But
$S^\alpha=W_\alpha^\perp=\Span(H_\al)$ is 
a copy of $\RR$ and the smooth even functions on this line are all smooth
functions of $\sigma = \al^2$, and so the operator $\al^{-1}H_\al = 2\frac{d\,}{d\sigma}$ 
certainly preserves the space of smooth even functions. Similarly, any
element $f \in \Cinf(\fraka)^{R_\beta}$ can be regarded
as a family of smooth even functions
$\tilde{f_x}$ on $S^\al$ too, as $x$ ranges over $W_\alpha$,
and the action of $\al^{-1}H_\al$ on
$f$ may be determined from the induced action on $\tilde{f}_x$.

We have proved that if $f \in \Cinf(\fraka)^W$, then there is a function 
$Lf \in \Cinf(\fraka)$ which agrees with $\Delta_\rad f$ on $\fraka_{\reg}$;
the $W$-invariance of $Lf$ follows from its $W$-invariance on the dense subset
$\fraka_\reg$.
\end{proof}

The actual identification of $\Cinf(M)^K$ with $\Cinf(\fraka)^W$ uses
this lemma, but also requires the ellipticity of $\Delta_M$, and so we
defer the proof until we have covered more preliminaries. However, 
we emphasize the conclusion, that the singularities of $\Delta_\rad$ 
are of the same nature as the singularities of the Laplacian on $\Real^n$ 
when written in polar coordinates. 

We conclude this subsection by exhibiting the many-body structure of $\Delta_\rad$ 
more plainly. Write  
\begin{equation}
\Delta_\rad = \Delta_{\fraka}+ 2H_\rho + E,
\label{eq:radlap3}
\end{equation}
where $H_\rho$ is as in (\ref{eq:casimir}), and 
\[
E = \sum_{\alpha\in\Lambda^+}m_\alpha(\coth\alpha-1) H_\alpha.
\]
The first terms, $\Delta_{\fraka}+2H_\rho$, are translation invariant, hence can be analyzed 
easily using Fourier analysis. On the other hand, each summand in $E$ is a first order 
operator which decays exponentially as the corresponding root $\al \to +\infty$. This 
rearrangement of the first order terms is only satisfactory in $C^+$, but the $W$ 
invariance of $\Delta_\rad$ implies that it is meaningful everywhere. The vectors 
$H_\al$ are not independent (except in the special, completely reducible case), and 
so (\ref{eq:radlap3}) shows that $\Delta_\rad$ has first order interaction terms of 
$N$-body type, where the finite intersections of Weyl 
chamber walls play the role of `collision planes'.

\subsection{Compactifications} Because of the many-body structure of $\Delta_\rad$, 
any thorough analysis of this operator and its resolvent must include some sort of 
delicate localization at infinity. As already explained in the introduction, the 
traditional approach of Harish-Chandra is most effective in sectors disjoint from 
the Weyl chamber walls, while uniformity of behaviour of various analytic objects 
on approach to these walls is more difficult to obtain; on the other hand, in our 
approach these walls are essentially `interior points', and create no difficulties. 
The main issue is to find and work in neighbourhoods which most effectively intermediate 
between these two types of behaviour. The use of compactifications to localize
at infinity, or at least to better visualize and control these localizations, is
well known. In the next subsections we shall introduce three main compactifications:
the first, $\hat\fraka$, is the geodesic, or radial, compactification; the second,
$\ofraka$, is known as the dual-cell compactification; the third, $\tfraka$, is the 
minimal compactification which dominates the other two. All of these have been used 
elsewhere, cf.\ \cite{Guivarch-Ji-Taylor:Compactifications}, \cite{Oshima:Realization}, 
but we shall emphasize their smooth structures; in particular our contention 
(born out by the conclusions of this paper) that $\tfraka$ is the most appropriate 
place to study $\Delta_\rad$, is a novel perspective.

As orientation for the remainder of \S 2, we sketch what lies ahead. The radial 
compactification $\hat\fraka$ is by far the simplest of the compactifications. It
is obtained either by `adding a point to the end of each geodesic', cf.\ \cite{Eberlein:Geometry},
or equivalently by completing the stereographic image of $\fraka \hookrightarrow
{\mathbb S}(\fraka \oplus \RR)$ as the closed upper hemisphere of $S^n$. 
This latter description
immediately equips $\hat\fraka$ with the structure of a smooth manifold with 
boundary. The monograph \cite{RBMGeo} contains an extended panegyric on the 
advantages of this space in the scattering analysis of the free Laplacian $\Delta_\fraka$
and its (short range) perturbations. However, the lifts of the first order terms 
in $\Delta_\rad$ to this space are not particularly simple, 
and this necessitates a slightly different approach. As a smooth manifold with
corners, the compactification $\ofraka$ is a slightly more complicated object,
but it accomodates these first order terms very nicely. It is obtained
essentially by requiring that the functions $e^{-\alpha}$ restricted to the positive
Weyl chamber extend to smooth functions on the closure of $C^+$. However, although 
the principal part $\Delta_\fraka$ lifts to a smooth $b$-operator on this space,
it does not have a product structure near the corners, even asymptotically, and
so its analysis here is still difficult. The space $\tfraka$ is the smallest
compactification for which there are smooth `blowdown maps' to both
$\hat\fraka$ and $\ofraka$, and it therefore has the property that 
{\it both} the principal part and the first order terms in $\Delta_\rad$ lift nicely
to this space. The precise sense in which we mean this will become apparent
in the discussion below. 

Through most of the ensuing discussion we tacitly assume that the root system 
$\Lambda$ spans $\fraka$. However, even if we start with a semisimple Lie algebra, 
where this is the case, we will always encounter situations in the overall induction 
on rank where $\fraka = \fraka' \oplus \fraka''$ and all roots vanish identically 
on the second summand. Therefore we must 
adapt all constructions and arguments to subsume this case too. Thus, to begin this 
generalization, the boundary of the radial compactification of $\fraka$
is a sphere, inside of which sit the boundaries of the radial compactifications
of the two summands as nonintersecting equatorial subspheres, and $\widehat{\fraka}$
is the simplicial join of these subspheres, i.e.
\begin{equation}
\pa \, \widehat{\fraka} = \pa\,\widehat{\fraka'}\  \# \ \pa\,\widehat{\fraka''}.
\label{eq:radjoin}
\end{equation}
Of course, we regard $\pa\, \widehat{\fraka}$ as a smooth (rather than a
combinatorial) manifold. 

\subsection{The compactification $\ofraka$}
The compactification $\ofraka$ is known elsewhere in the symmetric space literature 
as the polyhedral or dual-cell compactification,
see \cite[Section~3.22-3.33]{Guivarch-Ji-Taylor:Compactifications}.
It carries the natural structure of a 
polytope, i.e.\ is really a PL object, but for us it is only important that it is
a smooth manifold with corners. Briefly, $\ofraka$ is obtained by compactifying the 
positive Weyl chamber $C^+$ as a cube, $[0,1]^n$, to which the action of the Weyl group 
extends naturally; its translates by $W$ fit together affinely to generate the 
entire polytope. 

We now explain this more carefully. First fix an enumeration $\{\al_1, \ldots, \al_n\}$
of the set of positive simple roots $\Lambda^+_\ind$. This is a basis for $\fraka^*$, 
hence a maximal independent collection of linear coordinates on $\fraka$. For any $n$-tuple 
$T = (T_1, \ldots, T_n)\in\RR^n$, there is an affine isomorphism 
\begin{equation}
\calO(T) := \bigcap_{j=1}^n\alpha_j^{-1}((T_j,+\infty)) \longrightarrow 
\prod_{j=1}^n(T_j,+\infty).
\label{eq:orth}
\end{equation}
In particular, the positive Weyl chamber $C^+ = \calO((0,\ldots,0))$ corresponds to 
the standard orthant $(\RR^+)^n$. Now change variables, replacing $\al_j$ by 
$\tau_j := e^{-\al_j}$; the set $\calO(T)$ is compactified by adjoining
the faces where $\tau_j = 0$ and $\tau_j= e^{-T_j}$. Thus
\[
\calO(T) \subset \overline{\calO(T)} \equiv \prod_{j=1}^n [T_j,\infty]_{\al_j} \cong 
\prod_{j=1}^n [0, e^{-T_j}]_{\tau_j}.
\]
As already noted, $C^+ = \calO(\,\vec{0}\,)$, and so $\overline{C^+} = 
\overline{\calO(\,\vec{0}\,)}$.
By definition, the smooth structure on these sets is the minimal one which agrees 
with the standard smooth structure on $\fraka$ away from the outer boundaries and 
for which each $\tau_j$ is smooth. (Note, however, that $1/\al_j$ is {\it not}
$\Cinf$ on $\ofraka$!)

\begin{figure}[ht]
\begin{center}
\mbox{\epsfig{file=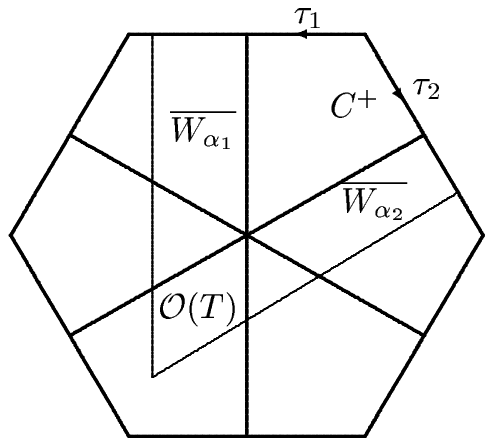}}
\end{center}
\caption{The compactification $\ofraka$ for $M=\SL(3,\RR)/\SO(3,\RR)$.
The thick lines indicate the boundary
faces and the Weyl chamber walls. The thin lines show the boundary of
$\calO(T)$ for $T_1<0$, $T_2<0$. The arrows indicate the coordinate
axes $\tau_1$ (i.e.\ $\tau_2=0$) and $\tau_2$ (i.e.\ $\tau_1=0$)
in the coordinate chart $\calO(T)$.}
\label{fig:acrp3}
\end{figure}

Any other Weyl chamber is the positive chamber for a different set of indecomposable 
roots, and so may be compactified similarly. These compactifications fit together to 
cover all of $\ofraka$. This shows that $\ofraka$ is a topological cell, and provides 
it with a smooth structure away from these patching regions at the walls. To exhibit
its structure as a smooth manifold with corners, observe that if all $T_j < 0$, then 
$\calO(T) \supsetneq C^+$, and so these neighbourhoods cover the entire space $\fraka$, 
and their completions patch together to cover all of $\ofraka$ with open overlaps. 
Thus it suffices to show that for any $w \in W$, the restriction
\begin{equation*}
w_T:w^{-1}(\calO(T))\cap \calO(T)\to \calO(T)
\end{equation*}
extends to a smooth map $\overline{w^{-1}(\calO(T))\cap \calO(T)}\to 
\overline{\calO(T)}$. For this, it is enough to prove that for any 
$\alpha_j \in \Lambda_\ind^+$, the function $w^*e^{-\alpha_j}$ extends 
smoothly to
\begin{equation*}
\overline{w^{-1}(\calO(T))\cap \calO(T)},
\end{equation*}
or equivalently,
that $w^* \tau_j$ is smooth on this set. Now, $w^*\alpha_j$ is either 
in $\Lambda^+$ or $\Lambda^-$. In the former case, it decomposes as
$\sum n_k\alpha_k$ where all $n_k$ are nonnegative integers, and so
\begin{equation*}
w^* \tau_j =\prod_k (e^{-\alpha_k})^{n_k}=\prod_k\tau_k^{n_k} \in 
\calC^\infty(\overline{\calO(T)}).
\end{equation*}
In the latter case, $w^*\alpha_j=-\sum n_k\alpha_k$, where the $n_k$ are again all
nonnegative. But the range of values of $w^*\al_j$ on $w^{-1}(\calO(T))$ matches that
of $\alpha_j$ on $\calO(T)$, i.e.\ $w^*\alpha_j\geq T_j$ here. In addition,
$\alpha_k\geq T_k$, on $\calO(T)$. These inequalities imply that for each $\ell$, 
\begin{equation*}
n_\ell \alpha_\ell = -\sum_{k\neq \ell} n_k\alpha_k- w^*\al_j \leq -\sum_{k\neq \ell}n_kT_k-T_j,
\end{equation*}
i.e.\ $n_\ell \al_\ell$ is bounded above on $w^{-1}(\calO(T))\cap \calO(T)$.
Hence either $n_\ell=0$, or else $\alpha_\ell$ is bounded above there. 
Writing $L=\{\ell:\ n_\ell\neq 0\}$,
\begin{equation*}
w^*e^{-\alpha_j}=\prod_{\ell\in L}(e^{\alpha_\ell})^{n_\ell}
=\prod_{\ell\in L}\tau_\ell^{-n_l},
\end{equation*}
which by the discussion above certainly extends smoothly to 
$\overline{w^{-1}(\calO(T))\cap \calO(T)}$.  

This proves that the transition maps are smooth, and hence that $\ofraka$ has the 
structure of a smooth manifold with corners. This completes the construction.

Following the arguments of the previous paragraphs, we see that 
this `bar compactification' construction commutes with taking products, i.e.\ if 
$\fraka = \fraka' \oplus \fraka''$, then
\begin{equation}
\ofraka = \overline{\fraka'} \times \overline{\fraka''}.
\label{eq:bcprod}
\end{equation}
Using this, we can directly adapt the construction to the reductive case, 
where the root system $\Lambda$ vanishes identically on the second factor,
once we have defined the appropriate compactification of an `unadorned' 
Euclidean space ${\mathfrak b}$, with trivial root system. In this case,
$\overline{\mathfrak b}$ is the `logarithmic blow-down' of the radial 
compactification $\widehat{\mathfrak b}$. Namely, it is the smooth 
manifold with boundary such that $\overline{\mathfrak b}_{\log} = 
\widehat{\mathfrak b}$; in other words, if $x$ is a smooth boundary 
defining function for $\widehat{\mathfrak b}$, then $\overline{\mathfrak b}$ 
is the same space as $\widehat{\mathfrak b}$, but with the smaller $\calC^\infty$ 
structure, where by definition $e^{-1/x}$ is a boundary defining function. With this
understanding, (\ref{eq:bcprod}) defines the bar compactification even in the 
reductive case.

Let us now examine the lift of $\Delta_\rad$ to $\ofraka$. It suffices for now 
to restrict to any $\overline{\calO(T)}$ where all $T_j > 0$ (to avoid the Weyl 
chamber walls). We can study the form of this operator near $\pa\, \ofraka$ 
by changing variables from $\{\al_1, \ldots, \al_n\}$ to $\{\tau_1, \ldots, \tau_n\}$. 
We have $\del_{\al_j} = -\tau_j \del_{\tau_j}$, and these latter vector fields 
generate $\calV_b(\ofraka)$, the space of smooth $b$ vector fields on $\ofraka$;
by definition $\calV_b$ consists of all smooth vector fields on $\ofraka$ which
are unconstrained in the interior but lie tangent to all boundaries.
Thus, all translation-invariant vector fields on $\fraka$ lift
to elements of $\calV_b(\ofraka)$, and indeed the latter is generated
by the lifts of these vector fields over $\Cinf(\ofraka)$. Hence,
all translation-invariant differential operators on $\fraka$
lift to elements of $\Diffb^*(\ofraka)$, the space of operators which
can be written locally as finite sums of elements of $\calV_b(\ofraka)$.

In particular, the
principal part $\Delta_\fraka$ is transformed to an elliptic, constant coefficient 
combination of these basic $b$ vector fields. In addition, $\coth\al-1$
is a $\Cinf$ function on $\fraka$ away from the Weyl chamber walls.
Indeed, $\coth \al - 1 = 2e^{-2\al}/(1-e^{-2\al})$, 
and so for $\al = \sum n_j \al_j 
\in \Lambda^+$, we have
\begin{equation*}
\coth \al - 1= \frac{\exp(-2\sum_{j=1}^n n_j \al_j)}{1-\exp(
-2\sum_{j=1}^n n_j \al_j)}
=\frac{\prod_{j=1}^n \tau_j^{2n_j}}{1-\prod_{j=1}^n \tau_j^{2n_j}},
\end{equation*}
which is certainly a $\Cinf$ function of the $\tau_j$ if $\tau_k<1$ for
all $k$. Since
\begin{equation*}
H_\al=\sum_{j=1}^n n_j \del_{\al_j} =\sum_{j=1}^n n_j (-\tau_j \del_{\tau_j})
\end{equation*}
is a translation-invariant vector field on $\fraka$,
we deduce that
away from the Weyl chamber walls, $\Delta_{\mathrm{\rad}}$ is indeed
an elliptic element of $\Diffb^2(\ofraka)$.

This may lead one to conclude that, except possibly having to
deal with some technicalities along the walls (which could be eliminated 
by working on the analogous compactification $\olM$ of $M$ which we define later), 
$\Diffb^*(\ofraka)$ is the appropriate setting to analyze $\Delta_{\rad}$.
However, this is not the case since the techniques of the so-called $b$-calculus 
on manifolds with corners only applies for operators which are asymptotically
of product type near the corners. This is unfortunately false for $\Delta_\rad$,
ultimately because the $\al_j$ are not orthogonal, but we now explain this more carefully.

The roots $\al_j$ are the linear coordinates for the dual basis $K_1, \ldots, 
K_n$ of $\fraka$ associated to $\Lambda^+_\ind$ (by $\al_i(K_j) = \delta_{ij}$ for all 
$i,j$). If $e_1, \ldots, e_n$ is any orthonormal basis for $\fraka$, then any vector $v 
\in \fraka$ can be expressed in terms of either basis:
\[
v = \sum_{j=1}^n y_j e_j = \sum_{\ell=1}^n x_\ell K_\ell. 
\]
Letting $\calK$ be the matrix with columns $K_1, \ldots, K_n$, then $y = \calK x$,
and so if $\calK^{-1} = \calH = (H_{rs})$, then we have
\[
\Delta_\fraka = \sum_{i=1}^n \frac{\del^2\ }{\del y_i^2} = \sum_{i,p,q=1}^n
\frac{\del x_p}{\del y_i}\frac{\del x_q}{\del y_i}\frac{\del^2\ }{\del x_p \del x_q}
=
\sum_{i,p,q=1}^n H_{pi}H_{qi}\frac{\del^2\ }{\del x_p \del x_q}.
\]

Next, associated to each $\al_j$ is the metrically dual vector $H_j$, i.e. 
$\al_j(w) = \langle H_j,w\rangle$ for all $w \in \fraka$. Then $\al_j(K_i) = 
\delta_{ij} = \langle H_j,K_i\rangle$, which means that the matrix $\calH = 
\calK^{-1}$ appearing above has columns equal to the vectors $H_1, \ldots, H_n$. 
We have thus shown that
\begin{equation}
\Delta_\fraka = \sum_{p,q=1}^n \gamma_{pq}\frac{\del^2\,}{\del x_p \del x_q},
\label{eq:fllap1}
\end{equation}
where $\Gamma = (\gamma_{pq}) = \calH \calH^t$. Finally, in terms of the coordinates
$\tau_j = e^{-\al_j}$, we have 
\begin{equation}
\Delta_\fraka = \sum_{p,q=1}^n \gamma_{pq}(\tau_p \del_{\tau_p})(\tau_q \del_{\tau_q}).
\label{eq:fllap2}
\end{equation}
However, the matrix $\Gamma$ is usually not diagonal, i.e.\ $\Delta_\fraka$ is not
`product-type'.

\subsection{The compactification $\tfraka$} 
We now describe the final, dominating, compactification $\tfraka$.
This is adapted from a compactification used in more general many-body settings,
as initially defined by the second author and employed in \cite{Vasy:Propagation-Many}. 
We first 
present this from the general point of view, not using the roots or the 
Weyl group action, but only the existence of a finite lattice $\calS$ of subspaces of 
the ambient space $\fraka = \RR^n$. This first construction of $\tfraka$ does not 
pass through $\ofraka$ as an intermediate space, but at the end of the section 
we discuss the relationship between the two spaces $\ofraka$ and $\tfraka$ and 
present a different construction of the latter space which does pass through the former.  

Let $\calS$ be the collection of all intersections of Weyl chamber walls 
$W_\al$ (as well as the `empty intersection' $\fraka$); this is a 
lattice, since it is closed under intersections and contains both $\{0\}$ 
and $\fraka$. We index this collection by a set $I$, so
$\calS = \{S_b: b \in I\}$; in particular, we suppose that
$\{0, *\} \subset I$, where $S_0 = \fraka$ and $S_* = \{0\}$. 
Finally, for any $S_b \in \calS$, write $S^b$ for the orthocomplement $S_b^\perp$. 

Now let us proceed with the construction. In the first step we pass to the 
radial (or geodesic) compactification $\hat\fraka$, which is obtained
by (hemispherical) stereographic projection, or alternatively, by compactifying 
each ray ${\mathfrak r} \cong [0,\infty)$ emanating from a fixed basepoint $o\in\fraka$ 
as a closed interval $[0,\infty]$. As described earlier, there is a natural topology 
and differential structure which makes $\hat\fraka$ into a smooth manifold with boundary.

Next, let $C_b$ be the boundary of the closure of $S_b$ in $\hat\fraka$; this
is a great sphere of dimension $\dim S_b - 1$. The collection of all such great 
spheres $\calC=\{C_b:\ b\in I\}$ is again a lattice. The singular and regular 
parts of $C_b$ are defined by 
\begin{equation*}
C_{b,\sing}=\bigcup\{C_c:\ C_c\subsetneq C_b\}, \qquad C_{b,\reg}=C_b\setminus C_{b,\sing},
\end{equation*}
and the singular and regular parts of $S_b$ are defined analogously.
The space $\tfraka$ is obtained by blowing up the collection $\calC$ inductively, 
in order of increasing dimension, as follows. $\calS$ is a union of subcollections
$\calS_j$, where $\dim S = j$ for any $S \in \calS_j$. We first blow up the
set of points $C_b$ corresponding to $S_b \in \calS_1$ to obtain a space 
$\hat{\fraka}^{(1)}$. Next, define the collection $\calC^{(1)}$ of submanifolds 
with boundary obtained by lifting the regular parts $C_{b,\reg}$ of each of the 
remaining sets $C_b$ and taking their closures in $\hat{\fraka}^{(1)}$. This 
is again a lattice, but the minimal dimension of its elements is now $1$, corresponding
to elements $S_b \in \calS_2$; furthermore, these $1$-dimensional submanifolds with
boundary are disjoint. We blow these up to form a space $\hat{\fraka}^{(2)}$. Continue
this process, obtaining a sequence of spaces $\hat{\fraka}^{(\ell)}$ and lattices
of submanifolds $\calC^{(\ell)}$ with components of dimension greater than or equal 
to $\ell$, and with all $\ell$-dimensional components disjoint submanifolds with
corners. We obtain after $n$ steps the space $\tfraka := \hat{\fraka}^{(n)}$. 
This compactification is a smooth manifold with corners, and is equipped with a 
smooth blow-down map $\beta:\tfraka \to \hat\fraka$. 

\begin{figure}[ht]
\begin{center}
\mbox{\epsfig{file=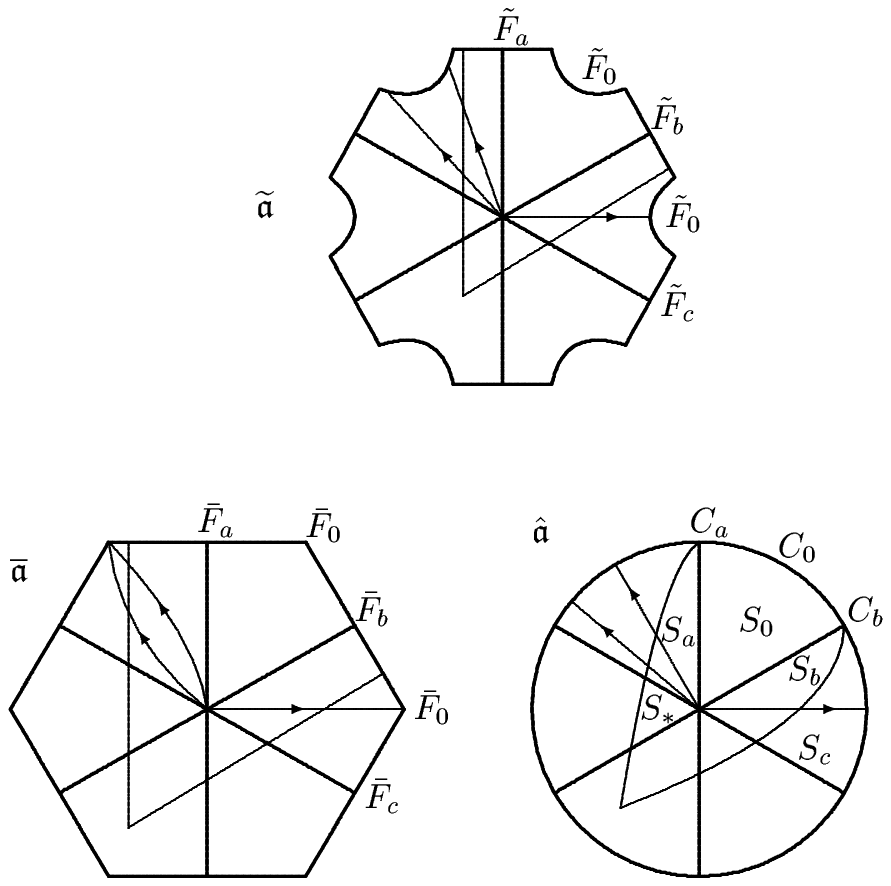}}
\end{center}
\caption{Representation of the compactifications $\ofraka$, $\hat\fraka$ and
$\tfraka$ for $M=\SL(3,\RR)/\SO(3,\RR)$.
The thick lines indicate the boundary
faces and the Weyl chamber walls. The thin lines without arrows
show the boundary of the closure of
$\calO(T)$, for $T_1<0$, $T_2<0$, in the various compactifications.
The thin lines with arrows are geodesic rays emanating from $0$; in particular
they bound {\em conic} regions. Geodesic rays in a single Weyl chamber
in $\ofraka$ hit the same point on $\pa\ofraka$, whereas in $\hat{\fraka}$, 
the boundary lines of $\calO(T)$ hit $C_a$ and $C_b$ for any $T$.}
\label{fig:acrp2}
\end{figure}

Notice that the indices $b \in I\setminus \{*\}$ are in bijective correspondence 
with the codimension one boundary faces of $\tfraka$, and also with the boundary 
faces of arbitrary codimension of $\ofraka$. Thus associated to any
$C_b$ is the (possibly disconnected) boundary hypersurface $\tilde{F}_b$
of $\tfraka$, and higher codimensional boundary face $\overline{F}_b$ of $\ofraka$. 
This suggests the alternate definition of $\tfraka$ as the logarithmic total 
boundary blow-up of $\ofraka$. More specifically, first replace each boundary 
defining function $\tau_j$ of $\ofraka$ by $\ot_j = -1/\log \tau_j$; then blow up 
the corners of $\ofraka$ inductively, in order of increasing dimension. This is 
essentially dual to the previous construction. In fact, the face $\tilde{F}_0$, 
corresponding to $S_0 = \fraka$ and $C_0 = S^{n-1}$, is the face obtained in this 
alternate definition by blowing up the highest codimension corners of $\ofraka$. 
Similarly, the faces $\tilde{F}_j$ created at the first stage in the first
definition of $\tfraka$ by blowing up the one dimensional elements $C_1$ 
correspond to the hypersurface faces of $\ofraka$.  All other faces of $\tfraka$ 
correspond to the various intermediate codimension corners in $\ofraka$. 
In any case, blowups of the boundary hypersurfaces of $\ofraka$ occur as boundary 
hypersurfaces of $\tfraka$, but that there are many other boundary hypersurfaces of 
this latter space, or in other words, $\tfraka$ distinguishes more directions
of approach to infinity. The replacement of each defining function by its 
logarithm here reflects the fact that in the ball model of hyperbolic
space, for example, the defining function $x$ is essentially $\exp(-\mbox{dist})$, 
while in the stereographic compactification of Euclidean space, the defining 
function $x$ is $1/\mbox{dist}$. We refer to \S 6 of \cite{Mazzeo-Vasy:Resolvents} 
for an extensive discussion of the role of smooth defining functions in 
compactification theory. 

The behaviour of this `tilde compactification' with respect to taking products
is a bit more complicated than for the bar compactification. First of all, 
if the root system of $\fraka$ is trivial, i.e.\ $\fraka$ is an unadorned 
Euclidean space, then $\tfraka = \hat\fraka = \ofraka_{\log}$. Secondly, if $\fraka = 
\fraka' \oplus \fraka''$, then $\tfraka$ is obtained by blowing up the 
closed ball $\hat\fraka$ along the collection of boundary submanifolds 
$\calC = \{C_a\} = \{\pa S_a\}$, where each $S_a$ is of the form $S_b' \times S_c''$ 
(including, of course, the cases $S_b'=\{0\}$ or $S_c''=\{0\}$). Hence $C_a$ is either the 
simplicial join $C_b' \# C_c''$ (regarded as a smooth great sphere in $\pa \hat\fraka$)
or else $C_b' \times \{0\}$ or $\{0\} \times C_c''$;
in particular, if all roots vanish on $\fraka''$, then each $C_a$ equals either
$C_b' \# \pa \hat\fraka''$ or $C_b' \times \{0\}$. Of course, we can also obtain 
$\tfraka$ as the total boundary blowup of $\ofraka$, i.e.\ as 
\begin{equation}
\tfraka = \left[ \left(\ofraka\right)_{\log}; \overline{\calF}\right]
= \left[ \left(\ofraka' \times \ofraka'' \right)_{\log}; \overline{\calF}\right]
= \left[ \left(\ofraka'\right)_{\log} \times \left(\ofraka'' \right)_{\log}; 
\overline{\calF}\right],
\label{eq:tprod1}
\end{equation}
where $\overline{\calF}$ is the collection of boundary faces of all codimension
in $\ofraka$. If all roots vanish on $\fraka''$, then
\begin{equation}
\tfraka = \left[ \left(\ofraka'\right)_{\log} \times \widehat{\fraka''};
\left(\overline{\calF'} \times \widehat{\fraka''}\right) \cup
\left(\left(\ofraka'\right)_{\log} \times \pa\,\widehat{\fraka''}\right)\right].
\label{eq:tprod2}
\end{equation}

\subsection{Compactifications of the full symmetric space} 
Before continuing with the more detailed description of $\Delta_\rad$ on $\tfraka$, we 
follow the train of thought from the past two subsections and define the 
compactifications $\Mb$ and $\Mt$ of the full symmetric space $M$, corresponding to 
$\ofraka$ and $\tfraka$, respectively. Their role in this paper is only minor since 
our emphasis is on the radial Laplacian. Nevertheless, many properties of the
operator $\Delta_\rad$, which has nonsmooth coefficients on $\fraka$, are proved
by appealing to its lift to $M$, which is just the operator $\Delta$, and 
which does have smooth coefficients; we also consider lifts of $\Delta_\rad$ to 
certain spaces intermediate between between the various compactifications of $M$ 
and $\fraka$.

As we have seen in \S 2.1, the Cartan decomposition $G = K \overline{\calC^+} K$
states that any $g \in G$ has a decomposition $k_1\cdot a \cdot k_2$, where 
$k_1,k_2 \in K$ and $a = \exp(H)$, $H=H(g) \in \overline{C^+}$, and with this 
normalization, $a$ is {\em unique}. Moreover, if $p\in M=G/K$ has
$H(p)\in\calC^+$ then $K^p$, the subgroup of $K$ that fixes $p$, is discrete;
the set of such $p$ is open and dense in $M$ and is diffeomorphic 
to $(K/K^{p_0})\times\calC^+$ (for any $p_0 \in \calC^+$).

As discussed in \S 2.6,
each (open) face $S_b^+$ of the closed positive Weyl chamber
$\overline{\calC^+}$ 
in $\fraka$ is an open set in a unique $S_b$, $b \in I$, and we
index the set of all such faces $S_b^+$ by a subset $I^+ \subset I$.

If $p\in \exp(S_{b,\reg}\cap \overline{C^+})$, $b\in I^+$,
let $\Lambda_b$ be the
set of roots vanishing at $p$. Since $S_b \subset \fraka \subset \frakg_0$,
there is an orthogonal splitting $\frakg_0 = S_b \oplus \frakg_0^b$, and we
then define
\begin{equation*}
\frakg^b = \frakg_0^b \oplus \sum_{\alpha\in\Lambda_b}\frakg_\alpha,\qquad
\mbox{and}\quad \frakp^b=\frakp\cap\frakg^b,
\end{equation*}
cf.\ \cite[Section~2.20]{Eberlein:Geometry}.
This is the Lie algebra of a Lie subgroup $G^b \subset G$,
which contains
the isotropy group of $p$ in $K$. Denoting this latter group by  
$K^b$, and its Lie algebra by $\frakk^b$, then $\frakg^b = 
\frakk^b \oplus \frakp^b$. There is a corresponding symmetric space 
$\Sigma^b = G^b/K^b$, which is identified with $\exp(\frakp^b)$.
Now, the image $N$ of a neighbourhood of $(S_{b,\reg}\cap \overline{C^+})
\times\{0\}$ in $(S_{b,\reg}\cap \overline{C^+})\times \frakp^b$ under $\exp$
is a submanifold of $M$, with $p$ lying on
it, and the $K$-action is transversal to $N$ at $p$.
Thus, a neighbourhood of the $K$-orbit of $p$ is diffeomorphic to the
$K$-orbit of the $K^b$-class of $(H(p),e,o)$,
where $e$ is the identity element in $K$ and 
$o$ the identity coset in $\Sigma^b$, in
\begin{equation*}
S_b\times (K \times \Sigma^b)/K^b \, , \qquad \mbox{where}\quad
k_1 \cdot (k,\sigma) = (kk_1^{-1}, \, k_1 \cdot\sigma)\quad \mbox{for any} \ 
k_1 \in K^b.
\end{equation*}
We can let $p$ vary in $\exp(S_{b,\reg}\cap \overline{C^+})$, and deduce
that a neighborhood of the $K$-orbit of $\exp(S_{b,\reg}\cap \overline{C^+})$
is diffeomorphic to the
$K$-orbit of the $K^b$-class of $(S_{b,\reg}\cap \overline{C^+})\times \{e\}
\times\{o\}$.
Reinterpreted, this says that the $K$-orbit of a neighbourhood of
$\exp(S_{b,\reg}\cap \overline{C^+})$ in $M$
is a $\Cinf$ bundle over $K/K^b \times \exp(S_{b,\reg}\cap \overline{C^+})$ 
with fibre (a neighbourhood of the origin in) $\Sigma^b$.

In fact, this argument shows more. Consider the action of $\Real^+$ by
dilations on $\frakp$: $\Real^+\times\frakp\ni (t,z)\mapsto tz\in\frakp$. 
A set is called conic if it is invariant under
the $\Real^+$-action. As remarked before, this $\RR^+$-action on $\frakp$ is
identified with dilations along the geodesic rays through $o$ via the
exponential map. Now, $k\cdot\exp(tX)=\exp(\Ad(k)(tX))=\exp(t\Ad(k) X)$
for $k\in K$, $X\in\frakp$, $t\in\RR^+$. Thus, under the identification
of a neighbourhood of $p$ as above with a neighbourhood of $(e,o,0) \in 
(K/K^b) \times\Sigma^b\times S_b$, 
the $\RR^+$-action is $(t,kK^b,q,x)\mapsto(t,kK^b,tq,tx)$, at first
for $t$ near $1$. Thus, we can extend the identification to a conic
neighbourhood of the $\RR^+$-orbit of $p$ via the dilation.
Letting $p$ vary in a bounded set, we deduce that there is
a conic neighbourhood $U_b$ of $S_{b,\reg}\cap \overline{C^+}$ in $\fraka$
such that $K\cdot \exp(U_b)$ can be identified with
a $\Cinf$ bundle over $K/K^b \times \exp(S_{b,\reg}\cap \overline{C^+})$ 
with fibre (a neighbourhood of the origin in) $\Sigma^b$.
We let $\Phi_b$ be this identification.

If $p\in \exp(S_{c,\reg}\cap\overline{C^+}) \cap \exp(U_b)$, then 
$S_b\subset S_c$ and $p\in \exp(U_c)$ as well, so there are two 
identifications of a conic neighbourhood of $p$: one
as a subset of $(K/K^b)\times\Sigma^b\times S_b$, and the other as
a subset of $(K/K^c)\times\Sigma^c\times S_c$. Since $K^c\subset K^b$,
we have $K/K^b\subset K/K^c$ and $\Sigma^c\subset\Sigma^b$. The map 
between these two identifications is
thus a diffeomorphism, and it commutes with the $\Real^+$-action.

We can now define $\Mb$; this is called the dual cell
compactification of $M$, see
\cite[Section~3.40]{Guivarch-Ji-Taylor:Compactifications}, where it is
defined as a topological space with a $G$-action.
Our construction proceeds by partially compactifying
part of the regions described in the preceeding paragraphs.
Thus, we fix a $K^b$-invariant bounded neighbourhood
$\calO_b$ of $o$ in each symmetric
space $\Sigma^b$; this has a $W^b$-invariant bounded
intersection $O_b$ with $S^b$.
Let $V_b$ be an open subset of $S_{b,\reg}$ such that $S_{b,\reg}
\setminus V_b$ is bounded and $V_b\times O_b\subset U_b$. Such a subset
exists since $U_b$ is a conic neighbourhood of $S_{b,\reg}\cap\overline{C^+}$.
Then, by the preceeding discussion,
$K\cdot\exp(V_b\times O_b)$ is a $\Cinf$ bundle
over $(K/K^b)\times V_b$ with fiber $\calO_b$. We partially compactify the base
of this bundle as $(K/K_b)\times \overline{V_b}$, where $\overline{V_b}$
is the closure of $V_b$ in $\overline{S_{b,\reg}}$, the regular part
of the bar-compactification of $S_b$.

If now $c$ is such that $S_b\subset S_c$, then we have seen that
on $K\cdot\exp((V_b\times O_b)\cap(V_c\cap O_c))$ the transition
maps between the identifications of the respective bundles is a
diffeomorphism.
It is now immediate that the same is true in these partial compactifications
since this amounts to showing that the identification map on the subset
$(V_b\times O_b)\cap(V_c\times O_c)$ of $\fraka$ extends to be
smooth on $(\overline{V_b}\times O_b)\cap({\overline{V_c}}\times O_c)$,
which is immediate from the definition of $\ofraka$.

\begin{figure}[ht]
\begin{center}
\mbox{\epsfig{file=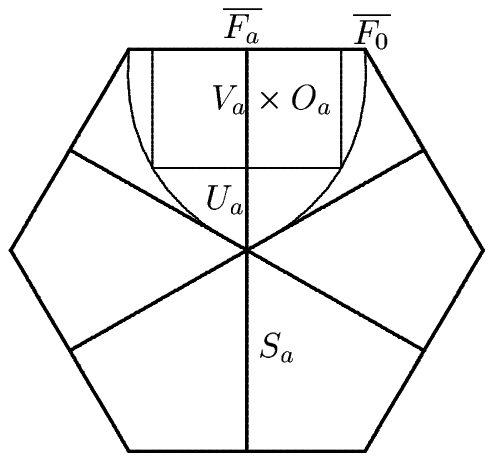}}
\end{center}
\caption{Subsets of $\ofraka$ used in the construction of $\Mb$
for $M=\SL(3,\RR)/\SO(3,\RR)$.
The thick lines indicate the boundary
faces and the Weyl chamber walls. The rectangular thin lines
show the boundary of $V_a\times O_a$. The curved ones indicate
the boundary of $U_a$; they are in particular geodesic rays from $o$.
The corresponding subsets for $b=0$ are $U_0=\calC^+$, the
positive Weyl chamber, $O_0=\{o\}$ and $V_0=\calC^+$. Thus, the $0$-chart
covers a neighbourhood of the corner, $\overline{F_0}$.}
\label{fig:acrp4}
\end{figure}

We can thus define $\Mb$ as the disjoint union of the $\calO_b$-bundles
over $(K/K_b)\times \overline{V_b}$, $b\in I^+$,
modulo the equivalence relation
corresponding to this identification.
Then $\Mb$ is a manifold with corners -- the corners arise from the
$\overline{V_b}$, i.e.\ from the compactification of the flat.

Even though we have remained in a bounded neighbourhood of $o$ in each
symmetric space $\Sigma^b$ to avoid a recursive definition of
the compactifications, it is now immediate that the boundary faces
$\overline{\calF_b}$,
$b\in I^+$,
of $\Mb$ are $\Cinf$ bundles over $K/K_b$ with fiber $\overline{\Sigma^b}$
(the bar-compactification of $\Sigma^b$). Indeed, this simply relies
on considering the closure of the conic set
$K\cdot\exp(U_b)$ in $\Mb$. Note, however,
that this closure does {\em not} include a neighbourhood of
$\overline{\calF_b}$.
Indeed, the issue is that the closure of $U_b$ in $\ofraka$ does not
include a neighbourhood of the face $\overline{F_b}$, though it {\em does}
contain a neighbourhood of the open face $F_b$.

This procedure may be modified easily for the construction of $\Mt$.
Indeed, in each step we simply replace
$\overline{V_b}$ by $\widetilde{V_b}$,
the closure of $V_b$ in $\widetilde{S_{b,\reg}}$, the regular part
of the bar-compactification of $S_b$.
By the naturality of all the steps, it is clear that we could also define 
$\Mt$ as the logarithmic total boundary blow-up of $\Mb$.

We recall that as a topological space, it is described in
\cite{Guivarch-Ji-Taylor:Compactifications} as
the smallest compactification that dominates both $\Mb$ and the
geodesic (or conic) compactification $\hat M$. Note that the latter does not
have a natural smooth structure: if it is defined by compactifying $\frakp$
radially and using the exponential map, the smooth structure depends on the
choice of the base point $o$. It is shown in
\cite[Theorem~8.21]{Guivarch-Ji-Taylor:Compactifications}
that, as a topological space, $\Mt$ is the Martin compactification of $M$.

\begin{rem}\label{rem:Mb-Mt}
Although we have defined $\Mb$ and $\Mt$, we never actually use them
in this paper. Rather, since we are working with $K$-invariant functions
and operators, the only reason to leave $\fraka$ (or $\ofraka$ and
$\tfraka$) is to make the differential operators have smooth
coefficients. For this purpose, the $K/K_b$ factor can be ignored, and
we may work instead on $V_b\times\calO_b$, etc, which is exactly
what we do in \S~4. However, it is
nice to know that there is a compactification of $M$ in the background,
rather than just an ad hoc collection of product spaces!
\end{rem}

\subsection{The lift of $\Delta_\fraka$ to $\tfraka$} In the remaining subsections of 
\S 2 we shall be examining the structure of $\Delta_\rad$ on $\tfraka$ in some detail, 
focusing specifically on its behaviour at and near the boundary. This involves several 
steps.  In this subsection we study the lift of the flat Laplacian $\Delta_\fraka$, 
and vindicate our earlier claim that this operator attains a product-type structure
near the corners of $\tfraka$. The results of this section are not used
elsewhere in the paper.

Recall the expression (\ref{eq:fllap2}), which exhibits $\Delta_\rad$ as an elliptic 
$b$-operator on $\ofraka$. We now introduce a singular change of variables on $\ofraka$. 
Using multi-index notation, set 
\[
\sigma = \tau^\theta,\qquad \mbox{i.e.} \quad \sigma_i = \tau_1^{\theta_{i1}}\ldots 
\tau_n^{\theta_{in}},
\]
where $\Theta = (\theta_{ij})$ is some $n$-by-$n$ matrix to be determined. We calculate
\[
\tau_s \del_{\tau_s} =  \sum_{r=1}^n \theta_{rs}\sigma_r \del_{\sigma_r},
\]
and so
\[
\Delta_\fraka = \sum \gamma_{pq}\theta_{ip}\theta_{jq}(\sigma_i \del_{\sigma_i})
(\sigma_j \del_{\sigma_j}) =
\sum \nu_{ij}(\sigma_i \del_{\sigma_i})(\sigma_j \del_{\sigma_j}),
\]
where $N = (\nu_{ij}) = \Theta\Gamma \Theta^t$. We wish to choose $\Theta$ so that $N$ 
is diagonal. We intend to study $\Delta_\fraka$ (and $\Delta_\rad$) near the closure of 
some face $F$, which we label for simplicity as $\tau_1 = 0$; the ordering of the other 
faces is then arbitrary. Relative to this ordering, since $\Gamma$ is positive definite, 
there is a factorization $\Gamma = LDU$, where $L$ and $U$ are lower and upper 
triangular, respectively, and $D$ is diagonal. Since this
factorization is unique, and $\Gamma = \Gamma^t$, we must have $U = L^t$. Hence if we 
define $\Theta = L^{-1}$, which is also lower triangular, then $L^{-1}\Gamma (L^{-1})^t = N$ 
is the diagonal matrix $D$ appearing in the decomposition, as desired. Somewhat more explicitly, 
this coordinate change has the form
\[
\sigma_1 = \tau_1,\ \sigma_2 = \tau_1^{\theta_{21}} \tau_2, \ldots 
\sigma_n = \tau_1^{\theta_{n1}} \cdots \tau_{n-1}^{\theta_{n\,n-1}} \tau_n.
\]
We have now shown that $\Delta_\fraka$ may be transformed to diagonal form near any corner 
of $\ofraka$, but at the expense of using a singular coordinate change. 

The other key step is to show that this singular coordinate change lifts to a smooth 
(local) diffeomorphism of $\tfraka$. Recall that this latter space is obtained 
by first introducing the logarithmic change of variables $\ot_i = -1/\log \tau_i$, and 
then blowing up the corners in order of increasing dimension. Defining $\osig_i = 
-1/\log \sigma_i$, then 
\[
\frac{1}{\osig_1} = \frac{1}{\ot_1}, \ \ \ldots\ \   , \ \ 
\frac{1}{\osig_j} = \frac{\theta_{j1}}{\ot_1} + \ldots + \frac{\theta_{j\, j-1}}{\ot_{j-1}} +
\frac{1}{\ot_j}, \ \ \ldots
\]
These formul\ae\ represent the lift of this map acting between $(\ofraka)_{\log}$, but it
is still not smooth. The passage to the total boundary blowup fixes this: 
to this end, first note that each $\osig_j$ is homogeneous of degree $1$ in the $\ot_i$, 
and so if we introduce polar coordinates $\ot = r \omega$, $\osig = r' \phi$ near 
$\ot = \osig = 0$, then we can identify the radial variables, $r=r'$. For simplicity, 
we examine this near the codimension $2$ corners of the blowup, i.e.\ near where exactly 
one of the $\omega_i$ vanish, and away from the higher codimension corners where two or more 
of these angular variables equal zero. Thus suppose we are working near $\omega_j = 0$. 
For every $k$ we have 
\begin{equation}
\frac{1}{\phi_k} = \frac{\theta_{k1}}{\omega_1} + \ldots + 
\frac{\theta_{k\, k-1}}{\omega_{k-1}} + \frac{1}{\omega_k}.
\label{eq:cch}
\end{equation}
Thus if $k < j$ then $\phi_k$ is obviously a smooth function of $\omega$ since all terms here
are nonvanishing (note that the whole right hand side cannot vanish, since otherwise we 
would reach the incorrect conclusion that $\osig_k$ itself would be undefined). Next, if 
$k = j$, then we can rewrite (\ref{eq:cch}) as 
\[
\phi_j = \frac{\omega_j}{
\theta_{j1}\frac{\omega_j}{\omega_1} + \ldots + \theta_{j\, j-1}\frac{\omega_{j}}{\omega_{j-1}} 
+ 1}, 
\]
which again is certainly smooth. Finally, if $k > j$, then
\[
\phi_k = \frac{\omega_k}{
\theta_{k1}\frac{\omega_j}{\omega_1} + \ldots + \theta_{kj} + \ldots 
+ \frac{\omega_j}{\omega_k}};
\]
if $\theta_{kj} \neq 0$, then this is smooth near $\omega_k = 0$, while if 
$\theta_{jk} = 0$, then $\phi_k$ is independent of $\omega_j$, hence 
again is smooth. The argument near the higher codimension corners is similar.

\subsection{Subsystems} 
We now consider the restrictions of $\Delta_\rad$ to the codimension one 
boundary faces of $\tfraka$; our goal is to show that each such restriction
is essentially the radial Laplacian on some lower rank symmetric space. 
To this end, we examine the geometry of $\pa\tfraka$ more closely. 

\subsubsection{Geometric and algebraic subsystems}
Any point $p\in\pa\hat\fraka$ belongs to a unique $C_{b,\reg}$ for some $b \in I$.
Note that $C_c \cap C_{b,\reg} \neq \emptyset$ only when $C_c\supset C_b$, 
or equivalently when $S_c\supset S_b$. Thus, in particular, for any root $\al$, the wall 
$W_\alpha$ equals $S_c$ for some $c \in I$, and the corresponding 
$C_c$ intersects $C_{b,\reg}$ only when $W_\alpha \supset S_b$. Thus $p$ has a neighbourhood 
$\calU$ in $\hat\fraka$ such that $\calU\cap W_\alpha \neq \emptyset$ only when $S_b\subset 
W_\alpha$. 

Next, the boundary hypersurfaces $F$ of $\tfraka$ are in one-to-one
correspondence with the indices $b\in I\setminus\{*\}$, where $F_b$ is the 
front face created by blowing up $C_{b,\reg}$. The interior of each $F_b$ 
has a (trivial) fibration induced by the blow-down map $\beta$, with base $C_{b,\reg}$ 
and fibre the orthocomplement $S^b$. We remark that
this extends to a fibration of the
closed face $F_b$, with fibre $\widetilde{S^b}$, the compactification of
$S^b$ obtained analogously to $\tfraka$
by regarding $S^b$ as a flat in the lower
rank symmetric space $\Sigma^b$,
and base the closure of the lift of $C_{b,\reg}$ in the partially
blown up space $\hat\fraka^{(\ell)}$, $\ell = \dim C_b$.
The base can also be identified with the lift of
$C_b$ to $\widetilde{S_b}=[\widehat{S_b};\{C_c:\ C_c\subsetneq C_b\}]$.
Indeed, this is description is {\em identical} to the geometry of
compactifications in N-body scattering; see
\cite[pp.~339-340]{Vasy:Propagation-Many} for a very detailed discussion of
the latter.

Translating by an element of the Weyl group, we can suppose that $p \in \overline{C^+}$.  
Let us then say that a root $\al$ is positive, negative, or zero at $p$ if $\al$ 
has this property on the ray in $\fraka$ corresponding to $p$. In particular, 
$\al$ vanishes at $p$ (and at every other $q\in C_{b,\reg}$ as well) if and only if 
$W_\al \supset S_b$. 

Let $\Lambda_b$ denote the subset of all roots $\al$ which vanish on $S_b$. We have
identifications
\[
\{\gamma \in \fraka^*: \gamma = 0\ \mbox{on}\ S_b\} \cong
(\fraka/S_b)^* \cong (S^b)^*;
\]
the first of these is tautological, while the second uses
the metric, but both are isometries. Hence we can also regard $\Lambda_b \subset (S^b)^*$,
with the same inner product relations as in $\fraka^*$, and clearly this is
a spanning set of covectors. In addition, $\al \in \Lambda_b$ 
if and only if $W_\al^\perp \subset S^b$, or equivalently $H_\al \in S^b$. 
It is now easy to check that $\Lambda_b$ satisfies all the axioms of a reduced 
root system on $\Span(\Lambda_b) \subset (S^b)^*$, cf.\ 
\cite[Section 9.2]{Humphreys:Lie}. We define $\Lambda_b^+=\Lambda_b\cap \Lambda^+$. 

In conclusion, we have shown that for each $b \in I\setminus \{*\}$,
$\fraka = S_b \oplus S^b$, where the latter summand is the Cartan subspace
for some symmetric space of rank less than $n$; furthermore, the face $F_b$ 
is the product of the base space, which is a compactification of $C_{b,\reg}$, 
and the radial compactification of the vector space $S^b$.
There is a more familiar geometric version of this statement. Fix $p \in 
C_{b,\reg}$ and let $\gamma$ be the geodesic in $M$ which is the exponential of 
the ray corresponding to $p$. We say that another geodesic $\gamma'$ is parallel 
to $\gamma$ if the two geodesics stay a bounded distance from one another in both 
directions. Following \cite{Eberlein:Geometry}, we define $F(\gamma)$ to be the union of all 
geodesics parallel to $\gamma$. This is a totally geodesic submanifold in $M$, and it
always admits a Riemannian product decomposition $\RR^k \times F_s(\gamma)$, where 
the second factor is a symmetric space of rank strictly less than $n$. The
correspondence is that the tangent space to these two factors are just
$S_b$ and $S^b$, respectively.

As noted earlier, the (interiors of the) faces $F_b$ which correspond to $1$-dimensional collision 
planes $S_b$ already appear as boundary hypersurfaces in the simpler 
compactification $\ofraka$. 

Even if $M$ itself is an irreducible symmetric space, the symmetric spaces 
$F_s(\gamma)$ which appear in these subsystems may well be reducible. On
the algebraic level, this occurs if there is an orthogonal decomposition 
$S^b=\oplus (S^b)_j$ so that each element of $\Lambda_b$ lie in one of the 
summands. An orthogonal partition of roots is the same as an orthogonal
partition of simple roots (see \cite[Section~10.4]{Humphreys:Lie}), and this
corresponds to the Dynkin diagram decomposing as a disjoint union. This phenomenon
occurs already in our standard examples $\SL(n+1)/\SO(n+1)$. In fact, to
every possible partition $m_1+\ldots+m_k = \ell \leq n$ one associates 
the subsystem 
\[
\RR^{n-\ell} \times \prod_{j=1}^k\SL(m_j+1)/\SO(m_j+1).
\] 
Thus, for example, the subsystems of $\SL(3)/\SO(3)$ are $\RR \times \HH^2= 
\RR \times \SL(2)/\SO(2)$, while the two different rank $2$ models $\RR \times 
\SL(3)/\SO(3)$ and $\RR \times \HH^2\times\HH^2$, and also the rank $1$ model 
$\RR^2 \times \HH^2$, comprise the subsystems of $\SL(4)/\SO(4)$.

\subsubsection{Analytic subsystems}
We now discuss the subsystem Hamiltonians, and the behaviour of
$\Delta_\rad$ near the faces of $\tfraka$. Set
\begin{equation}
\rho_b=\frac{1}{2}\sum_{\alpha\in\Lambda_b^+}m_\al\, \alpha \quad 
\left(\mbox{hence}\quad H_{\rho_b}\in S^b \right).
\label{eq:def-rho_b}
\end{equation}
The lifts of the roots $\al \in \Lambda^+ \setminus \Lambda^+_b$ to $\tfraka$
tend to $+\infty$ everywhere on the closed face $F_b$, so that the 
corresponding terms $(\coth \al -1)H_\al$ in $\Delta_\rad$ decay rapidly there 
and thus are negligible on that face. More precisely, we have the following
result.

\begin{lemma}\label{lemma:E_b-coeffs}
Let $Z_\alpha$ be the
closure of $\alpha^{-1}((-\infty,0])$ in $\hat\fraka$.
Then
\begin{equation*}
\coth \al -1\in\Cinf(\fraka\setminus
\alpha^{-1}(-\infty,0])
\end{equation*}
extends to an element of
$\Cinf(\hat\fraka\setminus Z_\alpha)$
that vanishes to infinite order at $\pa\hat\fraka\setminus
\pa Z_\alpha$.
Thus, if $\chi\in\Cinf(\hat\fraka)$ with $\supp\chi\cap Z_\alpha
=\emptyset$, then $\chi(\coth\al-1)\in\dCinf(\hat\fraka)$, i.e.\ it
vanishes to infinite order at $\pa\hat\fraka$.
\end{lemma}

\begin{proof}
The function $x\mapsto\alpha(x)/|x|$, $x\in\fraka\setminus\{0\}$,
is homogeneous degree zero,
so it extends to a smooth function on $\hat\fraka\setminus\{0\}$,
and its restriction to $\pa\hat\fraka\setminus
\pa Z_\alpha$ is strictly positive. It is immediate that
$e^{-\alpha(x)}=\exp\left(-\frac{\alpha(x)}{|x|}\,|x|\right)$ is smooth
and rapidly decreasing in $\hat\fraka\setminus Z_\alpha$, hence
the statements for $\coth\al-1=\frac{2e^{-2\al}}{1-e^{-2\al}}$ also follow.
\end{proof}

Note that if $\alpha\in\Lambda^+\setminus\Lambda^+_b$, then in particular
$C_{b,\reg}\subset\hat\fraka\setminus Z_\alpha$, so $\coth \al -1$
is Schwartz in a neighbourhood of $C_{b,\reg}$ in $\hat\fraka$.
In other words, there is a conic neighbourhood of $S_{b,\reg}$ in $\fraka$
on which $\coth\alpha-1$ is Schwartz.

We now return to $\Delta_\rad$. After subtracting
\begin{equation*}
E_b=\sum_{\alpha\in\Lambda^+\setminus\Lambda_b^+} (\coth\alpha-1)H_\alpha.
\end{equation*}
the remaining terms
\begin{equation}
L_b= \Delta_{S_b}+2(H_\rho-H_{\rho_b})+\Delta_{S^b}+2H_{\rho_b}
+\sum_{\alpha\in\Lambda^+_b}m_\alpha(\coth\alpha-1) H_\alpha.
\label{eq:decomp1}
\end{equation}

\begin{prop}\label{prop:product-model}
For each $b \in I \setminus \{*\}$ there is a decomposition
\[
L_b = T_b + \Delta_{b,\rad},
\]
where the first term is a constant coefficient elliptic operator on $S_b$ and 
the second is the radial Laplacian for the noncompact symmetric space
$\Sigma^b$, which has rank strictly less than $n$.
\label{pr:decomp-L_b}
\end{prop}

\proof 
The first summand, $T_b$, is the sum of the first two terms in
(\ref{eq:decomp1}), and $\Delta_{b,\rad}$ is the sum of the remaining three.
Since $\Lambda_b$ is a root system on $S^b$, it is clear that
\begin{equation}
\Delta_{\rad,b} := 
\Delta_{S^b}+2H_{\rho_b} +\sum_{\alpha\in\Lambda^+_b}m_\alpha(\coth\alpha-1) H_\alpha
\label{eq:decomp2}
\end{equation}
is indeed the radial part of the Laplacian on a symmetric space of lower rank. 
Thus it remains only to prove that the vector appearing as the first order term in $T_b$,
\begin{equation}
H_\rho-H_{\rho_b} = \frac{1}{2}\sum_{\alpha\in \Lambda^+ \setminus \Lambda_b^+}m_\al 
H_\al,
\label{eq:diffvec}
\end{equation}
is an element of $S_b$, as claimed. To prove this, note first that if $\beta$ is 
a simple root, with corresponding Weyl group element $w_\beta$ (the reflection 
across $W_\beta$) and $\al$ is a positive root which is linearly independent from 
$\beta$, then $w_\beta^*(\al)$ is again a positive root; for, $\al$ is nonnegative 
and not identically vanishing on $W_\beta \cap \overline{C^+}$, and $w_\beta$ 
fixes $W_\beta$ pointwise, hence $w_\beta^* \al$ is also nonnegative and not 
identically vanishing on this same set, hence must be positive on $C^+$, which 
is a characterization of positive roots. Next, clearly $H_{w_\beta^* \al} = 
w_\beta (H_\al)$ and so
\[
H_\al + H_{w_\beta^*(\al)} \in W_\beta.
\]
In addition, $m_{{w_\beta}^*\al} = m_\al$. Now let $\{\al_j: j \in J_b\}$ be 
an enumeration of the simple roots in $\Lambda^+_b$, and write $w_j = w_{\al_j}$. 
Then $w_j^*$ preserves the subsets $\Lambda_b$, hence also 
$\Lambda \setminus \Lambda_b$ and $\Lambda^+ \setminus \Lambda_b^+$ because $\al_j$
is linearly independent from any of the elements in these last two sets. Therefore 
\eqref{eq:diffvec}
is a sum over $w_j$ orbits, where each orbit consists of one or 
two elements: if it consists of just one element $\al$, then $H_\al \in W_{\al_j}$, 
and if it consists of two elements $\al$ and $\al'=w_j^* \al$, then 
$m_\al H_\al + m_{\al'} H_{\al'}$ also lies in $H_{\al_j}$. Hence
\eqref{eq:diffvec}
also lies in $W_{\al_j}$. This is true for every $j \in J_b$, and the claim follows. 
\endproof

In summary, we have made precise that $\Delta_\rad$ is locally
-- in a neighbourhood of the lift of $C_{b,\reg}$ to $\tfraka$ -- the
sum of a product model, $L_b$, and an error term $E_b$.

We remark that such a neighbourhood is diffeomorphic to an open subset
in the tilde-compactification of $\fraka$ with collision planes
given by $S_b\times (S_c\cap S^b)$
and $\{0\}$ as $S_c$ runs over all collision planes
satisfying $S_c\supset S_b$. In particular, if one studies the asymptotics
of the Green function, one can paste the asymptotics of the local model
operator Green functions directly from the
model space to $\tfraka$.

\section{Invariant smooth functions and \\ localization on the compactified
spaces}

\subsection{Invariant smooth functions}
As already discussed in \S 2.1, every $g \in G$ decomposes into a product
$g=k_1 a k_2$, where $k_1,k_2 \in K$ and $a \in \calA$; the middle factor
is determined up to translation by an element of $W$, and in particular is
unique if we require it to lie in $\overline{\calA^+}$. This defines a map
$\pi:M\to \overline{\calA^+}$. If $h$ is any (e.g.\ measurable) function
on $\overline{\fraka^+}$, or equivalently, a $W$-invariant function on $\fraka$,
then its pullback $\pi^* h$ is a $K$-invariant function on $G/K=M$. (As usual, 
we are identifying $\calA$ with $\fraka$.) Conversely, $K$-invariant functions 
on $M$ restrict to $W$-invariant functions on $\fraka$, and therefore $\pi^*$
induces an equivalence between these spaces.

It will be important for us to know whether $\pi^*$ yields an equivalence 
between functions with higher regularity.
Thus, for example, it is clear that $\pi^*$ induces an isomorphism
between continuous $W$- and $K$- invariant functions, and also between
$L^2_\loc$ invariant functions, though here we must use the degenerate
measure on $\fraka$ induced by pushforward by $\pi_*$ of a smooth 
invariant smooth measure on $M$. Somewhat more generally, $\pi$ is a 
Riemannian submersion since the $K$-orbits are orthogonal to $\calA$ and the metric
is invariant on both fibre and base. Hence it is distance-decreasing, 
i.e.\ $d(\pi(x),\pi(y))\leq d(x,y)$ for any $x,y\in M$; therefore $\pi$ is Lipschitz,
and $\pi^*$ gives an isomorphism between invariant functions which
are locally Lipschitz. The following result, however, is less obvious.

\begin{prop}\label{prop:A-M-W-K}
The map $\pi^*:\Cinf(\fraka)^W\to\Cinf(M)^K$ is an isomorphism.
\end{prop}

\begin{proof}
The easy direction is that the restriction of any $f \in \Cinf(M)^K$
to $\calA$ is in $\Cinf(\fraka)^W$. In fact, the inclusion map $\iota:
\calA \hookrightarrow M$ is smooth, so if $f \in \Cinf(M)$ then 
$\iota^*(f) \in \Cinf(\fraka)$. Moreover, since $W$ is the quotient
of the normalizer in $K$ of $\calA$ by its centralizer, $\iota$ commutes 
with the action of $W$, and so $\iota^*:\Cinf(M)^K\to\Cinf(A)^W$.

To prove the converse, we use induction on the rank $n$. Suppose the result 
has been proved for all symmetric spaces of rank strictly less than $n$. 
Fix $p\in \overline{\calC^+}\setminus\{0\}$, so $p\in S_{b,\reg}$ for some 
$b \in I\setminus \{*\}$. As explained in \S 2.6, there is a neighbourhood 
$\calU$ of $p$ in $\fraka$ such that the preimage $\pi^{-1}(\calU)$ in $M$ 
is a bundle over $K/K^b$ with fiber an open neighbourhood of $(o,p)$ in 
$\Sigma^b \times S_b$. The subgroup $W^b \subset W$ generated by roots $\al\in\Lambda_b$ 
is naturally identified with the Weyl group of $\Sigma^b$. Now suppose that 
$u\in\Cinf(\fraka)^W$. Then the restriction of $u$ to $\calU$ can be considered 
as a smooth $W^b$-invariant function on (some neighbourhood of a point
$(0,p)$ in) $S^b\oplus S_b$.  By the inductive hypothesis,  $\pi^* u$
can be identified with a smooth $K^b$-invariant function on a neighbourhood
of $(o,p) \in \Sigma^b \oplus S_b$. Since $b$ is arbitary, this proves that 
$\pi^* u \in \Cinf(M\setminus\{o\})^K$.

It remains to prove that $\pi^* u$ is also smooth near $o$. At the same
time we must also start the induction, proving that $\pi^* u$ is smooth on
$M$ for symmetric spaces of rank one, but since the only issue in that case
is to prove smoothness at $o$, this is the same argument.

We proceed as follows. Let $L$ be the operator on $\fraka$ induced by 
$\Delta_\rad$ on $\fraka_\reg$; according to Lemma~\ref{lemma:Delta-rad-L}, 
$L$ preserves $\Cinf(\fraka)^W$. We have already remarked that
since $u \in \Cinf(\fraka)^W$ is locally Lipschitz, the same is true
of $\pi^* u$. Moreover, $Lu\in\Cinf(\fraka)^W$, so $\pi^* (Lu)$
is also locally Lipschitz on $M$. By the induction, $\pi^*(Lu)$ agrees 
with the smooth function $f=\Delta (\pi^* u)$ away from $o$. Hence 
$\Delta (\pi^* u)$ is a distribution differing from the locally Lipschitz 
function $\pi^*(Lu)$ by a distribution supported at $o$. However,
$\nabla \pi^* u \in L^\infty_\loc$, so in particular $\pi^* u\in H^1_{\loc}$,
which implies that $\Delta (\pi^*u) \in H^{-1}_{\loc}$. Furthermore, since
it is locally Lipschitz, $\pi^*(Lu) \in H^1_{\loc}$ too. Therefore the
difference $g = \Delta (\pi^* u) - \pi^*(Lu) \in H^{-1}_{\loc}$. 
If $\dim M \geq 2$, no element of $H^{-1}_{\loc}$ can be supported at $o$, 
so $g=0$. If $\dim M = 1$, then the $K$ is finite and the same conclusion
is trivial. 

We have now proved that $\Delta \pi^* u$ is locally Lipschitz, and 
$\Delta\pi^* u=\pi^*(Lu)$. Now repeat the argument with $u$ replaced by 
$Lu$ to conclude that $\Delta^j \pi^* u$ is locally Lipschitz for every 
$j \geq 1$. By elliptic regularity, $\pi^* u \in \Cinf(M)$, and this
completes the proof.
\end{proof}

This result extends to the compactifications, as is easily seen from the
proof of Proposition~\ref{prop:A-M-W-K}: in the inductive step, we
merely need to compactify the base space $S_b$ of the family.

\begin{prop}
The map $\pi^*$ gives isomorphisms $\Cinf(\ofraka)^W \to \Cinf(\olM)^K$
and $\Cinf(\tfraka)^W \to \Cinf(\tilde{M})^K$. 
\end{prop}

\subsection{Invariant partitions of unity}
We now introduce $W$-invariant partitions of unity on $\fraka$ which
are compatible with the structures of the compactifications $\hat{\fraka}$
and $\ofraka$. The lifts of these partitions of unity are of course
$K$-invariant partitions of unity on $M$ compatible with the structures
of the corresponding compactifications.

Each (open) face $S_b^+$ of the closed positive Weyl chamber $\overline{\calC^+}$ 
in $\fraka$ is an open set in a unique $S_b$, $b \in I$, and therefore we
may index the set of all such faces $S_b^+$ by a subset $I^+ \subset I$.

We first consider invariant partitions of unity on $\fraka$:
\begin{Def}
A partition of unity $\{\chi_b:\ b\in I^+\}$ on $\overline{\calC^+}$ is 
$W$-adapted if each $\chi_b$ is the restriction to $\overline{\calC^+}$ 
of some $\chi_b' \in \Cinf(\fraka)^W$, and moreover if 
$\supp\chi_b\cap S_c^+=\emptyset$ except when $S_b^+\subset S_c$.
\end{Def}
\begin{rem}
Since $\sum \pi^* \chi_b' = \pi^*(\sum \chi'_b) = 1$, the
lifts $\pi^* \chi'_b$ are a smooth $K$-invariant partition of unity on $M$. 
\end{rem}

No conditions have been imposed on the $\chi_b$ at infinity, so this
partition of unity is only useful for studying local properties. 
To go further, let $\widehat{\calC^+}$ be the closure of $\overline{\calC^+}$ 
in the radial compactification $\hat\fraka$.
\begin{Def}\label{def:W-ahat-adapted}
A partition of unity $\{\chi_b:\ b\in I^+\}$ on $\widehat{\calC^+}$ is 
$(W,\hat\fraka)$-adapted if 
\begin{enumerate}
\item
each $\chi_b$ is the restriction to $\widehat{\calC^+}$ of an element
of $\Cinf(\hat\fraka)^W$,
\item
$\supp\chi_*$ is a compact subset of $\fraka$, and
\item
$\supp\chi_b\cap \widehat{S_c^+}=\emptyset$ unless $S_b^+\subset S_c$; 
here $\widehat{S_c^+}$ is the closure of $S_c^+$ in $\hat\fraka$.
\end{enumerate}
\end{Def}
The restriction that $\chi_b$ be supported sufficiently near to $S_b^+$,
i.e.\ (iii), ensures that $L_{b,\rad}$ is a good model for 
$\Delta_{\rad}$ on its support. On the other hand, (ii) guarantees
that the partition of unity is not trivial: i.e.\ that $\chi_*\not\equiv 1$. 

\begin{lemma}
There exists a $(W,\hat\fraka)$-adapted partition of unity.
\end{lemma}
\begin{proof}
We first construct a partition of unity on $\hat\fraka$ with the appropriate 
support properties, then average it over $W$. 

For any root $\alpha$, let $\widehat{W_\alpha}$ denote the closure of 
the wall $W_\alpha$ in $\hat\fraka$. Also, set
\begin{equation*}
\widehat{W_{\alpha,\pm}}= 
\overline{\alpha^{-1}(\RR^{\pm})}\setminus \widehat{W_\alpha};
\end{equation*}
this is the closure in $\hat\fraka$ of the set where $\alpha > 0$, 
respectively $\alpha < 0$, minus the closure of the wall. We say that 
$\alpha>0$ on $\widehat{W_{\alpha,+}}$ and $\alpha<0$ on $\widehat{W_{\alpha,-}}$ 
and $\alpha=0$ on $\widehat{W_\alpha}$.

Each face of each Weyl chamber is defined by a map $\mu: \Lambda \to \{0,+,-\}$, 
corresponding to whether each root is $>0$, $<0$ or $=0$ on that face. Denote
the space of all such maps by $\calP$. Certain $\mu \in \calP$ correspond to
empty faces (for instance if one requires that both $\alpha$ and $-\alpha$ 
are positive), so we let $\calP_0$ be the subset of $\mu$ for which the
corresponding face is nonempty. To any $\mu \in \calP_0$ such that 
$\mu(\alpha)\neq 0$ for at least one $\alpha$ we associate the 
relatively open set 
\begin{equation*}
U_\mu=\left(\bigcap\{\widehat{W_{\alpha,+}}:\ \mu(\alpha)>0\}\right)
\cap \left(\bigcap\{\widehat{W_{\alpha,-}}:\ \mu(\alpha)<0\}\right) \subset
\hat{\fraka};
\end{equation*}
with $*$ corresponding to the map $\mu\equiv 0$ we also set $U_*=\fraka$. 

The collection $\calU = \{U_\mu\}$ is an open cover of $\hat\fraka$, and 
we choose a partition of unity $\{\psi_\mu\}$ subordinate to it. Every
$w\in W$ is an endomorphism of $\fraka$, and extends to a diffeomorphism
of $\hat\fraka$. To each such $w$, if $\alpha \in \Lambda$, then $w_*\mu$ 
is the map which assigns to $w^*\alpha$ the value $\mu(\alpha)$. Finally, let
\begin{equation*}
\phi_\mu=\frac{1}{|W|}\sum_{w\in W}w^*\psi_{\mu}.
\end{equation*}
Then $\sum_\mu\phi_\mu=1$ and each $\phi_\mu$ is clearly $W$-invariant.

If the face corresponding to some $\mu$ is not contained in $\overline{C^+}$,
then $U_\mu\cap\widehat{C^+}=\emptyset$. Indeed, for any such $\mu$ there
is a positive root $\alpha$ such that $\mu(\alpha)<0$, so $\alpha<0$ on $U_\mu$,
which means that $U_\mu$ does not intersect the closed positive chamber.

Note also that for any $\mu \in \calP_0$, there is a unique $\mu_+ = w_* \mu$
which is $\geq 0$ on all positive roots. Since $w^*\psi_\mu$ is supported in 
$w^{-1}(U_\mu)=U_{w_*\mu}$, we have $\supp w^*\psi_\mu\cap \widehat{\calC^+}
=\emptyset$ unless $w_*\mu=\mu_+$.

Now suppose that $S_b^+$ is a face of $\overline{C^+}$. Clearly 
$S_b\subset S_c$ if and only if for every root $\alpha$, $\alpha\equiv 0$
on $S_c$ implies $\alpha\equiv 0$ on $S_b$. Thus if $S_b\not\subset S_c$,
then there is a root $\alpha$, which we may assume is positive, which 
vanishes identically on $S_c$ but not on $S_b$. In particular,
if $\nu$ is the map corresponding to $b\in I^+$, then $\nu(\alpha)$ is 
positive (since $b\in I^+$), hence non-zero, and so $U_\nu\cap 
\widehat{S_c^+}=\emptyset$ by the definition of $U_\nu$.

Finally, combine each $W$-orbit of $\phi_\mu$ into a single term
\begin{equation*}
\chi_b=\chi_\nu=\sum_{w\in W}\phi_{w_*\nu}.
\end{equation*}
Now, for $w\in W$, $\supp w^*\psi_{v_*\nu}\cap \widehat{\calC^+}=\emptyset$
unless $w^*v^*\nu=(v^*\nu)_+=\nu_+=\nu$ since $\nu$ is $\geq 0$ on positive
roots. On the other hand,
if $w^*v^*\nu=\nu$, and $c$ is as in the previous paragraph,
then $\supp w^*\psi_{v^*\nu}\cap \widehat{S_c^+}\subset U_\nu
\cap \widehat{S_c^+}=\emptyset$. Therefore, for every
$v,w\in W$, $\supp w^*\psi_{v^*\nu}\cap \widehat{\calC^+}=\emptyset$.
This shows that $\supp\chi_b\cap \widehat{S_c^+}=\emptyset$, which
finishes the proof.
\end{proof}

\begin{Def}\label{def:W-aol-adapted}
A partition of unity $\{\chi_b:\ b\in I^+\}$ on $\overline{\calC^+}$ 
is $(W,\ofraka)$-adapted if
\begin{enumerate}
\item
each $\chi_b$ is the restriction to $\overline{\calC^+}$ of an element
in $\Cinf(\ofraka)^W$,
\item
$\supp\chi_b\cap \overline{S_c^+}=\emptyset$
unless $S_b^+\subset S_c$ (where $\overline{S_c^+}$ is the closure
of $S_c^+$ in $\ofraka$), and
\item
$\supp\chi_b\subset S_b\cap \Omega_b$, where $\Omega_b$ is a compact subset
of $S^b$ (and in particular, $\chi_*$ has compact support since $S_*=\{0\}$). 
\end{enumerate}
\end{Def}

\begin{lemma}
There exists a $(W,\ofraka)$-adapted partition of unity.
\end{lemma}
The proof proceeds just as for the $(W,\hat\fraka)$-adapted case,
and so we omit it.

\begin{rem}\label{rem:chi_0-supp}
Note that if $\{\chi_b:\ b\in I^+\}$ is a $(W,\ofraka)$-adapted partition
of unity, then there exists $T=(T_1,\ldots,T_n)$ with all
$T_j>0$ such that $\supp\chi_0\subset\calO(T)$, so $\chi_0$ localizes
away from the walls. In addition, all the $N$-body features of the analysis
are already present on $\supp\chi_0$.
\end{rem}

\section{Differential operators, function spaces and mapping properties}
In this section we explain the appropriate spaces of differential
operators and functions
of finite regularity that are used later.

We start with differential operators, or more specifically, $K$-invariant
operators acting on $K$-invariant function spaces. If $P$ is such an operator
and $P_\rad$ its radial part, then since $\Cinf_c(M)^K$ is identified with 
$\Cinf_c(\fraka)^W$, and $\Cinf_c(M)^K$ is dense in every function space 
we wish to study, we can regard $P_\rad$ either as a map $\Cinf_c(M)^K\to\Cinf_c(M)^K$ 
(i.e.\ as the restriction of $P$), or as a map $\Cinf_c(\fraka)^W\to\Cinf_c(\fraka)^W$. 
In the former case, $P_\rad$ is a differential operator on $M$ with $\Cinf$ coefficients,
while in the latter case, $P_\rad$ is a differential operator whose coefficients
on $\fraka_{\reg}$ are smooth, hence gives a map $\Cinf(\fraka_\reg)^W\to\Cinf(
\fraka_\reg)^W$, which restricts to a map $\Cinf_c(\fraka)^W\to\Cinf_c(\fraka)^W$. 
One could define the appropriate space of differential operators directly
on $\fraka$, but one must take care to see their uniformity near the walls.
We proceed instead by identifying functions on neighbourhoods of the walls
in $\fraka$ with neighbourhoods in a product model.

Let $\{\chi_b: b\in I^+\}$ be a $(W,\ofraka)$-adapted partition of unity, and 
fix diffeomorphisms 
\begin{equation*}
\Psi_b: \mbox{supp}\, \chi_b \hookrightarrow S_b \oplus S^b.
\end{equation*}
Then to any $W^b$-invariant function $u$ on $S_b \oplus S^b$, we can associate 
a $K^b$-invariant function $\tilde{u}_b$ on $S_b \times \Sigma^b$, and
conversely the restriction of such a $K^b$ invariant function
to $S_b \oplus S^b$ is $W^b$-invariant.
If $\mbox{supp}\, u
\subset S_b \times \calV_b$, then $\mbox{supp}\,(\tilde{u}_b) \subset
S_b \times \widetilde{\calV}_b$, where $\widetilde{\calV}_b$ is a bounded set
containing the origin in $\Sigma^b$.

The operators we shall single out are generated by translation 
invariant operators on $S_b$ and arbitrary differential operators on 
the bounded set $\widetilde{\calV}_b \subset \Sigma^b$. We also require
the operators of multiplication by functions in both $\Cinf(\hat\fraka)^W$ and
$\Cinf(\ofraka)^W$, since elements of the former are required in the partition 
of unity patching the local models, while the form of the Laplacian requires 
the latter; these requirements suggest that we allow multiplication by functions in
$\Cinf(\widetilde M)^K\equiv\Cinf(\tfraka)^W$.

\begin{Def}\label{Def:Diffss}
The space $\Diffso^m(M)$ consists of all differential operators
$P:\Cinf_c(M)\to\Cinf_c(M)$ of order $m$ which are $K$-invariant, and 
such that for each $b\in I^+$, the $K$-radial part $P_\rad$ of $P$, restricted to
functions supported in $\pi^{-1}(\supp\chi_b)$, is the $K^b$-radial part $Q_\rad$
of a differential operator $Q$ on $S_b\times\widetilde{\calV}_b$ which is a 
linear combination of products of translation invariant operators on $S_b$ and
differential operators on $\widetilde{\calV}_b$, with
coefficients in $\Cinf(\widetilde{S_b}\times\widetilde{\calV}_b)$.
\end{Def}

\begin{rem}\label{rem:diff-ops-0-1}
The subscript $o$ has been included in this notation because this
space of operators depends on the choice of origin in $M$. We note  
also that this definition only restricts the behaviour of these operators
near infinity. Finally, recall from Remark~\ref{rem:chi_0-supp}
that $\chi_0$ (considered as a function on $\overline{\calC^+}$)
is supported in $\calO(T)$ for some $T$ with all $T_j>0$. Thus, for $b=0$ 
the requirement is that $P_\rad$ restricted to $\calO(T)$ is a linear combination 
of translation invariant differential operators in $\fraka$ with coefficients
in $\Cinf(\tfraka)$. Apart from the localization to $\calO(T)$,
this is {\em exactly} the definition of $N$-body differential
operators $\Diff_{\mathrm{sc}}^*(\hat\fraka,\calC)$, $\calC=\calS\cap\pa
\hat\fraka$, in \cite{Vasy:Propagation-Many}.

The use of the product spaces $S_b\times\widetilde{\calV}_b$ is motivated
by the results of \S~2.6; see in particular Remark~\ref{rem:Mb-Mt}.
\end{rem}

We now discuss the associated $L^2$-based Sobolev spaces. The basic
$L^2$ space is, of course, $L^2(M,dg)^K$, which is identified with an 
$L^2$-space with respect to the degenerate measure on $\fraka$,
$dg_0 = \pi_* dg := \eta\, da$ where $\pi:M\to\calC^+$; note that
$\eta$ extends to be $W$-invariant function on $\fraka$.
There is an explicit formula \cite[Ch. 1, Theorem 5.8]{Helgason:Groups}
\begin{equation}\label{eq:eta-def}
\eta(a) = \prod_{\alpha \in \Lambda^+} (\sinh \alpha(a) )^{m_\alpha}, \qquad
a \in C^+.
\end{equation}
Notice that $\eta(a)$ is $\Cinf$ and
strictly positive on $\fraka_\reg$, but 
degenerates like various powers of the distance function along the Weyl
chamber walls,
i.e.\ where various roots $\alpha$ vanish. Then
\begin{equation*}
L^2(M,dg)^K\equiv L^2(\calC^+,dg_0)\equiv L^2(\fraka,\frac{1}{|W|}\,dg_0)^W
\end{equation*}
as Hilbert spaces; of course, the norms of the last terms are equivalent
without the constant factor $|W|^{-1}$.

As $M$ is a non-compact space, there are various spaces of $K$-invariant
Sobolev functions that we can associate to it.
We need the spaces that correspond to $\Diffso(M)$,
which was in turn constructed to accommodate both the Laplacian and
multiplication by cutoffs in $\Cinf(\hat\fraka)$.
For $b\in I^+$, we let
\begin{equation*}
\eta_b(a)=\prod_{\alpha \in \Lambda_b^+} (\sinh \alpha(a) )^{m_\alpha}, \qquad
a \in C^+,
\end{equation*}
note that on $\supp\chi_b$ we can identify $\eta_b\,da^b$ with the push-forward
of the Riemannian measure $dg_b$ on $\Sigma^b$ to the positive
chamber of $S^b$. Moreover, the other
positive roots $\alpha\in\Lambda^+\setminus\Lambda_b^+$ tend to $+\infty$
on $\supp\chi_b$, so $e^{-2(\rho-\rho_b)}
\prod_{\alpha \in \Lambda^+\setminus\Lambda_b^+}
(\sinh \alpha(a) )^{m_\alpha}$ is bounded from below and above by positive
constants. Correspondingly, for functions in $L^2(M,dg)^K$ supported in
$\supp\chi_b$, the $L^2(M,dg)$-norm is equivalent to the
$L^2(S_b\times\Sigma^b;e^{2(\rho-\rho_b)}\,da_b\,dg_b)$-norm;
here $da_b$ is the Euclidean density on $S_b$.
We now define the Sobolev spaces as follows.

\begin{Def} The space $\Hso^s(M)^K$ is the set of distributions 
$u \in \calD'(M)^K\equiv\calD'(\fraka)^W$ with the property that
$e^{\rho-\rho_b}\left((\Psi_b)_*(\chi_b u)
\right)^{\widetilde{}}_b \in H^s(S_b \times \Sigma^b)$. (Because the support
is bounded in the second factor, there are no subtleties involving noncompact
supports in this condition.)
\end{Def}

\begin{rem}
Continuing Remark~\ref{rem:diff-ops-0-1}, note that for $b=0$ the requirement
is simply that $e^{\rho}\chi_0 u\in H^s(\fraka)$, i.e.\ $\chi_0 u$
is in the weighted Sobolev space $e^{-\rho}H^s(\fraka)$ (where $H^s(\fraka)$
is the standard Sobolev space on the vector space $\fraka$).
\end{rem}

\begin{rem}
We could have equally well defined these adapted classes of differential operators
and Sobolev spaces using the identification of neighbourhoods of the supports
of elements of a $(\hat\fraka,W)$-adapted partition of unity, i.e.\ by working
on conic neighbourhoods of the $S_b$. This would require that definitions 
be made inductively on the rank, since we would no longer be working in compact subsets 
of the subsystems $\Sigma^b$. 
\end{rem}

If $s\geq 0$ is an integer, this means that for any
$A\in\Diffso^k(M)$ with $k\leq s$,
\begin{equation*}
Au\in L^2(M,dg)^K.
\end{equation*}
Indeed,  by the definition of
$\Diffso(M)$,
the latter statement is equivalent to requiring that for any
translation invariant differential operator $P$ of order $k\geq 0$
on $S_b$ and for any differential operator $Q$ of order $l\geq 0$ on
$\Sigma^b$, with $k+l\leq s$,
\begin{equation*}
PQ \left((\Psi_b)_*(\chi_b u)
\right)^{\widetilde{}}_b\in L^2(e^{2(\rho-\rho_b)}\,da_b\,dg_b).
\end{equation*}
Since commuting the weight through $P$ introduces lower order differential
operators, this is easily seen to be equivalent to
\begin{equation*}
PQ e^{\rho-\rho_b}\left((\Psi_b)_*(\chi_b u)
\right)^{\widetilde{}}_b\in L^2(da_b\,dg_b),
\end{equation*}
for all $P$ and $Q$ as above,
which is the definition of the Sobolev spaces.

A key property that a parametrix $G$ for $\Delta_\rad-\ev$
should have is that its error $F
=(\Delta_\rad-\ev)G-\Id$ should
be a compact operator, say on $L^2(M,dg)^K$. We can achieve this
by showing that $F$ maps into a positive order Sobolev space
with additional decay at infinity.
Thus, we also consider spaces of functions on $\tfraka$ with
some specified rate of decay at the boundary. To this end, we introduce
the total boundary defining function 
\[
\bdef = \prod_{b\in I\setminus\{*\}} \bdef_b,
\]
where $\bdef_b$ is a defining function for the face $\widetilde{F_b}$ of $\tfraka$. 
Note that $\hat x$ agrees with $\bdef$ up to a smooth non-vanishing positive factor,
as follows by considering $\tfraka$ as a blow-up of $\hat\fraka$.

Supposing that $\bdef$ is $W$-invariant, we then define 
\[
\bdef^\delta \Hso^s(M)^K = \{u = \bdef^\delta v: v \in \Hso^s(M)^K\}
\]
(which by the remark above is the same as $\hat x^\delta H^s_W(\tfraka)$). 

\begin{prop}\label{prop:ss-diff-op-bds}
For any $s,\delta\in\Real$, $\Diffso^m(M):
x^\delta \Hso^s(M)^K\to x^\delta \Hso^{s-m}(M)^K$.
\end{prop}

\begin{proof}
Both the Sobolev spaces and the differential operators are defined
by localization to $S_b\times\widetilde{\calV}_b$, and on these the claims
are clear.
\end{proof}

It is crucial for us that parametrix constructions can be localized on
$\hat\fraka$. This is reflected by the following proposition.

\begin{prop}\label{prop:a-hat-commutes}
The multiplication operators $\phi\in\Cinf(\hat\fraka)^W$ commute
with operators $P\in \Diffso^k(M)$ to top order,
i.e.\ $[P,\phi]\in x \Diffso^{k-1}(M)$. Thus,
$[P,\phi]:x^\delta \Hso^{s+m-1}(M)^K\to x^{\delta+1} \Hso^s(M)^K$.
\end{prop}

\begin{rem}
The analogue of this result has been widely used in $N$-body scattering.
There is a much larger class of (pseudo-)differential
operators which commute to top order with every $P\in \Diffso^k(M)$, 
and which can be used to microlocalize, see \cite{Vasy:Propagation-Many}.
\end{rem}

\begin{proof}
Using a partition of unity, we assume that $P$ is supported in
$\pi^{-1}(\supp\chi_b)$. Valid local coordinates on $\hat\fraka$ near
$\hat S_b$ are of the form $\frac{\alpha_j(a)}{|a|}$, $a\in\fraka$,
where the $\alpha_j$
are linearly independent simple roots that
vanish on $S_b$, as well as coordinates on
$\hat S_b$. Thus, in a neighbourhood of $\hat S_b$ (which includes $\supp P$)
\begin{equation*}
\phi=\phi|_{\hat S_b}+\sum_j \frac{\alpha_j(a)}{|a|}\phi_b,
\end{equation*}
with $\phi_b$ smooth in this open subset of $\hat\fraka$. In particular,
its commutator with $P$ is in $\Diffso^{k-1}(M)$. Using this expansion now
it is straightforward to complete the proof.
\end{proof}

Specializing these results to the Laplacian, we deduce that
for any $s,\delta \in \RR$ and $\ev \in \CC$,
\[
\Delta_\rad - \ev: \bdef^\delta \Hso^{s+2}(M)^K \longrightarrow 
\bdef^\delta \Hso^s(M)^K.
\]

Ultimately, of course, we are interested in inverting this operator,
and as usual, this will rely on its ellipticity.

\begin{Def}
We say that $P\in\Diffso(M)$ 
is {\it radially elliptic} if for every $b\in I^+$,
there is
an operator $Q=Q_b\in\Diffso^m(S_b\times\Sigma^b)$
as in Definition~\ref{Def:Diffss} that is symbol-elliptic.
\end{Def}

\begin{rem}
We emphasize that symbol-ellipticity in $\Diffso^m(S_b\times\Sigma^b)$ is 
a {\em uniform} condition near infinity in $S_b$.

In particular, for $b=0$, such a differential operator has the form
$\sum_{|\gamma|\leq m} p_\gamma(a) D^\gamma$, with $p_\gamma$ smooth
on the closure of $\calO(T)$ in $\tfraka$, $T_j>0$ for all $j$.
Symbol ellipticity then is the requirement that
$\sum_{|\gamma|=m} p_\gamma(a) \xi^\gamma$ never vanish for
$(a,\xi)$ in the closure of $\calO(T)\times\fraka^*$ in $\tfraka
\times\fraka^*$.
\end{rem}

Clearly, $\Delta_\rad$ is radially elliptic. Indeed, we can take
$Q_b=T_b+\Delta_{\Sigma^b}$. Thus, one can use the standard parametrix
construction for $\Delta_\rad-\ev$; indeed, even the standard {\em large
spectral parameter} construction works, i.e.\ we can precisely analyze
$|\ev|\to\infty$.

\section{Complex scaling}
As explained in the introduction, there are two main tools in our proof of 
the analytic continuation of $\Delta_\rad$: construction of the parametrix,
which takes place in the $b$-calculus on $\tfraka$, and the method of complex
scaling. In this section we focus on the second of these, and shall
review this method, which produces a holomorphic family of operators
for which the essential spectrum is shifted away from the positive real axis.

The ingredients needed in this procedure are a family of (possibly unbounded)
operators $U_\theta$ acting on $L^2(\fraka)^W$, for $\theta$ lying in some
contractible domain $D\subset\Cx$, and a dense subspace of `analytic vectors' 
$\AAA \subset L^2(\fraka)^W$, such that:
\begin{enumerate}
\item
$U_0=\Id$ and for $\theta\in D\cap\Real$, $U_\theta$ is unitary on $L^2(\fraka)^W$ and 
bounded on all Sobolev spaces; 
\item
For $f\in\AAA$, the map $\theta \to U_\theta f$ extends analytically 
from $D \cap \Real$ to all of $D$ with values in $L^2(\fraka)^W$;
\item
For each $\theta \in D$, the subspace $U_\theta\AAA$ is dense in $L^2(\fraka)^W$.
\end{enumerate}

By (i), we can define $\Delta_{\rad,\theta}= U_{\theta}\Delta_\rad 
U_{\theta}^{-1}$ directly when $\theta\in\Real$. We shall show below
that the coefficients of this operator extend analytically in $\theta$ to 
the sector $|\im\theta|<\frac{\pi}{2}$; hence for fixed $f \in \Cinf_c(M)$, 
$\theta \to \Delta_{\rad,\theta} f$ is analytic in this same region. 
We must actually prove that the family $\Delta_{\rad,\theta}$ is analytic of 
type A, see Proposition~\ref{prop:type-A} below. 
The resolvent of the scaled radial Laplacian, $(\Delta_{\rad,\theta}-\ev)^{-1}$,
will be constructed by parametrix methods. From this we can deduce the
meromorphic continuation of $R(\ev)$ from the equality $(\Delta_{\rad,\theta} 
- \ev)^{-1} = U_\theta R(\ev)U_{\theta}^{-1}$, which is initially valid when 
$\ev$ is in the resolvent set common to both operators and $\theta$ is real. 
In fact, we prove only that the matrix element $\langle f,R(\ev)g\rangle$ continues 
meromorphically to $D$ whenever $f,g \in \AAA$; this is sufficient for purposes 
of spectral theory. 

\subsection{Complex dilations} Let $\frakp_{\Cx}$ denote the complexification
of $\frakp$ and $D$ some domain in $\Cx$ containing $0$, and define 
\[
\Phi: D \times \frakp \longrightarrow \frakp_{\Cx};
\qquad \Phi(\theta,X) = e^{\theta}X.
\]
We also denote that restriction of $\Phi$ to $D \times \fraka \longrightarrow 
\fraka_{\Cx}$ by $\Phi$, and often write $\Phi_\theta(X) = \Phi(\theta,X)$.
Identifying $\frakp$ and $M$ by the exponential map, for $\theta\in\RR\cap D$
$\Phi_\theta$ is the diffeomorphism on $M$ given by dilating by the factor
$e^\theta$ along geodesic rays emanating from $o$.  

When $\theta \in \Real$, the
induced family of unitary operators $U_{\theta}$ on
\begin{equation*}
L^2(M)^K
\equiv L^2(\fraka,|W|^{-1}\pi_*dg)^W
\end{equation*}
is defined by
\begin{equation}\label{eq:U_theta-def}
(U_{\theta}f)(a) = (\det\,D_\theta\Phi)^{\frac12}f(e^\theta a)
=J_\theta^{\frac12}(\Phi_\theta^* f)(a),\ a\in\fraka;
\end{equation}
the Jacobian prefactor, which is calculated with respect to
the density $\pi_* dg=\eta\,da$ in \eqref{eq:eta-def}, makes this map unitary.
Explicitly, with $n=\dim\fraka$,
\begin{equation*}
J_\theta(a)=(\det\,D_\theta\Phi)(a)=w^n\frac{\eta(wa)}{\eta(a)}=w^n
\prod_{\alpha \in \Lambda^+}
\left(\frac{\sinh (w\alpha(a))}{\sinh (\alpha(a))} \right)^{m_\alpha},
\ a\in C^+.
\end{equation*}
Note that $J_\theta$ does not vanish for $|\im\theta|<\frac{\pi}{2}$.
The product can be replaced by one over $\Lambda$, if $m_\alpha$ is
replaced by $m_\alpha/2$, and then the formula is valid on all of $\fraka$;
this also shows that $J_\theta$ is $\Cinf$ on $\fraka$.

While the use of $U_\theta$ fits nicely into the Aguilar-Balslev-Combes theory, 
one could also work with $\Phi_\theta^*$ directly, which would be closer
in spirit to the microlocal complex deformations of Sj\"ostrand and Zworski
\cite{Sjostrand-Zworski:Complex}.

\begin{lemma}
For $\theta\in\Real$, $\Phi_\theta$ extends to a `conormal diffeomorphism' 
of $\ofraka$, in the sense that $\Phi^*_\theta:S^m(\ofraka)
\mapsto S^{m w}(\tfraka)$, where $w=e^\theta$ and $S^m(\ofraka)$ 
denotes the symbol space. In addition, it extends to a diffeomorphism of
$\tfraka$.
\end{lemma}

\begin{proof}
The first claim is easy to check since the effect of dilations is that
roots $\alpha$ are multiplied by $e^\theta$: $\Phi_\theta^*\alpha(a)
=\alpha(e^\theta a)=e^\theta\alpha(a)$, and the negative exponentials
of the simple roots define
the smooth structure of $\ofraka$ in a neighbourhood of $\overline{C^+}$.

The second claim follows from either description of $\tfraka$. Indeed,
$\Phi_\theta$ extends to a diffeomorphism of $\hat\fraka$, and then
lifts to its blow-up $\tfraka$. Alternatively, the logarithmic total
boundary blow-up replaces the defining functions $e^{-\alpha_j}$ of
$\ofraka$ in $C^+$ by $\alpha_j^{-1}$, so
$\Phi_\theta$ extends to a diffeomorphism of the this blow-up, which then
lifts to $\tfraka$.
\end{proof}

\begin{lemma}
The Jacobian determinant $J$ extends to an analytic nonvanishing function in the region
\begin{equation*}
D = \{\theta\in\Cx:\ |\im\theta|<\frac{\pi}{2}\}.
\end{equation*}
In addition, $J$, $J^{1/2}$ and $J^{-1/2}$ are conormal $K$-invariant
functions on $\olM$, equivalently, conormal $W$-invariant functions on $\ofraka$.
\end{lemma}

We shall need a slight generalization of this definition later. Let $\Phi_{\theta,T}$ 
be a $W$-invariant diffeomorphism of $\fraka$ which is the identity on the ball 
$B_T(0)$ and equals the dilation by $e^\theta$ outside a larger ball, and which 
depends analytically on $\theta$. For example, fix $T>0$ and a nondecreasing cutoff 
function $\phi\in\Cinf(\Real;[0,1])$ which equals $1$ near $\infty$ and vanishes 
on $[0,T]$, and define 
\begin{equation*}
\Phi_{\theta,T}(a)= e^{\phi(r)\theta}a;
\end{equation*}
then $\Phi_{\theta,T}(a)=a$ if $|a| \leq T$, and $\Phi_{\theta,T}(a)=
e^\theta a$ for $|a| \geq T' > T$, and $\theta \mapsto \Phi_{\theta,T}(a)$ is
analytic. It is clear that $\Phi_{\theta,T}$ is a diffeomorphism when 
$\theta$ is real and near $0$, and that it extends analytically to complex $\theta$. 

\begin{lemma}
There exists $\delta>0$ such that $\Phi_{\theta,T}:M\to M$ is a
diffeomorphism when $\theta \in \Real$, $e^\theta>1-\delta$. In addition, 
$(\det D\Phi_{\theta,T})^{1/2}$ extends analytically to
the region 
\begin{equation*}
\{\theta\in\Cx:\ |\im\theta|<\frac{\pi}{2},\ e^\theta\nin(-\infty,1-\delta)\}.
\end{equation*}
\end{lemma}

Now set 
\begin{equation}
(U_{\theta,T}f)(a)=(\det D\Phi_{\theta,T})^{1/2} f(\Phi_{\theta,T}(a)).
\label{eq:defuo}
\end{equation}

Because of the simple geometric nature of the transformations 
$U_\theta$ and $U_{\theta,T}$, we may define the families of differential operators 
\[
\Delta_{\rad,\theta}=U_{\theta}\Delta_\rad U_{\theta}^{-1}, 
\qquad \Delta_{\rad,\theta,T}=U_{\theta,T}\Delta_\rad U_{\theta,T}^{-1},
\]
without worrying about functional analytic issues of domain. These are 
$W$-invariant on $\fraka$, with coefficients depending analytically on $\theta$ 
in the region $D=\{\theta:\ |\im\theta|<\pi/2\} \subset \Cx$.

Indeed, we have already seen that $J_\theta^{1/2}$ extends to be analytic
and nonvanishing on $D$. Since
\begin{equation*}
U_\theta\Delta_{\rad}U_\theta^{-1}=J_\theta^{1/2}\Phi_\theta^*
\Delta_{\rad}(\Phi_\theta^{-1})^*J_\theta^{-1/2},
\end{equation*}
we only need to consider $\Phi_\theta^*\Delta_{\rad}(\Phi_\theta^{-1})^*$.
Now, the $\Phi_\theta^*$-conjugates of the principal part
$\Delta_\fraka$ (as well
as the first order constant coefficient terms) continue to $\Cx\setminus \Real^-$
(and even to a larger Riemann surface). For example, $\Phi_\theta^*
\Delta_\fraka (\Phi_\theta^{-1} )^*
= e^{-2\theta}\Delta_\fraka$.
However, the coefficients $\coth \al$ only continue 
up to $|\im\theta|=\frac{\pi}{2}$, 
and genuine singularities appear in these continuations on this ray.

The coefficients of $\Delta_{\rad,\theta}$ are thus smooth on $\fraka$ 
when $|\im\theta|<\frac{\pi}{2}$, but we also require information about their
behaviour at $\pa \tfraka$. 

\begin{prop}\label{prop:type-A}
If $\theta\in\Cx$ has $|\im\theta|<\frac{\pi}{2}$, then
$\Delta_\theta$ is a (polyhomogeneous) conormal
b-differential operator on $\olM$. Its radial part
$\Delta_{\rad,\theta}$ is radially elliptic. The operators
\begin{equation*}
L_{b,\theta}=T_{b,\theta}+\Delta_{b,\rad,\theta},\ b\in I^+,
\end{equation*}
on $L^2(S_b\times\Sigma^b;e^{2(\rho-\rho_b)}\,da_b\,dg_b)$,
are product models for $\Delta_{\theta,\rad}$ in the sense that if
$\chi_b\in\Cinf(\hat\fraka)$ satisfies (i) and (iii) of
Definition~\ref{def:W-ahat-adapted} then
\begin{equation*}
E_{b,\theta}\chi_b=(\Delta_{\rad,\theta}-L_{b,\theta})\chi_b
\in x^\infty\Diffso^1(M).
\end{equation*}
Also, $\theta \to \Delta_{\rad,\theta}$ is an analytic type-A family 
on $L^2(\tfraka)^W$ with domain $\Hso^2(M)^K$.
\end{prop}

\begin{proof}
The first part is easy from the explicit formula. We remark that
$L_{b,\theta}$ is defined using the dilations on
$S_b\times\Sigma^b$ and the Jacobian corresponding to the $L^2$-space
\begin{equation*}
L^2(S_b\times\Sigma^b;e^{2(\rho-\rho_b)}\,da_b\,dg_b).
\end{equation*}
Thus, $\Delta_{b,\theta}$ is indeed the complex scaled $\Delta_b$,
defined by \eqref{eq:U_theta-def} with $M$ replaced by $\Sigma^b$.
Moreover, with $w=e^\theta$, $\tilde\rho=\rho-\rho_b$,
\begin{equation*}\begin{split}
&T_{b,\theta}={\mathcal J}_\theta^{1/2}
(w^{-2}\Delta_{S_b}+2w^{-1}H_{\tilde\rho}){\mathcal J}_\theta^{-1/2},
\qquad {\mathcal J}_\theta=w^{n}e^{2(w-1)\tilde\rho},
\end{split}\end{equation*}
so
\begin{equation}\label{eq:T_b-theta}
T_{b,\theta}=e^{-\tilde\rho}(w^{-2}\Delta_{S_b}+|\tilde\rho|^2)e^{\tilde\rho}.
\end{equation}

Now, since $\Delta_\theta$ is radially elliptic, the domain of 
$\Delta_{\rad,\theta}$ is $\Hso^2(M)^K$. For any 
$f\in \Hso^2(M)^K$, the map $\theta\mapsto\Delta_{\rad,\theta} f\in 
L^2(M,dg)$ is strongly analytic, and this is what it means for 
$\Delta_{\rad,\theta}$ to be an analytic family of type A.
\end{proof}

\subsection{Analytic vectors}
A general abstract theorem due to Nelson, cf.\ \cite[Volume~2]{Reed-Simon}, 
uses the functional calculus to construct a dense set of analytic vectors for 
the generator of a group of unitary operators. We shall instead define an
explicit subspace of analytic vectors $\AAA$, which is meant to demonstrate 
the essentially elementary nature of this result in our context. We ultimately 
wish to employ the operators $\Delta_{\rad,\theta}$ for $\theta \in D=
\{\theta:\ |\im\theta| <\frac{\pi}{2}\}$, and using Nelson's theorem we could do 
this directly. A slight disadvantage with our more concrete approach is that 
this must be done in two steps now, first letting $\theta \in D'=\{\theta:\ 
|\im\theta|<\frac{\pi}{4}\}$, and then extending to $\theta \in D$, but only
a minor extra argument is needed for this.

The action of the Weyl group $W$ extends naturally to $\fraka_{\Cx}$. 
Define $\AAA$ to be the space of restrictions to $\fraka$ of entire functions 
$f$ on $\fraka_{\Cx}$ which are $W$-invariant and which decay faster than any 
power of $e^{-|z|}$ in every cone $\{z\in \fraka_\Cx:\ |\im z|\leq 
C|\re z|\}$, $0<C<1$. In other words, denoting both the entire function and
its restriction to $\fraka$ by $f$, we have $f \in \AAA$ if, for every 
$0<C<1$ and $N>0$, 
\begin{equation*}
\sup_{|\im z|\leq C|\re z|} |f(z)|e^{N|z|}<+\infty.
\end{equation*}
Clearly, for any $\theta \in D'$ and $f \in \AAA$, $U_\theta f$
is rapidly decreasing on $\fraka$.  

\begin{prop}\label{prop:dense}
For $\theta \in D'$, i.e.\ $|\im\theta\,|<\frac{\pi}{4}$, $U_{\theta}\AAA$ 
is dense in $L^2(\fraka)^{W}$.
\end{prop}
\begin{proof} 
Since $\calC^0_c(\fraka)^W$ is dense in $L^2(\fraka)^W$ (with respect to
the singular measure $dg_0 = \eta\, dx$ on $\fraka$ -- in this
proof we use $x$ for points in $\fraka$), it suffices 
to show that any $f \in \calC^0_c(\fraka)^W$ can be approximated by
functions $f_t \in \AAA$. To this end, set 
\begin{equation*}
f_t(x)= c_n t^{-n/2} \int f(y) e^{-|x-y|^2/t}\,dy,
\end{equation*}
where $n = \dim \fraka$ and $c_n$ is chosen so that $\int f_t(x)\, dx = \int f(x)\, dx$
for all $t>0$, i.e. so that $c_n t^{-n/2}e^{-|x|^2/t}$ is the Euclidean heat kernel. 
We claim first that $f_t\in\AAA$ when $t > 0$. Indeed, $f_t(x)$ is the restriction 
to $\fraka$ of $f_t(z) = \int c_n t^{-n/2} e^{-(z-y)^2/t}f(y)\, dy$ and $\exp(-(z-y)^2)$ is 
entire in $z$ and decreases faster than any power of $e^{-|z|}$ in $|\im z| 
< C|\re z|$ whenever $C < 1$, and this decay is preserved even after the integration over a 
compact set in $y$. Moreover, the action of $W$ is by Euclidean isometries and 
hence commutes with the heat kernel, so each $f_t(x)$ is $W$-invariant. This proves 
the claim. 

Now let us show that $U_\theta\AAA$ is dense in $L^2(\fraka)^W$ when $\theta \in D'$.
For the case $\theta = 0$, note that for $f \in \calC^0_c(\fraka)^W$, 
$e^{|x|^2}f_t$ is uniformly bounded when $t<1$, and 
$\sup e^{|x|^2}|f(x)-f_t(x)|\to 0$ as $t\to 0$. Since $e^{-|x|^2}\in 
L^2(\fraka;dg_0)^W$, we have $f_t\to f$ in this space. In the general
case, for any $\theta \in D'$, define
\begin{equation*}
\tilde{f}_t(x)= c_n e^{n\theta} t^{-n/2}\int f(y) e^{-e^{2\theta}|x- y|^2/t}\,dy. 
\end{equation*}
We must show that $\tilde{f}_t \to f$ in $L^2(\fraka)^W$ and $f_t \in U_\theta \AAA$.
For the former, note that $\tilde{f}_t(x)$ is just the function $f_t(x)$
analytically continued to complex time $\tau = e^{-\theta}t$, and the
same proof as above shows that $f_\tau \to f$ in $L^2$. Finally,
\[
U_{-\theta}\tilde{f}_t(x) = c_n e^{n\theta/2}t^{-n/2}\int f(y)e^{-|x-e^\theta y|^2/t}\,dy
\]
and as in the first part of the proof, this is certainly in $\AAA$.
\end{proof}

\begin{cor}\label{cor:dense}
For $|\im\theta\,|<\frac{\pi}{4}$, $U_{\theta}\AAA$ is dense in $\Hso^s(M)^K$
for any $s \geq 0$. 
\end{cor}
\begin{proof}
Implicit in the definition of these Sobolev spaces, i.e.\ using radial
ellipticity and the positivity of the Laplacian, cf.\ 
\cite{Mazzeo-Vasy:Sl3} for an explanation, 
\[
(\Delta_\rad +1)^{s/2}:\Hso^s(M)^K\to \Hso^0(M)^K\equiv L^2(M,dg)^K
\equiv L^2(\fraka,dg_0)^W
\]
is an isomorphism. Thus, $f_t\to f$ as $t\to 0$ in $\Hso^s(M)^K$ if and only 
if $(\Delta+1)^{s/2}f_t\to (\Delta+1)^{s/2}f$ in $L^2(\fraka,dg_0)^W$. 
So given $f\in \Hso^s(M)^K$, let $k=(\Delta+1)^{s/2}f$. Since $\AAA$ is 
dense in $L^2(\fraka;dg_0)^W$, there exists a family $k_t\in\AAA$ with $k_t\to k$ 
as $t\to 0$ in $L^2(\fraka;dg_0)^W$. Now let $f_t=(\Delta+1)^{-s/2}k_t$ and note 
that $f_t\in\AAA$. Thus, $f_t\to f$ in $\Hso^s(M)^K$ as desired.
\end{proof}

For functions or distributions $k$ which do not lie in $\AAA$, 
$U_\theta k$ may still have a continuation. For example, 
if $k=\delta_o$, the delta distribution at $o$, then using its homogeneity
we see that for $\theta$ real, $U_{\theta}\delta_o=(\det D_o\Phi_{\theta})^{-1/2}\delta_o$. 
Hence $U_{\theta}\delta_o$ extends to be analytic in $\theta$ (e.g.\ with values 
in some Sobolev space of sufficiently negative order), and so the Green function, 
$R(\ev)\delta_o$ also extends via $\langle f,R(\ev)\delta_o\rangle$ for $f\in\AAA$.

\subsection{The domain of continuation}
We now describe the Riemann surface $\widetilde{\calY}_{\pi/2}$ to which $R(\ev)$ continues.
We expect that $\widetilde{\calY}_{\pi/2}$
should be very simple, specifically either $\Cx$ or the Riemann
surface for $\sqrt{z}$ or, at worst, for $\log z$, and in
particular should be ramified at only one point. However, we only 
consider the continuation up to angle $\pi$ ($\im\theta = \pm \pi/2$),
and in particular omit the ray where $\ev$ makes an angle of $\pm \pi$
with the spectral axis, and on which it is known that there exist poles of $R(\ev)$
in many cases (e.g.\ on even dimensional hyperbolic spaces). 

In addition, the $N$-body methods by themselves
cannot rule out the existence of other
poles in the nonphysical half-plane of $\sqrt{z}$. 
These poles are more serious than they 
might seem at first because in the inductive scheme, poles for
the resolvent on spaces of rank less than $n$ give rise to ramification
points in the continuation for spaces of rank $n$. In the present paper
we only describe the `worst case scenario', and allow for the existence
of these poles. We expect that the precise analysis of $\im\theta\to\pm\pi/2$
will exclude their existence, see the discussion at the end of the
last section.

Recall the symmetric space of lower rank, $\Sigma^b$, associated to $S_b$, $b \in 
I\setminus \{*\}$. Denote by ${\mathcal P}_{b,\theta}$ the pure point spectrum of 
$\Delta_{\Sigma^b,\rad,\theta}$, and also assume that the set ${\mathcal T}_{b,\theta}$ 
of thresholds for $\Delta_{b,\rad,\theta}$ has been defined inductively. Now define
the set of thresholds for $\Delta_{\rad,\theta}$, ${\mathcal T}_\theta$, by
\begin{equation*}
{\mathcal T}_\theta = \bigcup_{b\neq *}\{|\rho-\rho_b|^2+ \gamma:\ \gamma
\in {\mathcal P}_{b,\theta} \cup {\mathcal T}_{b,\theta}\}.
\end{equation*}
Note that for $b=0$, $\Sigma^b$ is a point, and so $\rho_b = 0$ and 
${\mathcal P}_{0,\theta}=\{0\}$ for all $\theta$; this means that we always have 
$|\rho|^2 \in {\mathcal T}_\theta$ for any $\theta$. In addition,
since $\rho-\rho_b \in S_b$ and $\rho_b \in S^b$ are orthogonal, this again contributes 
the value $|\rho-\rho_b|^2 + |\rho_b|^2 = |\rho|^2$ to ${\mathcal T}_\theta$. 
Presumably, ${\mathcal T}_\theta$ consists of the single element  $|\rho|^2$, 
but this would rely on knowing that all spaces of rank less than $n$ have
no point spectrum and no thresholds except at $|\rho_b|^2$; in any case, this is 
true when $n=2$. 

We shall prove later, in Theorem~\ref{thm:ess-spec}, that 
as an operator on $L^2(\fraka;dg_0)^W$,
\begin{equation}\label{eq:ess-spec-8}
\spec_{\mathrm{ess}}(\Delta_{\rad,\theta})=\{\gamma+e^{-2i\im\theta}[0,+\infty):\ 
\gamma\in {\mathcal T}(\theta)\}
\end{equation}
when $|\im\theta|<\pi/2$. In other words, every eigenvalue and threshold of the 
scaled radial Laplacian of each subsystem $\Sigma^b$ contributes a ray 
to the essential spectrum of $\Delta_{\rad,\theta}$ making an angle $-2\im\theta$ 
with the positive real axis and emanating from that point. This ray
is, in fact, the essential spectrum of the scaled `tangential operator'
$T_{b,\theta}=U_\theta^{-1}(\Delta_{S_b}+2H_{\rho-\rho_b})U_\theta$.
Granting this result, we now proceed with the rest of 
the complex scaling argument.

Normalize so that $\arg(z)\in(-2\pi,0)$ for $z\in\Cx\setminus[0,+\infty)$,
and let $\sqrt{z}$ be the branch of the square root function with
$\im\sqrt{z}<0$ on $\Cx\setminus[0,+\infty)$. Let $S$ be the Riemann
surface of $\sqrt{\ev-\ev_0}$, with the ray with
$\arg\sqrt{\ev-\ev_0}=\frac{\pi}{2}$ removed. The map
\begin{equation*}
F:S\ni z=\sqrt{\ev-\ev_0}
\mapsto \ev=z^2+\ev_0
\end{equation*}
gives a double cover of $\Cx\setminus(-\infty,\ev_0]$;
the ray $(-\infty,\ev_0)$ is only covered once. We call the
part $S_0$ of $S$ with $\im\sqrt{\ev-\ev_0}<0$, i.e.\ $\arg\sqrt{\ev-\ev_0}
\in(-\pi,0)$, the `physical half-plane'.

We define Riemann surfaces $\calY_{\beta}$, $\beta\in[0,\pi/2]$,
that are open subsets of $S$ and
such that $S_0\subset\calY_{\beta}$.
The part $S_-$ of $S$ with $\arg\sqrt{\ev-\ev_0}\in(-\pi/2,\pi/2)$
can be identified with $\Cx\setminus (-\infty,\ev_0]$ via $F$.
Then by definition
\begin{equation}\begin{split}\label{eq:calY-beta-def}
\calY_{\beta}\cap S_-\equiv& \{\ev\in\Cx:\ \arg\sqrt{\ev-\ev_0}
\in (-\pi/2,\beta)\}\\
&\qquad\setminus 
\{\gamma+ e^{2i\beta}[0,+\infty):\ \gamma\in
{\mathcal T}(\beta)\},\ \beta\in[0,\pi/2).
\end{split}\end{equation}
Note that $\{\gamma+ e^{2i\beta}[0,+\infty):\ \gamma\in
{\mathcal T}(\beta)\}$ is exactly the right hand side of
\eqref{eq:ess-spec-8} if we let $\im\theta=-\beta$.
With $S_+$ denoting the part of $S$ with
$\arg\sqrt{\ev-\ev_0}\in(-3\pi/2,-\pi/2)$, we define
\begin{equation*}\begin{split}
\calY_{\beta}\cap S_+\equiv& \{\ev\in\Cx:\ \arg\sqrt{\ev-\ev_0}
\in (-\pi-\beta,-\pi/2)\}\\
&\qquad\setminus 
\{\gamma+e^{-2i\beta}[0,+\infty):\ \gamma\in
{\mathcal T}(\beta)\},\ \beta\in[0,\pi/2).
\end{split}\end{equation*}
Note that with this definition, $\calY_0$ is the `physical half plane' $S_0$.

\begin{rem}\label{rem:identify}
Although each $\calY_\beta$ can be considered as a subset of $S$, it is
important to realize that even in the overlap of these regions for different
values of $\beta$, the $\calY_\beta$ should not be identified
with each other. Rather, two points $p\in\calY_\beta$
and $q\in\calY_\gamma$ with $\gamma\leq\beta$
with the same image $\ev'$ in $S_-$, say,
should only be identified if
\begin{equation*}
\ev'\nin\{\gamma+ e^{2i\theta}[0,+\infty):\ \gamma\in
{\mathcal T}(\theta),\ \theta\in[\gamma,\beta]\}.
\end{equation*}
An equivalent formulation would be that the two points should be
identified if
there is a path
in $S_-$ connecting $\ev'$ to `physical region'
$\arg\sqrt{\ev-\ev_0}\in(-\pi/2,0)$ which stays entirely in the intersection
of $S_-\cap\calY_\beta$ and $S_-\cap\calY_\gamma$.
\end{rem}

For this reason we make the following definition.

\begin{Def}\label{Def:calY-pi-2}
For $\beta\in(0,\pi/2]$,
we define $\widetilde{\calY}_\beta$ as the disjoint union
of $\calY_\gamma$, $\gamma\in[0,\beta)$, modulo the equivalence
relation described above. We define the topology of $\widetilde{\calY}_\beta$
by requiring that open subsets of $\calY_\gamma$ to be open in
$\widetilde{\calY}_\beta$, and taking these as a base for the topology
of $\widetilde{\calY}_\beta$ as $\gamma$ runs over $[0,\beta)$.
Letting the $\calY_\gamma$ be coordinate charts, we make
$\widetilde{\calY}_\beta$ into a Riemann surface.
\end{Def}

\begin{rem}
In this definition, if $\beta<\frac{\pi}{2}$,
we could replace $\gamma\in[0,\beta)$ by
$\gamma\in [0,\beta]$; the resulting Riemann surface would be the same.
\end{rem}

Denote by $R(\ev,\theta)$ the operator $(\Delta_{\rad,\theta} - \ev)^{-1}$.
To be definite, we consider only the analytic continuation of $R(\ev) = 
R(\ev,0)$ from the lower right quadrant $\im(\ev-\ev_0)<0$ through the ray 
$(\ev_0,+\infty)$), i.e.\ to $S_-\cap\calY_\beta$;
the continuation from $\im(\ev-\ev_0)>0$ is handled 
nearly identically. 

The main point, roughly speaking, is that when $-\frac{\pi}{2}<\im\theta<0$, 
$\Delta_{\theta}-\ev$ is a holomorphic family of operators (in $\ev$)
with values in the space of radially elliptic operators on $M$.
Thus $R(\ev,\theta)$ is meromorphic in $\ev$ outside $\mbox{spec}_{\mathrm{ess}}
(\Delta_{\rad,\theta})$ with values in bounded operators on
$L^2(\fraka;dg_0)^W$. 
This family has only finite rank poles, and these are the poles of the continuation of 
$R(\ev)_{\mathrm{rad}}$ in $\calY_\beta\cap S_-$ if we choose $\theta$ so that $\beta=-\im\theta<\frac{\pi}{2}$.

\subsection{Continuation of the resolvent}
We finally indicate the proof of the analytic continuation of the resolvent, 
which is simply an application of the theorem of Aguilar-Balslev-Combes in our setting. 
\begin{thm*} {\rm (\,}\cite[Theorem~16.4]{Hislop-Sigal:Spectral}\,{\rm )}
Suppose that $U_\theta$ and $\AAA$ satisfy the hypotheses (i)-(iii) listed
in the beginning of \S 4, and that $\Delta_\theta$ is a type-A analytic 
family in the strip $D'=\{\theta:\ |\im\theta|<\frac{\pi}{4}\}$, and
\eqref{eq:ess-spec-8} holds for $\theta\in D$. Then
\begin{enumerate}
\item
For $f,g\in\AAA$, $\beta<\frac{\pi}{4}$, the function
$\langle f,R(\ev)_{\mathrm{rad}}g\rangle$ has a meromorphic 
continuation to $\calY_{\beta}$.
\item
The poles of the continuation of $\langle f,R(\ev)g\rangle$ to $\calY_\beta$,
$\beta<\frac{\pi}{4}$, are the eigenvalues of $\Delta_{\rad,\beta}$.
\item
The poles are independent of the choice of $U_\theta$ in the sense
that if $U'_\theta$ and $\AAA'$ also satisfy (i)-(iii) and if
$\AAA\cap\AAA'$ is dense, then the eigenvalues of $U'_\theta
\Delta_{\rad} (U'_\theta)^{-1}$ are the same as those of 
$\Delta_{\rad,\theta}$.
\end{enumerate}
\end{thm*}

All of the hypotheses have already been discussed and verified. We shall
briefly outline the proof of the first part since the idea is simple.
To relate $R(\ev,\theta)$ and $R(\ev)$, fix $\ep>0$, and suppose that 
\begin{equation*}
\theta \in \Omega_{\ep}=\{-\ep<\im\theta<\frac{\pi}{4}\}
\qquad \mbox{and}\qquad \arg(\ev-\ev_0)\in(-\pi,-\ep).
\end{equation*}
When $\theta$ is real, $U_{\theta}$ is unitary and so
\begin{equation}
\langle f,R(\ev)g\rangle
=\langle U_{\bar\theta}f,(U_{\theta}R(\ev)U_{\theta}^{-1})
U_{\theta}g\rangle =\langle U_{\bar\theta}f,R(\ev,\theta)U_{\theta}g\rangle
\label{eq:anal-cont-16}
\end{equation}
since $U_{\theta}R(\ev)U_{\theta}^{-1}=R(\ev,\theta)$. The left side of 
this equation is independent of $\theta$, while the expression on the (far) 
right is analytic in $\theta$ on $\Omega_{\ep}$, hence is also constant on 
this domain. This holds when $\arg(\ev-\ev_0)\in(-\pi,-\ep)$. 

To extend $\langle f,R(\ev)g\rangle$ to $\calY_\beta$, take $\theta$ with 
$\im\theta=-\beta$. Then for $\ev\in\Cx$ with $\im(\ev-\ev_0)<0$,
$\langle f,R(\ev)g\rangle$ is given by the right hand side of
\eqref{eq:anal-cont-16}. But this right side is analytic in $\ev$
away from the spectrum of $\Delta_{\rad,\theta}$, and meromorphic away from 
its essential spectrum, hence is meromorphic on $\calY_\beta$, as claimed.

This continuation is clearly independent of the choice of $\theta$ with
$-\im\theta=\beta$ since any such continuation is
a meromorphic function of $\ev$ that agrees with a given function on
an open set. In addition, the continuation is independent of $\beta$
in the sense that if $p\in\calY_\beta$ and $q\in\calY_\gamma$ are
identified in the sense of Remark~\ref{rem:identify}, so there
is a path connecting them to the physical region that does not
intersect the cuts in either $\calY_\beta$ or in $\calY_\gamma$,
then $\langle f,R(\ev)g\rangle$ is the same whether the $\beta$ or
$\gamma$ is used to define it.

Note that this does not yet quite say that $R(\ev)\delta_o$ continues
as a distribution, since that would require that the right hand side of 
\eqref{eq:anal-cont-16} be defined for any $f\in\Cinf_c(\fraka)^W$,
while for most $f$, $U_{\theta} f$ does {\em not} have an analytic extension 
from the real axis. This is where we require the deformed group of unitary 
operators, $U_{\theta,T}$, defined in (\ref{eq:defuo}). 
Recall that the associated diffeomorphisms $\Phi_{\theta,T}$ fixes $B_T(o)$
pointwise and equals $\Phi_{\theta}$ when $|a|$ is sufficiently large.
We use precisely the same arguments as above to establish the density of 
$U_{\theta,T}\AAA$. Hence by the 
uniqueness part of the Aguilar-Balslev-Combes theorem, the induced analytic 
extensions agree with one another no matter the value of $T$, and also agree
with the extension associated to $U_\theta$. But if $f\in\Cinf_c(B_T(o))^W$,
then $U_{\theta,T}f=f$ and so $U_{\theta,T}f=f$ has an analytic extension
to $\theta\in\Cx$. Arguing as before, the formula 
\begin{equation}\label{eq:anal-cont-32}
\langle f,R(\ev)\delta_o\rangle 
=\langle U_{\bar\theta,T}f,R(\ev,\theta,T)U_{\theta,T}\delta_o\rangle
=\langle f,R(\ev,\theta,T)\delta_o\rangle
\end{equation}
shows that $R(\ev)\delta_o$ does indeed extend analytically 
as a distribution to $\calY_{\beta}$, $\beta\in(0,\frac{\pi}{4})$,
since the right hand side has this
property.

Although we have only constructed a subset $\AAA \subset L^2(\fraka;dg_0)^W$
for which $U_{\theta}\AAA$ is dense in $L^2(\fraka;dg_0)^W$ when
$|\im\theta|<\pi/4$, we can still continue $R(\ev)$ to 
$\widetilde{\calY}_{\pi/2}$,
for which the formula \eqref{eq:anal-cont-16} requires larger $\im\theta$.

\begin{thm*}[Theorem~\ref{thm:main}]
The Green function $G_o(\ev)$ continues meromorphically to 
$\widetilde{\calY}_{\pi/2}$ as a distribution.
\end{thm*}

\begin{proof}
We have shown that the hypotheses of the Aguilar-Balslev-Combes theorem are
satisfied for $D'=\{\theta:\ |\im\theta|<\frac{\pi}{4}\}$ (for either
$U_\theta$ or $U_{\theta,T}$) (except for the proof of \eqref{eq:ess-spec-8}).
Hence $R(\ev)$ continues meromorphically to $\calY_{\beta}$ , $\beta\in
(0,\pi/4)$, in the
precise sense of the theorem. In particular, $G_o(\ev)$ continues 
meromorphically to $\calY_{\beta}$ as a distribution. However, at first we
ignore the continuation itself, i.e.\ restrict to $\ev$ with
$\arg\sqrt{\ev-\ev_0}\in(-\pi/2,0)$, and extend the scaling argument instead.

Namely, we use the semigroup property $U_\theta U_{\theta'}=
U_{\theta+\theta'}$, which implies the analogue of \eqref{eq:anal-cont-16}:
\begin{equation}\label{eq:anal-cont-24}
\langle f,R(\ev,\theta')g\rangle=\langle U_{\bar\theta}f,R(\ev,\theta+\theta')
U_\theta g\rangle
\end{equation}
for $f,g\in\AAA$, $|\im\theta|<\frac{\pi}{4}$,
$\arg\sqrt{\ev-\ev_0}\in(-\pi/2,0)$. Hence $U_\theta R(\ev,\theta') 
U_\theta^{-1}=R(\ev,\theta+\theta')$ for $\theta\in\RR$,
and so \eqref{eq:anal-cont-24} gives the
continuation of $R(\ev,\theta')$ to $\ev\in\calY_{-\im\theta'-\im\theta}$. For
$\beta\in(0,\pi/2)$, we may take $\theta,\theta'$ with
$\im\theta=\im\theta'=-\beta/2$, so we conclude that
$R(\ev)$ continues analytically to $\calY_\beta$.

This also gives the extension of $R(\ev)\delta_o$ to $\calY_\beta$ as
a distribution. Indeed, this extension exists in ${\mathcal D}'(B_T(o))$ 
for any $T > 0$, and the density of $\AAA$ implies that these extensions 
are all the same.

Finally, by the very definition of $\widetilde{\calY}_{\pi/2}$, the
analytic continuation of $G_o(\ev)=R(\ev)\delta_o$ to $\calY_\beta$ for every
$\beta\in(0,\pi/2)$ gives the desired analytic continuation to
$\widetilde{\calY}_{\pi/2}$.
\end{proof}

\begin{rem}
We emphasize that although the analytic extension to $\calY_\beta$, $\beta
\in[\pi/4,\pi/2)$ is defined in two steps, the analytic extension of
$\delta_o$ as a distribution on $B_T(o)$ can be done at once. Indeed,
both $U_{\theta,T}\delta_o$ and $U_{\theta,T}f$, $f\in\Cinf_c(B_T(o))$,
have an analytic extension
to $\{\theta:\ |\im\theta|<\pi/2\}$, so \eqref{eq:anal-cont-32}
defines the extension (in $\dist(B_T(o))$)
of $R(\ev)\delta_o$ directly in the region $\calY_\beta$, $\beta\in(0,\pi/2)$.
\end{rem}

\section{The parametrix construction}
Our final goal is to identify the essential spectrum of $\Delta_{\rad,\theta}$
when $|\im \theta| < \pi/2$. As usual, the strategy is to construct a 
parametrix for the scaled resolvent $(\Delta_{\rad,\theta}-\ev)^{-1}$ with 
compact remainder when $\ev$ is outside the putative essential spectrum. 
We shall approach this in a series of steps. The procedure is inductive,
and the parametrix is built up from the resolvents of the scaled model 
operators $L_{b,\theta} = T_{b,\theta} + \Delta_{\Sigma^b,\rad,\theta}$, 
$b \in I$, localized to neighbourhoods of $S_b \times \{0\} \subset S_b 
\times S^b$ (for $b=*$, $L_{b,\theta}=\Delta_{\rad,\theta}$ and we localize 
to a compact neighbourhood of $0 \in \fraka$). In the first step, we use the `softest' form 
of this induction, employing only radial ellipticity, to obtain an exact 
inverse to $\Delta_{\rad,\theta}-\ev$ when $\ev$ is sufficiently large 
and lies outside any small cone surrounding the essential spectrum. 
We also obtain decay estimates for the norm of the resolvent as $|\ev| \to 
\infty$. The point is that we are able to get a parametrix with remainder
which has small norm, which can then be inverted away using a Neumann
series. This involves the use either of the associated semiclassical calculus
or, perhaps more familiarly, a pseudodifferential calculus with spectral parameter, 
as described for example in \cite{Shubin}; see also \cite{Vasy-Wang:Smoothness}
where this is used in the $N$-body setting. These decay estimates are 
necessary in the next step, where we use the convolution formula for the 
resolvent on a product space from \cite{Mazzeo-Vasy:Resolvents} to describe 
the resolvents $(L_{b,\theta}-\ev)^{-1}$ in terms of the resolvents for 
$T_{b,\theta}$ and $\Delta_{\Sigma^b,\rad,\theta}$; here we use the
induction hypothesis, specifically the estimates from the first step,
for the latter factor. A slight technical twist is that we need to modify 
this formula to handle sums of nonselfadjoint operators. This would follow
from a more general abstract theorem (Ichinose's lemma), but we also
indicate a direct proof. In the third and final step we use the resolvents
of the model operators obtained in the previous step to obtain
a parametrix for $(\Delta_{\rad,\theta}-\ev)^{-1}$ with a {\it compact}
remainder, for all $\ev$ outside the essential spectrum. After this
we can finish the whole construction by applying the analytic Fredholm
theorem. 

\medskip
\noindent{\bf Step 1: The parametrix for large spectral parameter}

\noindent As described above, the first task is to construct and obtain estimates 
on the resolvent $(\Delta_{\rad,\theta}-\ev)^{-1}$ when $\ev$ tends to 
infinity and remains outside some sector. More precisely, we show that for 
any $\ep>0$, and $R = R_\ep > 0$ sufficiently large, depending on $\ep$, 
\begin{equation*}
\spec(\Delta_{\rad,\theta}) \cap \{|\ev| > R\} \subset e^{-2i[\im\theta-\ep,
\im\theta+\ep]}[0,+\infty) \cap \{|\ev| > R\} := D_{R,\ep}^c,
\end{equation*}
and for $\ev$ large and outside this latter set we estimate the norm of 
$(\Delta_{\rad,\theta}-\ev)^{-1}$ on $L^2(M)^K$ in terms of 
powers of $1/|\ev|$. This is proved by constructing a parametrix with error 
term which tends to zero in operator norm as $\ev \to \infty$, 
and which then be inverted away. This step is `soft' inasmuch as we only use 
radial ellipticity in this argument, but we emphasize that this error term is
small, but not necessarily compact.

One could proceed rather abstractly at this stage by showing that $\Delta_{\rad,\theta}$
is m-sectorial, cf.\ \cite[Volume~II, Section~VIII.6]{Reed-Simon}. This would 
involve considering the quadratic form $\langle \phi, \Delta_{\rad,\theta} \phi\rangle$ 
for $\phi\in\Cinf_c(M)^K$. The point here is that the difference between 
$\Delta_{\fraka,\theta}$ and $\Delta_{\rad,\theta}$ is a first order differential 
operator, and the form corresponding to this difference can be estimated via 
Cauchy-Schwartz. However, the fact that we must use a nontrivial measure on
$\fraka$ because of the identification $L^2(M)^K \cong L^2(\fraka,dg)^W$
makes this not entirely trivial. 

However, in keeping with the other steps, we construct the parametrix by piecing 
together the simplest of parametrices for the model operators $L_{b,\theta}$ 
using a $(W,\ofraka)$-adapted partition of unity, maintaining control on 
large $\ev$ behaviour. 

\begin{prop}\label{prop:res-est-infty}
For any $\ep>0$ there exist $R,C>0$ such that when $|\ev|>R$ and
$|\arg\ev + 2\im\theta|>\ep$, we have 
\begin{equation*}\begin{split}
&R(\ev,\theta)=(\Delta_{\rad,\theta}-\ev)^{-1}\in\bop(L^2(M)^K),\\
&\|R(\ev,\theta)\|_{\bop(L^2(M)^K)}\leq \frac{C}{|\ev|}.
\end{split}\end{equation*}
\end{prop}

\begin{proof}
Recall that, for any $b \in I$, $L_{b,\theta}-\ev = T_{b,\theta} + 
\Delta_{b,\rad,\theta}-\ev$ is an operator on $S_b \times \Sigma^b$ which is 
constant coefficient on the first factor and radial on the second; moreover, we 
are only interested in its restriction to a fixed bounded neighbourhood in $\Sigma^b$. 
For $\ev$ outside this sector, this is an elliptic element of the pseudodifferential calculus 
with large spectral parameter (satisfying uniform estimates in the $S_b$ factor), 
as defined in \cite{Shubin}. Choose two different sets of 
cutoffs, $\{\phi_b\}$ and $\{\psi_b\}$, $b \in I$, each satisfying (i)-(iii) of 
Definition~\ref{def:W-aol-adapted}, and such that $\psi_b$ is identically $1$ 
on a neighbourhood of $\supp\phi_b$ and $\supp\psi_b$ is sufficiently
close to $\overline{S_b}$; the smallness of the
support ensures that $\Delta_{\Sigma^b,\theta}$ is elliptic on it.
There exists a parametrix in this calculus, 
$G_{b,\theta}(\ev)$, which we may as well assume is $K^b$-invariant
(by averaging it over $K^b$), which is supported
near $\supp\psi_b$. This satisfies the analogues of the bounds in the statement of this 
proposition, and in addition,
\[
(L_{b,\theta} - \ev)G_{b,\theta}(\ev)\phi_b = \phi_b + F_{b,\theta}(\ev),
\]
where $\|F_{b,\theta}(\ev)\|_{\bop(L^2(M)^K)}\leq C_{N,\ep}/|\ev|^N$ for
any $N,\ep > 0$, by virtue of the properties of residual elements in this large
parameter calculus. Finally, define 
\begin{equation*}
G_\theta(\ev)=\sum_b \psi_b G_{b,\theta}(\ev)\phi_b.
\end{equation*}

We have 
\begin{equation*}
(\Delta_{\rad,\theta}-\ev)G_\theta(\ev)=
\Id+\sum_b \left([\Delta_{\rad,\theta},\psi_b] G_{b,\theta}(\ev)\phi_b +
\psi_b F_{b,\theta}(\ev)\right) = \Id+F_\theta(\ev).
\end{equation*}
Since $\supp \, [\Delta_{\rad,\theta},\psi_b]$ is disjoint from $\supp\phi_b$, this
error term also satisfies 
\[
\|F_\theta(\ev)\|_{\bop(L^2(M)^K)}\leq \frac{C_N}{|\ev|^N}
\]
for any $N,\ep > 0$. Thus $\Id + F_\theta(\ev)$ is invertible when $|\ev| > R$
(still outside this sector), so
\begin{equation*}
(\Delta_{\rad,\theta}-\ev)G_\theta(\ev)(\Id+F_\theta(\ev))^{-1}=\Id,
\end{equation*}
and standard arguments also show that this is a left inverse too. This
means that 
\[
(\Delta_{\rad,\theta}-\ev)^{-1} = R(\ev,\theta)=G_\theta(\ev)(\Id+F_\theta(\ev))^{-1}.
\]
The estimates for $R(\ev,\theta)$ follow directly from those for $G_{b,\theta}(\ev)$.
\end{proof}

\medskip

\noindent{\bf Step 2: Resolvents of the model operators} 

\noindent We now use the convolution formula from \cite{Mazzeo-Vasy:Resolvents}
and the decay estimates obtained in the previous step to express the resolvent 
for each model operator 
\begin{equation}
L_{b,\theta}=T_{b,\theta}+\Delta_{\Sigma^b,\theta}
\label{eq:declb}
\end{equation}
in terms of the resolvents of the two summands.  We assume now that $b \neq *$,
since the analysis of $L_{*,\theta} = \Delta_{\rad,\theta}$ is what we are 
ultimately trying to understand. Note also the other extreme case $b=0$,
where $L_{0,\theta} = (\Delta_\fraka)_\theta = e^{-2\theta}\Delta_\fraka$. 

The first summand is a constant coefficient operator on $S_b$
which is the rescaling of 
\begin{equation*}
T_b=\Delta_{S_b}+2(H_\rho-H_{\rho_b}).
\end{equation*}
Recall that if $M_f$ is the operator of multiplication by a function $f>0$, then 
\begin{equation*}
M_f:L^2(S_b,f^2\,da_b)\to L^2(S_b,da_b)
\end{equation*}
is a unitary isomorphism. Thus choosing $f=e^{\rho-\rho_b}$, then we see 
that $T_b$ acting on 
$L^2(S_b,e^{2(\rho-\rho_b)}\,da_b)$ is unitarily equivalent to
\begin{equation}
\tilde T_b=f^{-1}(\Delta_{S_b}+2H_{\rho-\rho_b})f 
=\Delta_{S_b}+(\rho-\rho_b)\cdot(\rho-\rho_b) = \Delta_{S_b} + |\rho-\rho_b|^2,
\end{equation}
acting on $L^2(S_b,da_b)$, and correspondingly, using the same $f$,
see \eqref{eq:T_b-theta},
$T_{b,\theta}$ is unitarily equivalent to
\begin{equation*}
\tilde T_{b,\theta}=\Delta_{S_b,\theta}+|\rho-\rho_b|^2,
\end{equation*}
also on $L^2(S_b,da_b)$. In particular, since 
$\Delta_{S_b,\theta} = e^{-2\theta}\Delta_{S_b}$, it follows immediately that
\begin{equation}
\spec(T_{b,\theta})=|\rho-\rho_b|^2+e^{-2i\im\theta}[0,+\infty).
\label{eq:spectb}
\end{equation}
In addition, from the Fourier transform representation of this operator we deduce
that 
\begin{equation}
\|(T_{b,\theta}-\ev)^{-1}\| \leq C/|\ev|
\label{eq:desttb}
\end{equation}
as $\ev \to \infty$ away from $D_{R,\ep}^c$. 

Since the rank of $\Sigma^b$ is strictly less than $n$, the spectrum of the other
summand in (\ref{eq:declb}) is understood by induction. 
Because these rescaled operators are not self-adjoint, it is not completely trivial
that the spectrum of $L_{b,\theta}$ is the sum of spectra of the two operators on 
the right. This follows from an abstract lemma due to Ichinose, cf.\ \cite[Volume~IV, 
Section~XIII.9, Corollary~2]{Reed-Simon}, but also follows directly from the
existence of the resolvent when $\ev$ is outside the sum of these two spectra:

\begin{cor}\label{cor:tensor-product}
For any $b \in I\setminus\{*\}$, as an operator on 
$L^2(\Sigma^b\times S_b,e^{2(\rho-\rho_b)}\,da_b\,dg_b))$, 
\begin{equation}\label{eq:tens-prod-spec}
\spec(L_{b,\theta})=\{\ev'+\ev'':\ \ev'\in\spec(\Delta_{\Sigma^b,\theta}),
\ \ev''\in |\rho-\rho_b|^2+e^{-2i\im\theta}[0,+\infty)\}.
\end{equation}
In particular, outside this set, 
\begin{equation*}
R_{b,\theta}(\ev)=(L_{b,\theta}-\ev)^{-1}\in \bop
(L^2(\Sigma^b\times S_b,e^{2(\rho-\rho_b)}\,da_b\,dg_b)).
\end{equation*}
\end{cor}

\begin{proof}
The convolution formula states that
\begin{equation}\label{eq:sa-product-formula}
R_{b,\theta}(\ev)=\frac{1}{2\pi i}\int_\gamma (\Delta_{\Sigma^b,\theta}-\mu)^{-1}
\otimes (T_{b,\theta}-(\ev-\mu))^{-1}\,d\mu,
\end{equation}
where $\gamma$ is a path in $\Cx$ which avoids $\spec(\Delta_{\Sigma^b,\theta})$ 
and $\ev - \spec(T_{b,\theta})$, and which diverges linearly from these rays.
The decay estimates 
\begin{equation*}
\|(\Delta_{\Sigma^b,\theta}-\mu)^{-1}\|\leq |\im\mu|^{-1},
\ \|(T_{b,\theta}-(\ev-\mu))^{-1}\|\leq |\im(\ev-\mu)|^{-1}
\end{equation*}
from Proposition~\ref{prop:res-est-infty} and (\ref{eq:desttb}) show that this integral 
converges as a bounded operator. Note that the operator defined by this integral
agrees with the scaled resolvent follows by first varying $\theta$ while keeping $\gamma$ 
fixed, and then everywhere outside the set \eqref{eq:tens-prod-spec} by virtue of the 
analytic dependence on $\ev$.
\end{proof}

\medskip
\noindent{\bf Step 3: The parametrix with compact remainder}

\noindent We now prove the main
\begin{thm}\label{thm:ess-spec} 
The operator $\Delta_{\rad,\theta}$ has essential spectrum
\begin{equation}
\mathrm{ess\,spec}(\Delta_{\rad,\theta}) = \bigcup_{b \in I^+ \setminus \{*\}} 
\mathrm{spec}(L_{b,\theta}).
\label{eq:ess-spec}
\end{equation}
The map 
\begin{equation*}
\ev\mapsto R(\ev)=(\Delta_{\rad,\theta}-\ev)^{-1}
\end{equation*}
is meromorphic on $\Cx\setminus\cup_{b\neq *}\spec(L_{b,\theta})$ with residues of finite rank.
\end{thm}

The inclusion of the set on the right side of (\ref{eq:ess-spec}) into the set on the
left is immediate because $\Delta_{\rad,\theta}$ is well approximated by each of the
$L_{b,\theta}$ in appropriate neighbourhoods of infinity. To prove the inclusion of the
set on the left into the set on the right, it suffices to prove that when $\ev$
is outside the spectrum of $L_{b,\theta}$ for every $b \neq *$, then there is
a parametrix for the operator $(\Delta_{\rad,\theta} - \ev)^{-1}$ with compact remainder.

As before, choose a $(W,\hat\fraka)$-adapted partition of unity $\{\phi_b\}$,
$b\in I^+$, on the geodesic compactification $\hat\fraka$ of
$\fraka$, and let $\{\psi_b\}$, $b\in I^+$, be a corresponding
collection of cutoff
functions on $\hat\fraka$, so $\psi_b\in\Cinf(\hat\fraka)$ satisfies
(i)-(iii) of Definition~\ref{def:W-ahat-adapted} and
such that $\psi_b$ is identically $1$ in a neighbourhood of $\supp\phi_b$. 

Denote by $\pi:M\to\overline{\calC^+}$ and $\pi^b:\Sigma^b\to S^b_+$ the projections 
induced by the Cartan decompositions on $M$ and $\Sigma^b$. On a neighbourhood $U_b$ of 
$\supp\psi_b$, $L^2(\pi^{-1}(U_b),dg)^K$ may be identified with $L^2(\pi_b^{-1}(U_b),
e^{2(\rho-\rho_b)}\,da_b\,dg_b)^{K^b}$.  

We assume, by induction, that the spectrum of $L_{b,\theta}$ is known for every 
$b \in I^+ \setminus \{*\}$. As above, for every such $b$ let $R_{b,\theta}(\ev)=
(L_{b,\theta}-\ev)^{-1}$ for $\ev\nin\spec(L_{b,\theta})$. When $b=*$, let $R_{*,\theta}$
denote an ordinary $K$-invariant parametrix for $\Delta_{\rad,\theta}$ on some
large
ball in $\fraka$. The restriction of every $\psi_b R_{b,\theta}(\ev)\phi_b$
to $K^b$-invariant functions may be regarded as acting on $K$-invariant functions on $M$,
and with this identification we define the parametrix
\begin{equation*}
P_\theta(\ev)=\sum_{b} \psi_b R_{b,\theta}(\ev)\phi_b.
\end{equation*}

\begin{prop}
For any $k,l,r,s\in\Real$ and $\ev\nin\spec(L_{b,\theta})$, and
$\bdefh$ a defining function for $\pa\hat\fraka$, 
\begin{equation}
R_{b,\theta}(\ev): \bdefh^k \Hso^s(M)^K \longrightarrow \bdefh^k \Hso^{s+2}(M)^K,
\label{eq:R_b-supp2}
\end{equation}
is bounded; moreover, if $\chi,\phi \in\Cinf(\hat\fraka)^W$ have disjoint support, then
\begin{equation}\label{eq:R_b-supp}
\chi R_{b,\theta}(\ev)\phi:\bdefh^k \Hso^s(M)^K
\to \bdefh^l \Hso^{r}(M)^K.
\end{equation}
\end{prop}

\begin{proof}
The argument below does not depend on $\theta$ at all, so we suppress
the scaling in the already cumbersome notation. Also, assume $b \in 
I^+\setminus \{*\}$, since the result is straightforward when $b=*$. 

We first show that \eqref{eq:R_b-supp2} implies \eqref{eq:R_b-supp}. 
In fact, since the supports of $\chi$ and $\phi$ are disjoint, 
\begin{equation*}
\chi R_b(\ev)\phi=[\chi,R_b(\ev)]\phi=R_b(\ev)[L_b,\chi]R_b(\ev)\phi.
\end{equation*}
Certainly $[L_b,\chi]\in\bdefh\Diffso^1(S_b\times\Sigma^b)$ by the
Proposition~\ref{prop:a-hat-commutes},
hence is bounded as a map $\bdefh^k \Hso^{s+2}(M)^K \to \bdefh^{k+1} \Hso^{s+1}(M)^K$
due to Proposition~\ref{prop:ss-diff-op-bds}.
Using \eqref{eq:R_b-supp2}, we deduce that
\begin{equation*}
\chi R_b(\ev)\phi:\bdefh^k \Hso^s(M)^K\to \bdefh^{k+1}\Hso^{s+3}(M)^K;
\end{equation*}
iterating this proves the claim. 

Let us now prove \eqref{eq:R_b-supp2}.
The case $k=0$ follows from elliptic regularity and the definition of the
spaces $\Hso^s(M)^K$. For general $k$, we must show that 
\[
\bdefh^k R_b(\ev)\bdefh^{-k}: \Hso^s(M)^K \longrightarrow \Hso^{s+2}(M)^K.
\]
Assume that $k>0$ since the case $k<0$ then follows by applying the argument 
below to the adjoint. Using the identity
\[
[R_b(\ev),\bdefh^{-k}] = R_b(\ev)[\bdefh^{-k},L_b]R_b(\ev),
\]
we have 
\begin{equation*}
\bdefh^k R_b(\ev)\bdefh^{-k}=R_b(\ev)+\bdefh^k[R_b(\ev),\bdefh^{-k}]
=R_b(\ev)+\bdefh^k R_b(\ev)[\bdefh^{-k},L_b] R_b(\ev).
\end{equation*}
Obviously the first term on the right is bounded from $\Hso^s(M)^K$ to $\Hso^{s+2}(M)^K$.
Next, $[\bdefh^{-k},L_b]: \Hso^r(M)^K \to \Hso^{r-1}(M)^K$ is bounded provided 
$0 \leq k \leq 1$. Applying this with $r=s+2$, and using that multiplication by $\bdefh^{k}$
is bounded on $\Hso^s(M)^K$, we see that the second term on the right is
bounded from $\Hso^s(M)^K$ to $\Hso^{s+3}(M)^K$, so altogether 
$R_b(\ev): \bdefh^k \Hso^s(M)^K \to \bdefh^k \Hso^{s+2}(M)^K$ is bounded when $|k|\leq 1$.

In general, if it is known that $R_b(\ev): \bdefh^l \Hso^s(M)^K \to \bdefh^l 
\Hso^{s+2}(M)^K$ is bounded for some $l>0$, then the identity
\[
\bdefh^{k-l} R_b(\ev)\bdefh^{-k+l}= R_b(\ev)+\bdefh^{k-l}R_b(\ev)[\bdefh^{-k+l},L_b] R_b(\ev)
\]
shows that it is true for any $k$ with $l<k\leq l+1$. (This uses the boundedness of
$[\bdefh^{k-l},L_b]: \bdefh^l \Hso^{s+2}(M)^K \to \bdefh^l \Hso^{s+1}(M)^K$.)
This proves the result for all $k$. 
\end{proof}

\begin{prop}\label{prop:parametrix}
For $\ev\in\Cx\setminus \cup_{b\neq *}\spec(L_{b,\theta})$,
\begin{equation*}
P_\theta(\ev)(\Delta_{\rad,\theta}-\ev)-\Id,
(\Delta_{\rad,\theta}-\ev)P_\theta(\ev)-\Id:\bdefh^k \Hso^s(M)^K
\to \bdefh^l \Hso^{s+1}(M)^K,
\end{equation*}
for any $s,k,l\in\Real$. 
\end{prop}

\begin{proof}
Again $\theta$ plays no role, so we drop it from the notation.

For $\ev$ in the specified domain, each $R_b(\ev)$ is bounded on $L^2(M)^K$, 
by Corollary~\ref{cor:tensor-product}. Now
\begin{equation*}
(\Delta-\ev)P(\ev)=\sum_{b \in I^+} (\Delta-\ev)\psi_b R_b(\ev)\phi_b.
\end{equation*}
On $\supp \psi_b$, $b\neq *$, $\Delta=L_b+E_b$. Here
\begin{equation}\label{eq:E_b-map}
E_b\psi_b:\bdefh^k \Hso^s(M)^K
\to \bdefh^l \Hso^{s-1}(M)^K
\end{equation}
for any $k,l,s$ since $E_b\psi_b\in x^\infty\Diffso^1(M)$
by Lemma~\ref{lemma:E_b-coeffs} (and Proposition~\ref{prop:type-A} for
$\theta\nin\RR$). Hence
\begin{equation*}
(\Delta-\ev)P(\ev)= 
\sum_{b\neq *} E_b \psi_b R_b(\ev)\phi_b
+\sum_b [L_b,\psi_b] R_b(\ev)\phi_b+\sum_b \psi_b (L_b-\ev)R_b(\ev)\phi_b
\end{equation*}
By \eqref{eq:E_b-map}, the first term on the right maps 
$\bdefh^k\Hso^s(M)^K\to \bdefh^l \Hso^{s+1}(M)^K$. The third term equals
$\sum_{b} \psi_b\phi_b + Q=\Id + Q$, where $Q$ is a compactly supported pseudodifferential
operator of order $-\infty$. Finally, $[L_b,\psi_b]$ is a differential operator with 
coefficients supported in a set disjoint from $\supp\phi_b$ in $\hat\fraka$. 
The result now follows from the previous proposition.
\end{proof}

Theorem~\ref{thm:ess-spec} now follows from Proposition~\ref{prop:parametrix} and the analytic 
Fredholm theorem. 

When $\theta = 0$, there is an even stronger conclusion:

\begin{thm}
The spectrum of $\Delta_{\rad}$ is the half-line $[\,|\rho|^2,\infty)$; in other words,
there is no point spectrum below the continuous spectrum.
\end{thm}

\begin{proof}
Suppose that $\Delta_{\rad}$ has an eigenvalue $\ev_1<|\rho|^2$. Then $\ev_1$
is also an eigenvalue of $\Delta$, the Laplacian on the symmetric
space $M$. By a theorem of Sullivan \cite[Theorem~2.1]{Sullivan},
the existence of a positive 
solution to $(\Delta-\ev)u = 0$ is equivalent to $\ev \leq \inf\spec(\Delta)$,
so to prove the theorem we only need provide such a positive solution with $\ev>\ev_1$. 

To do this, recall the decomposition $G=NAK$, so that $M = G/K$ is identified
with $NA$. We consider the $N$-invariant solutions of $(\Delta - \ev)u = 0$.
The radial part of $\Delta$ with respect to the $N$-action
(i.e.\ $\Delta$ acting on $N$-invariant functions) has the form
$e^\rho \Delta_{\fraka}e^{-\rho}+|\rho|^2$, see
\cite[Chapter~II, Proposition~3.8]{Helgason:Groups}; the discrepancy in signs
arises because our Laplacian is the one with positive spectrum.
It is thus natural to consider `plane wave solutions', i.e.\ those of
the form $u(H) = \exp((\rho - \beta)(H))$, $H \in \fraka$, 
where $\beta \in \fraka^*_{\Cx}$ satisfies $-\beta\cdot\beta+|\rho|^2=\ev$.
When $\ev \in \Real$, $\ev < |\rho|^2$, then we can take $\beta \in 
\fraka$, and so $u$ is real-valued and everywhere positive.
Choosing $\ev>\ev_1$, so $\ev\in(\ev_1,|\rho|^2)$,
completes the proof as noted above.
\end{proof}

We also claim that there are no eigenvalues embedded in the continuous
spectrum, i.e.\ in the ray $(\ev_0,\infty)$. This may be proved using
$N$-body techniques, i.e.\ positive commutator techniques as in 
\cite{Vasy:Propagation-Many}. Indeed, \cite{Vasy:Exponential} proves the
corresponding result for first order $N$-body
perturbations of $\Delta_\fraka$. Unfortunately, while the method requires
only trivial modifications, the result does not apply directly due to
the apparent singularities at the Weyl chamber walls. Since setting up this
approach would require a substantial detour, we postpone this to elsewhere.

It is natural to conjecture that there are no eigenvalues in the resolvent
set of $(\Delta_{\rad,\theta}-\ev)^{-1}$ for any $\theta$ with $\im \theta < \pi/2$,
or in other words, one does not encounter poles of the
continued resolvent until one rotates a full angle of $\pi$. Furthermore,
the poles on the negative real axis should correspond
to a spectral problem on the compact dual of $M$. This can be checked directly
when $M = {\mathbb H}^n$. We expect to prove this conjecture using
purely analytic arguments, i.e.\ without resorting to representation theory.
The main point is to analyze the limiting operators
$\Delta_{\rad,\theta}$ when $\im \theta \to\pm\pi/2$; this is nontrivial since
the coefficients of this operator develop a number of new singularities
in this limit. Roughly, the limiting operators are the radial parts
of Laplacians on infinitely many copies of the compact dual, connected by 
linking `boundary conditions'. More precisely, $\im \theta \to\pm\pi/2$ is 
an analytic surgery limit, as described and studied 
in \cite{Mazzeo-Melrose:Analytic} and \cite{McDonald}: $M$ becomes pinched along the submanifolds where 
roots $\alpha$ assume values which are non-zero integer multiples of $\pi$. 
This is already seen in the expression 
\eqref{eq:eta-def} for the density $\eta\,da$.
Such a result would imply the very pleasant consequence
that the the domain of analytic continuation has only the single
ramification point $|\rho|^2$, and does not inherit the thresholds and
eigenvalues from lower rank cases as Regge poles, i.e.\ new thresholds.
Unfortunately but necessarily, the proof would be
rather involved, and it has seemed prudent
to defer it to another paper.

\bibliographystyle{plain}
\bibliography{sm}

\end{document}